\documentclass[11pt]{amsart}

\usepackage{geometry}
\usepackage{amsfonts, adjustbox, amssymb, amscd}
\numberwithin{equation}{section}

\usepackage{bm}
\usepackage{verbatim}
\usepackage{mathrsfs}
\usepackage{graphicx}
\usepackage{tikz-cd}
\usepackage{subcaption}
\usepackage{listings}
\usepackage{subfiles}
\usepackage[toc,page]{appendix}
\usepackage{mathtools}
\usepackage{comment}
\usepackage{enumerate}
\usepackage{enumitem}
\usepackage[all]{xy}

\usepackage{graphicx}
\usepackage{appendix}
\usepackage{hyperref}
\hypersetup{
    colorlinks=true,
    citecolor=red,
    linkcolor=blue,
    filecolor=magenta,      
    urlcolor=red,
}
\newcommand{\bb}{\bm{b}}
\newcommand{\Mm}{{\bf{M}}}
\newcommand{\Nn}{{\bf{N}}}

\newcommand{\Pp}{{\bf{P}}}

\newcommand{\Qq}{\mathbb{Q}}

\newcommand{\Rr}{\mathbb{R}}

\newcommand{\Exc}{\operatorname{Exc}}

\newcommand{\rk}{\operatorname{rank}}

\newcommand{\Supp}{\operatorname{Supp}}

\newcommand{\Diff}{\operatorname{Diff}}

\newcommand{\mult}{\operatorname{mult}}

\newcommand{\cont}{\operatorname{cont}}

\newcommand{\Aa}{{\bf{A}}}

\newcommand{\Bb}{{\bf{B}}}
\newcommand{\Ff}{\mathcal{F}}

\newcommand{\Oo}{\mathcal{O}}
\newcommand{\Ii}{\Gamma}

\newcommand{\Pic}{\mathrm{Pic}}


\newcounter{parentnumber}

\newtheorem{thm}{Theorem}[section]
\newtheorem{conj}[thm]{Conjecture}

\newtheorem{lem}[thm]{Lemma}
\newtheorem{prop}[thm]{Proposition}

\newtheorem{claim}[thm]{Claim}

\theoremstyle{definition}
\newtheorem{defn}[thm]{Definition}
\newtheorem{ques}[thm]{Question}
\theoremstyle{definition}
\newtheorem{rem}[thm]{Remark}

\newtheorem{deflem}[thm]{Definition-Lemma}

\newtheorem{nota}[thm]{Notation}

\theoremstyle{definition}

\begin{document}

\title{Minimal model program for algebraically integrable foliations on klt varieties}
\author{Jihao Liu, Fanjun Meng, and Lingyao Xie}

\date{\today}

\begin{abstract}
For lc algebraically integrable foliations on klt varieties, we prove the base-point-freeness theorem, the contraction theorem, and the existence of flips. The first result resolves a conjecture of Cascini and Spicer, while the latter two results strengthen a result of Cascini and Spicer by removing their assumption on the termination of flips.

Moreover, we prove the existence of the minimal model program for lc algebraically integrable foliations on klt varieties and the existence of good minimal models or Mori fiber spaces for lc algebraically integrable foliations polarized by ample divisors on klt varieties. As a consequence, we show that $\mathbb{Q}$-factorial klt varieties with lc algebraically integrable Fano foliation structures are Mori dream spaces. We also show the existence of a Shokurov-type polytope for lc algebraically integrable foliations.
\end{abstract}

\address{Department of Mathematics, Peking University, No. 5 Yiheyuan Road, Haidian District, Peking 100871, China}
\email{liujihao@math.pku.edu.cn}

\address{Department of Mathematics, Johns Hopkins University, 3400 N. Charles Street, Baltimore, MD 21218, USA}
\email{fmeng3@jhu.edu}

\address{Department of Mathematics, University of California, San Diego, 9500 Gilman Drive \# 0112, La Jolla, CA 92093-0112, USA}
\email{l6xie@ucsd.edu}

\maketitle

\pagestyle{myheadings}\markboth{\hfill Jihao Liu, Fanjun Meng, and Lingyao Xie \hfill}{\hfill Minimal model program for algebraically integrable foliations on klt varieties\hfill}

\tableofcontents

\section{Introduction}\label{sec:Introduction}
We work over the field of complex numbers $\mathbb C$.

The goal of this paper is to prove the existence of the minimal model program (MMP) for lc algebraically integrable foliations on varieties with mild singularities. One of our main theorems is the following.

\begin{thm}\label{thm: main theorem smooth projective}
Let $X$ be a smooth projective variety and $\Ff$ an algebraically integrable foliation with at worst log canonical singularities on $X$. Then we may run a $K_{\Ff}$-MMP.
\end{thm}

\noindent\textbf{History on MMP for foliations.} The minimal model program for foliations has been extensively studied in the past several years not only due to its importance on the characterization of the ambient variety and its tangent bundle, but also due to its close connection with the major conjectures of the classical minimal model program. For example, foliations play a crucial role in the proof of several key cases of the abundance conjecture for threefolds (cf. \cite{Miy87}).

To prove the existence of the minimal model program, we need at least three ingredients: the cone theorem, the contraction theorem, and the existence of flips. We can only run MMP when all of them are known. For foliations in low dimensions, these three ingredients have been established for surfaces \cite{McQ08,Bru15} and threefolds \cite{CS20,Spi20,CS21,SS22}. Although it is difficult to achieve any of these results in dimension $\geq 4$, there are some developments on the minimal model program for foliations induced by dominant rational maps recently, i.e. \emph{algebraically integrable foliations}. For example, for algebraically integrable foliations, \cite{ACSS21} proved the cone theorem in full generality, and \cite{CHLX23} proved the contraction theorem and the existence of flips when the foliations have at worst $\Qq$-factorial foliated dlt singularities. These two results together imply the existence of the minimal model program for algebraically integrable foliations with at worst $\Qq$-factorial foliated dlt singularities.

\medskip

\noindent\textbf{MMP for lc foliations on klt varieties.} It is known that we can run minimal model program for algebraically integrable foliations with at worst $\Qq$-factorial foliated dlt singularities. $\Qq$-factorial foliated dlt singularities are usually considered as a foliated version of $\Qq$-factorial dlt singularities for usual pairs \cite{CS21,CHLX23}, and foliated log smooth singularities are always $\Qq$-factorial foliated dlt.

However, Cascini and Spicer \cite{CS25a} pointed out that it is necessary to consider the minimal model program for foliations with singularities which are worse than $\Qq$-factorial foliated dlt. One major motivation is that Fano foliations (i.e. foliations with ample anti-canonical divisor $-K_{\Ff}$) are never foliated dlt (cf. \cite[Theorem 5.1]{AD13}), and they form an important topic in the theory of foliations. This makes \cite{CHLX23} not applicable to Fano foliations. 

To resolve this issue, we should consider the minimal model program for algebraically integrable foliations with at worst log canonical (lc) singularities on klt varieties, which is natural and necessary. In this paper, we prove the existence of the minimal model program under this setting.

\begin{thm}\label{thm: main}
Let $(X,\Ff,B)/U$ be an lc algebraically integrable foliated triple such that $(X,\Delta)$ is klt for some $B\geq\Delta\geq 0$. Let $R$ be a $(K_{\Ff}+B)$-negative extremal ray$/U$. Then:
    \begin{enumerate}
        \item (Contraction theorem) There exists a contraction$/U$ $\cont_R: X\rightarrow T$ of $R$.
        \item (Existence of flips) If $\cont_R$ is a flipping contraction, then the flip$/U$ $X^+\rightarrow T$ associated to $R$ exists.
    \end{enumerate}   
\end{thm}
Theorem \ref{thm: main} is known by \cite[Theorem 3.2]{CS25a} under the following additional assumptions:
\begin{itemize}
        \item The termination of $\Qq$-factorial klt flips in dimension $r=\rk\Ff$. 
        \item $B=\Delta$ with rational coefficients.
        \item $X$ is projective, and $\Qq$-factorial for statement (2).
\end{itemize}

Theorem \ref{thm: main} implies that we can run minimal model program for algebraically integrable foliations with at worst lc singularities on klt varieties. In fact, under the assumptions of Theorem \ref{thm: main}, we can show that $(X,\Delta)$ remains klt after each step of the minimal model program. Therefore, with the help of the (relative) cone theorem \cite{ACSS21,CHLX23}, we prove the following result on the existence of minimal model program:

\begin{thm}[Existence of minimal model program]\label{thm: mmp can run}
   Let $(X,\Ff,B)/U$ be an lc algebraically integrable foliated triple such that $(X,\Delta)$ is klt for some $B\geq\Delta\geq 0$. Then we may run a $(K_{\Ff}+B)$-MMP$/U$. Moreover, for any birational map $\phi: X\dashrightarrow X^+$ that is a sequence of steps of a $(K_{\Ff}+B)$-MMP$/U$, $(X^+,\Delta^+:=\phi_*\Delta)$ is klt.
\end{thm}

We remark that when $X$ is $\Qq$-factorial, which is the most natural setting when we run MMP, the condition that $(X,\Delta)$ is klt for some $B\geq\Delta\geq 0$ is equivalent to the condition that $X$ is klt. In particular, Theorem \ref{thm: main theorem smooth projective} is a direct consequence of Theorems \ref{thm: main} and \ref{thm: mmp can run}.

We also remark that Theorem \ref{thm: mmp can run} is known when $(X,\Ff,B)/U$ is $\Qq$-factorial foliated dlt by \cite[Theorem 2.1.1]{CHLX23}. Although \cite[Theorem 2.1.1]{CHLX23} does not require $X$ to be klt, $X$ is automatically klt by \cite[Theorem 2.1.9]{CHLX23}. Therefore, \cite[Theorem 2.1.1]{CHLX23} can be viewed as a special case of Theorem \ref{thm: mmp can run}.

\medskip

\noindent\textbf{Base-point-freeness theorem, good minimal models, and Mori fiber spaces.} After the establishment of the cone theorem, the contraction theorem, and the existence of flips, our next goal is to show the existence of good minimal models or Mori fiber spaces. First, we prove the existence of Mori fiber spaces for lc algebraically integrable foliations on klt varieties.

\begin{thm}[Existence of Mori fiber spaces]\label{thm: eomfs}
Let $(X,\Ff,B)/U$ be an lc algebraically integrable foliated triple such that $(X,\Delta)$ is klt for some $B\geq\Delta\geq 0$. Assume that $K_{\Ff}+B$ is not pseudo-effective$/U$.
       
Then we may run a $(K_{\Ff}+B)$-MMP$/U$ with scaling of an ample$/U$ $\Rr$-divisor and any such MMP terminates with a Mori fiber space of $(X,\Ff,B)/U$.
\end{thm}

Next, we deal with the existence of good minimal models. Unfortunately, since we do not know the existence of minimal models for smooth projective varieties in dimension $\geq 5$, we cannot prove the existence of minimal models for lc algebraically integrable foliations on klt varieties unconditionally. Nevertheless, we can prove the existence of good minimal models when the boundary divisor contains an ample $\Rr$-divisor, or when the numerical dimension is $0$.

\begin{thm}[Existence of good minimal models with polarizations]\label{thm: eolmm+A 1} 
Let $(X,\Ff,B)/U$ be an lc algebraically integrable foliated triple such that $(X,\Delta)$ is klt for some $B\geq\Delta\geq 0$. Let $A$ be an ample$/U$ $\Rr$-divisor on $X$ such that $K_{\Ff}+B+A$ is pseudo-effective$/U$. Then:
\begin{enumerate}
    \item We may run a $(K_{\Ff}+B+A)$-MMP$/U$ with scaling of an ample$/U$ $\Rr$-divisor and any such MMP terminates with a minimal model of $(X,\Ff,B+A)/U$.
    \item The minimal model in (1) is a good minimal model.
\end{enumerate}
\end{thm}

We remark that similar statements for threefold foliations in \cite{CS20,CS21,SS22} usually require that $(X,\Ff,B+A)$ is lc as Bertini-type theorems generally fail. In comparison, we do not need $(X,\Ff,B+A)$ to be lc in Theorem \ref{thm: eolmm+A 1}.

An interesting fact is that we use Theorem \ref{thm: eolmm+A 1}(1) to prove the following base-point-freeness theorem, which in return gives us Theorem \ref{thm: eolmm+A 1}(2):

\begin{thm}[Base-point-freeness theorem]\label{thm: bpf intro}
Let $(X,\Ff,B)/U$ be an lc algebraically integrable foliated triple such that $(X,\Delta)$ is klt for some $B\geq\Delta\geq 0$. Let $A$ be an ample$/U$ $\Rr$-divisor on $X$ such that $K_{\Ff}+B+A$ is nef$/U$. Then:
  \begin{enumerate}
      \item $K_{\Ff}+B+A$ is semi-ample$/U$.
      \item If $K_{\Ff}+B+A$ is Cartier, then $\mathcal{O}_X(n(K_{\Ff}+B+A))$ is globally generated$/U$ for any integer $n\gg 0$.
  \end{enumerate}
\end{thm}

In particular, Theorem \ref{thm: bpf intro} solves \cite[Conjecture 4.1]{CS25a} which further assumes that $B=\Delta$ and $(X,\Ff,B)$ is foliated dlt. We remark that \cite[Theorem A]{CHLX23} proved \cite[Conjecture 4.1]{CS25a} when $(X,\Ff,B)$ is $\Qq$-factorial foliated dlt but the non-$\Qq$-factorial case is much more difficult to prove.

An immediate consequence of Theorems \ref{thm: eolmm+A 1} and \ref{thm: bpf intro} is the finite generation of the log canonical ring for lc polarized algebraically integrable foliations on klt varieties:

\begin{thm}[Finite generation of the log canonical rings with polarizations]\label{thm: fg intro}
Let $(X,\Ff,B)/U$ be an lc algebraically integrable foliated triple such that $(X,\Delta)$ is klt for some $B\geq\Delta\geq 0$. Let $A$ be an ample$/U$ $\Rr$-divisor on $X$ such that $B+A$ is a $\Qq$-divisor. Then the log canonical ring
$$R(X,K_{\Ff}+B+A):=\oplus_{m=0}^{+\infty}\pi_*\mathcal{O}_X(\lfloor m(K_{\Ff}+B+A)\rfloor)$$
is a finitely generated $\mathcal{O}_U$-algebra.
\end{thm}

Another case when we have the existence of good minimal models is when $X$ is projective and the numerical dimension of the foliated triple is $0$.

\begin{thm}\label{thm: nt eogmm}
    Let $(X,\Ff,B)$ be a projective lc algebraically integrable foliated triple such that $(X,\Delta)$ is klt for some $B\geq\Delta\geq 0$. Assume that $\kappa_{\sigma}(K_{\Ff}+B)=0$.

    Then we may run a $(K_{\Ff}+B)$-MMP with scaling of an ample $\Rr$-divisor and any such MMP terminates with a minimal model $(X_{\min},\Ff_{\min},B_{\min})$ of $(X,\Ff,B)$ such that $K_{\Ff_{\min}}+B_{\min}\sim_{\mathbb R}0$.
\end{thm}

\medskip

\noindent\textbf{Fano foliations and Mori dream spaces.} As a consequence of Theorems \ref{thm: eomfs} and \ref{thm: eolmm+A 1}, we show that we can run MMP for any $\Rr$-Cartier $\Rr$-divisor on any klt projective variety with an lc algebraically integrable Fano foliation structure, and any such MMP terminates with either a good minimal model or a Mori fiber space. It implies that the ambient variety of an lc algebraically integrable Fano foliation is a Mori dream space if the ambient variety is $\Qq$-factorial klt.

\begin{thm}\label{thm: weak mfs}
Let $\Ff$ be an lc algebraically integrable Fano foliation on a klt projective variety $X$. Let $D$ be an  $\Rr$-Cartier $\Rr$-divisor on $X$. Then: 
\begin{enumerate}
    \item We may run a $D$-MMP which terminates with either a good minimal model of $D$ or a Mori fiber space of $D$.
    \item $X$ is a Mori dream space if it is $\Qq$-factorial.
\end{enumerate}
In particular, $D$ is semi-ample if it is nef, and the section ring $R(X,D)$ is a finitely generated $\mathbb{C}$-algebra if $D$ is a $\Qq$-divisor.

\end{thm}

Theorem \ref{thm: weak mfs} also holds in the relative setting and for the foliated Fano type case. We refer the readers to Theorem \ref{thm: weak mfs gfq} for the most general version of Theorem \ref{thm: weak mfs}.

\medskip

\noindent\textbf{Minimal models in the sense of Birkar-Shokurov.} We have studied the MMP for lc algebraically integrable foliations on klt varieties in details in Theorems \ref{thm: main}, \ref{thm: mmp can run}, \ref{thm: eomfs} and \ref{thm: eolmm+A 1}. However, it is also interesting and important to consider the MMP for lc algebraically integrable foliations when ambient varieties are not necessarily klt. They appear in birational geometry naturally. For example, for a locally stable family $f: X\rightarrow Z$, the foliation induced by $f$ is lc, but $X$ may not be klt or even lc (the singularities of $X$ can actually be very bad). Recently, MMP for locally stable families has been established in \cite{MZ23} when they study the wall crossing for moduli of stable pairs.

Unfortunately, we are unable to prove the contraction theorem or the existence of flips without any assumption on the ambient variety. Therefore, it can be difficult to talk about minimal models or Mori fiber spaces in this setting. One way to resolve this issue is to study minimal models or Mori fiber spaces \emph{in the sense of Birkar-Shokurov}, i.e. minimal models or Mori fiber spaces which allow extraction of lc places (cf. \cite{Sho96,Bir12}). We refer the readers to Definition \ref{defn: minimal model} for details.

In this paper, we prove the following two results. First, an lc foliated triple polarized by an ample divisor always has a minimal model or a Mori fiber space in the sense of Birkar-Shokurov.

\begin{thm}\label{thm: eolmm+A} 
Let $(X,\Ff,B)/U$ be an lc algebraically integrable foliated triple and $A\geq 0$ an ample$/U$ $\Rr$-divisor on $X$. Then $(X,\Ff,B+A)/U$ has either a minimal model or a Mori fiber space in the sense of Birkar-Shokurov.
\end{thm}

Second, when the ambient variety is klt, we show that the existence of a minimal model or a Mori fiber space in the sense of Birkar-Shokurov is equivalent to the termination of minimal model program with scaling of ample divisors. While the latter obviously implies the former, to prove the reverse is highly non-trivial, even for usual pairs (cf. \cite[Theorem 1.9(3)]{Bir12}).

\begin{thm}\label{thm: foliation emm equal to tof scaling}
Let $(X,\Ff,B)/U$ be an lc algebraically integrable foliated triple. Assume that $(X,\Ff,B)/U$ has a minimal model or a Mori fiber space in the sense of Birkar-Shokurov, and $X$ is potentially klt. Let $A$ be an ample$/U$ $\Rr$-divisor on $X$. Then:
\begin{enumerate}
    \item Any $(K_{\Ff}+B)$-MMP$/U$ with scaling of $A$ terminates provided that the $(K_{\Ff}+B)$-MMP$/U$ with scaling of $A$ exists.
    \item If there exists a klt pair $(X,\Delta)$ such that $B\geq\Delta\geq 0$, then $(X,\Ff,B)/U$ has a minimal model or a Mori fiber space.
\end{enumerate}
\end{thm}

\medskip

\noindent\textbf{A Shokurov-type polytope.} Finally, as an important ingredient in the proof of our main theorems, we prove the existence of a Shokurov-type polytope (cf. \cite[Corollary 1.1.5]{BCHM10}) for algebraically integrable foliations.

\begin{thm}\label{thm: shokurov polytope foliation}
    Let $(X,\Ff,B:=\sum_{i=1}^mv^0_iB_i)/Z$ be an lc algebraically integrable foliated triple such that $K_{\Ff}+B$ is nef$/Z$. Let $\bm{v}_0:=(v^0_1,\dots,v^0_m)$. Then there exists an open subset $U$ of the rational envelope of $\bm{v}_0$ in $\mathbb R^m$ such that $(X,\Ff,\sum_{i=1}^mv_iB_i)$ is lc and $K_{\Ff}+\sum_{i=1}^mv_iB_i$ is nef$/Z$ for any $(v_1,\dots,v_m)\in U$.
\end{thm}

\medskip

\noindent\textbf{Generalized foliated quadruples.} Generalized foliated quadruples play an important role in the proofs of the main theorems in many recent works (cf. \cite{LLM23,LMX24,CHLX23}). In this paper, generalized foliated quadruples are also crucially used due to the failure of Bertini-type theorems for foliations. The generalized foliated quadruple version of all our main theorems holds, although some of them require the nef part of the generalized foliated quadruple to be NQC. We refer the readers to Appendix \ref{sec: gfq case} for details. 

For the convenience of the readers, we avoid using generalized foliated quadruples in the statements and proofs of most of our results. We only essentially use this structure in Theorem \ref{thm: cont and flip with detail} and its proof, and make remarks on why we need this structure in footnotes therein. 

\medskip

\noindent\textbf{Main difficulties in the proof of the main theorems.} Roughly speaking, \cite[Theorem A]{CHLX23} established the MMP for algebraically integrable foliations that are $\Qq$-factorial foliated dlt (which implies that the ambient variety is klt). However, $\Qq$-factorial foliated dlt singularities might be too good to hope for in practice. Actually, \cite{CHLX23} showed that (as conjectured in \cite{ACSS21}) such a foliation is induced by an equidimensional morphism $f\colon X\to Z$ (not just a rational map). Moreover, (even though highly non-trivial) there exists an lc pair $(X,G)$ such that $K_\Ff\sim_{f} K_X+G$. Thanks to the cone theorem, the global $K_\Ff$-MMP turns out to be over $Z$, hence is equivalent to a $(K_X+G)$-MMP$/Z$ whose theory is well-established. Thus the authors of \cite{CHLX23} were able to deduce many corresponding results.

In general, things become much more complicated if the foliation is only induced by a rational map, in which case we do not have an associated auxiliary pair to work with. One natural idea is to consider the ``dlt" modification (whose existence is proved in \cite{ACSS21,CHLX23}) $g:X'\to X$ such that $g^{-1}\Ff$ is induced by a fibration $f':X'\to Z'$ with some desired properties. As we have explained, it is very promising that corresponding MMP on $X'$ can be run. The key point is how to descend them back to $X$. This process is subtle because the modification $g$ is not easy to control. For example, if we want to show the semi-ampleness of $K_\Ff+A$ with ample $A$, the ampleness of $A$ will not be preserved under the modification since perturbation is not allowed due to the failure of Bertini-type theorems. The best we could hope for is that the pullback of $A$ under $g$ stays ample on the leaves of $g^{-1}\Ff$ which are exactly (the reduction of) the fibers of $f':X'\to Z'$. In this case, the restrictions of $g$ on the leaves are finite morphisms, but to descend semi-ampleness under finite morphisms is quite subtle when the schemes (leaves of $g^{-1}\Ff$ and $\Ff$) are not normal (actually the properties of our leaves are even worse). Instead, it turns out that we need to apply deep and complicated MMP techniques (on $X'$) to prove the desired semi-ampleness results (e.g. base-point-free theorem) on $X$.

It is also important to notice that Fano foliations are never foliated dlt, and do not even satisfy Property $(*)$. Actually if $(X,\Ff,B)$ is a projective lc foliated triple such that $-(K_{\Ff}+B)$ is ample, then $\Ff$ is not induced by a contraction (cf. \cite[Theorem 5.1]{AD13}). Moreover, many Fano foliations are lc, and a lot of work contributed to the classification of lc Fano foliations on smooth projective varieties (cf. \cite{AD13,AD16}). Therefore, it is natural to consider the minimal model program for lc algebraically integrable Fano foliations on smooth projective varieties (e.g. Theorem \ref{thm: weak mfs}). 

\medskip
\noindent\textbf{Idea of the proof.} The key idea of the proof is the following observation. Given two contractions $h: X'\rightarrow X$ and $f: X'\rightarrow Z$ between normal quasi-projective varieties, there exists a unique normal quasi-projective variety $\bar X$ satisfying the following:
\begin{enumerate}
    \item $h$ and $f$ factor through $\bar X$.
    \item $\bar X$ is ``minimal" among all varieties satisfying $(1)$. In other words, if there exists a variety $X''$ such that $h$ and $f$ factor through $X''$, then the induced contractions $X''\rightarrow X$ and $X''\rightarrow Z$ factor through $\bar X$.
\end{enumerate}
We call $\bar X$ the \emph{core model} of $(h,f)$. It is natural and not difficult to observe the existence of such a $\bar X$. For example, given two contractions $X\rightarrow Z_0$ and $Z\rightarrow Z_0$, then the normalization of the main component of $X\times_{Z_0}Z$, $\bar X$, is automatically the core model of $(\bar X\rightarrow X,\bar X\rightarrow Z)$. In many scenarios, $\bar X\rightarrow Z$ is viewed as a ``base change" of $X\rightarrow Z_0$. For general contractions $h$ and $f$, the core model of $(h,f)$ can be viewed as an analogue of such a base change but without a base. 

Pairs of contractions $h: X'\rightarrow X$ and $f: X'\rightarrow Z$ are very common in the study of algebraically integrable foliations. Given an algebraically integrable foliation $\Ff$ on $X$, there are many cases when $h$ is a $\Qq$-factorial ACSS modification of $\Ff$ (or a $(*)$-modification of $\Ff$) and $f$ is a contraction which induces $\Ff':=h^{-1}\Ff$. It is usually easier to study $\Ff'$ due to its connection with an lc pair structure $(X',B'+G)$ (cf. \cite[Proposition 3.6]{ACSS21}). We will use $\Ff'$ to study $\Ff$.

The problem is that, when we study the minimal model program for $\Ff$ (e.g. contraction theorem, existence of flips), we usually need to consider $\Ff$ together with a polarization by an ample divisor $A$ on $X$. However, $A':=h^*A$ is only big and nef and not ample. Moreover, $(X',B'+G+A')$ is only an lc pair polarized by a big and nef divisor, and we do not know the existence of good minimal models for such pairs yet. It causes troubles for us to use the minimal model program for $(X',B'+G+A')$ to study the minimal model program for $K_{\Ff}+A$. 

Nevertheless, by using the core models we have introduced, we can resolve this problem. Let $\bar X$ be the core model of $(h,f)$ with $\bar h: \bar X\rightarrow X$, $\bar f: \bar X\rightarrow Z$, and $g: X'\rightarrow\bar X$. Then:
\begin{itemize}
    \item Although $\bar A:=\bar h^*A$ is no longer ample, it is ample$/Z$.
    \item Let $\bar B:=g_*B'$ and $\bar G:=g_*G$. Since $g$ is a contraction over $Z$, $(\bar X,\bar B+\bar G)$ is crepant to $(X',B'+G)$ by the negativity lemma and thus it is lc.
\end{itemize}
In particular, $\bar X$ not only has an lc pair structure associated to the foliation $\bar\Ff:=\bar h^{-1}\Ff$ (while $X$ does not necessarily have) but also preserves some information of $X$ via the divisor $\bar A$, at least over $Z$. The core model $\bar X$, instead of $X'$, seems to be a more natural object for us to study due to its uniqueness with respect to the base. Moreover, it preserves more core information of $X$ comparing to an arbitrary model $X'$. This is why we call it as a ``core model".  This is also inspired by \cite[Proof of Theorem 1.5]{MZ23}. We study the basic properties of core models and its relationship with foliations in Section \ref{sec: core model}.

In our case, the MMP$/Z$ for $(\bar X,\bar B+\bar G+\bar A)$ with scaling of ample divisors behaves nicely. Moreover, by our construction, $(\bar X,\bar B+\bar G)$ has a close connection with the induced foliation $\bar\Ff$ on $\bar X$, and we can show that the MMP$/Z$ for $K_{\bar\Ff}+\bar B+\bar A$ with scaling of ample divisors behaves nicely. By the cone theorem for algebraically integrable foliations (\cite[Theorem 3.9]{ACSS21}), any $(K_{\bar\Ff}+\bar B+\bar A)$-MMP is a $(K_{\bar\Ff}+\bar B+\bar A)$-MMP$/Z$. Thus any $(K_{\bar\Ff}+\bar B+\bar A)$-MMP with scaling of ample divisors behaves nicely. Then we can use this fact to study the minimal model program for $K_{\Ff}+A$ with scaling of ample divisors, and hence for $K_{\Ff}$ by adopting the some ideas from \cite{CS25a,CHLX23}.

More precisely, the key idea of \cite[Proof of Theorem 3.2]{CS25a} is the following:
\begin{itemize}
    \item We take a $(*)$-model"\footnote{They are called ``Property $(*)$-models" in \cite{ACSS21} and the arXiv version of \cite{CS25a}. Cascini suggested us that the name ``$(*)$-models" is better.} of $(X,\Ff,B)$ which satisfies good properties. Achieving this requires us to run the ``first MMP" for a foliated log smooth model.
    \item We run the ``second MMP" which contracts the strict transforms of the $(K_{\Ff}+B)$-negative extremal rays. Termination of flips is needed here.
    \item We run the ``third MMP" which contracts the strict transform of the exceptional divisor of the  $(*)$-modification. Termination of flips is again needed here.
\end{itemize}
We follow the same idea to prove Theorems \ref{thm: main} and \ref{thm: mmp can run}, but many arguments are different. We do not need to run the ``first MMP" and we shall directly use our core model $\bar X$. We use properties of the core model $\bar X$ to show that the ``second MMP" could terminate. Finally, we can use the termination of MMP with scaling for klt generalized pairs of log general type to get the termination of the ``third MMP".

There are several additional things to remark for our proof of Theorems \ref{thm: main} and \ref{thm: mmp can run}.
\begin{enumerate}
    \item We need the minimal model program for generalized foliated quadruples because Bertini-type theorems generally fail for algebraically integrable foliations. We need the concept of generalized foliated quadruples to consider MMP with scaling of ample divisors in more details. We need results on the minimal model program for generalized foliated quadruples in \cite{CHLX23}.
    \item When the boundary $B$ has irrational coefficients, we need to establish the existence of a Shokurov-type polytope (Theorem \ref{thm: shokurov polytope foliation}) for algebraically integrable foliations in order to show that the ``third MMP" terminates. This is done in Section \ref{sec: sho polytope}. Moreover, if we consider the contraction theorem and the existence of flips for (non-NQC) generalized foliated quadruples instead of foliated triples, the Shokurov-type polytope does not exist. We need to resolve this issue by introducing and studying a special class of nef $\Rr$-divisors, namely ``$\epsilon$-nef $\Rr$-divisors". See Appendix \ref{sec: enef} for details.
    \item When $X$ is not $\Qq$-factorial, it can be tricky to run MMP on $\bar X$ since $\bar X$ may not be $\Qq$-factorial either. In this case, we need to prove some results on the minimal model program on $X'$. Nevertheless, we can still use $\bar X$ as an auxiliary variety to help us establish the MMP on $X'$, hence the ``second MMP" on $X'$. See Section \ref{sec: eolmm} for details. Also, when $X$ is not $\Qq$-factorial, the arguments in \cite{CS25a} no longer work for the existence of flips, and we need an alternative argument. Our choice is to consider the ample model of a special $\Rr$-divisor over the base of the contraction. In this case, some basic properties on different types of models of foliations need to be proved. See Section \ref{sec: models} for details.
\end{enumerate}

Finally, we say a few more words about the proof of other main theorems. First, by using core models, we establish an MMP on $\bar X$ with scaling of the pullback of an ample divisor on $X$ in Section \ref{sec: eolmm}, and show that such a MMP terminates in some cases (Theorem \ref{thm: take strict simple model run mmp}). Theorem \ref{thm: eolmm+A} follows from the establishment of such a MMP. Theorem \ref{thm: nt eogmm} follows from the establishment of such a MMP, the fact that movable divisors with zero numerical Iitaka dimension are numerically trivial, and the abundance for numerically trivial algebraically integrable foliations.

To prove Theorems \ref{thm: eomfs} and \ref{thm: eolmm+A 1}, we need to lift MMP to $\Qq$-factorial ACSS models. \cite[Remark 3.3]{CS25a} briefly mentions such lifting for flips. We discuss the lifting of the MMP in more details in Section \ref{sec: lifting}, which allows divisorial contractions and non-$\Qq$-factorial MMP to be lifted. More importantly, we show that MMP with scaling can also be lifted. This, together with the results in Section \ref{sec: eolmm}, implies Theorems \ref{thm: eomfs} and \ref{thm: eolmm+A 1}. Theorem \ref{thm: weak mfs} is a direct consequence of Theorems \ref{thm: eomfs} and \ref{thm: eolmm+A 1}, since any $D$-MMP is a $(K_{\Ff}+(-K_{\Ff}+\epsilon D))$-MMP for some $0<\epsilon\ll 1$ such that $-K_{\Ff}+\epsilon D$ is ample. We remark that we need generalized foliated quadruples again to prove these theorems due to the failure of Bertini-type theorems.

Finally, to prove Theorem \ref{thm: foliation emm equal to tof scaling}, we need to show that the existence of minimal models in the sense of Birkar-Shokurov for $(X,\Ff,B)$ is equivalent to the existence of minimal models of a pair $(X',B'+G)$ which is related to $(X,\Ff,B)$. Therefore, we can use the known results on the existence of minimal models for $(X',B'+G)$ to deduce Theorem \ref{thm: foliation emm equal to tof scaling}. Such a relationship is automatic when we have an equidimensional Property $(*)$ structure, but it is more complicated when $(X,\Ff,B)$ does not satisfy Property $(*)$. The task is done in Section \ref{sec: models}. A key observation is to reinterpret an MMP$/U$ which is also an MMP$/Z$ as an MMP$/Z_U$, where $Z_U$ is the core model of $(X\rightarrow U, X\rightarrow Z)$. The use of the auxiliary variety $Z_U$ will greatly help us transform the $(K_{\Ff}+B)$-MMP into the $(K_{X'}+B'+G)$-MMP and lead to the proof of Theorem \ref{thm: foliation emm equal to tof scaling}.

\medskip
\noindent\textit{Sketch of the paper.} In Section \ref{sec: preliminaries} we recall some preliminary results on algebraically integrable foliations and results in \cite{CHLX23}. In Section \ref{sec: core model} we introduce the concept of core models and study its basic properties. In Section \ref{sec: models} we study different types of models for algebraically integrable foliations. In Section \ref{sec: eolmm} we use the concept of core models to study the MMP for $\Qq$-factorial ACSS foliated triples polarized by the pullback of an ample divisor and prove Theorem \ref{thm: eolmm+A}. In Section \ref{sec: sho polytope} we construct a Shokurov-type polytope for algebraically integrable foliations and prove Theorem \ref{thm: shokurov polytope foliation}. In Section \ref{sec: proof flip} we prove Theorems \ref{thm: main} and \ref{thm: mmp can run}.  In Section \ref{sec: lifting} we study the lifting of the minimal model program for algebraically integrable foliations to $\Qq$-factorial ACSS models. In Section \ref{sec: proof of main theorem} we prove the rest of our main theorems for foliated triples. In Section \ref{sec: problem} we propose and discuss some remaining open problems on the minimal model program for algebraically integrable foliations and prove some results that might be useful for future applications. Finally, Appendices \ref{sec: gfq case} and \ref{sec: enef} focus on generalized foliated quadruple. In Appendix \ref{sec: gfq case}, we state and prove the generalized foliated quadruple version of our main theorems. In Appendix \ref{sec: enef}, we introduce the concept of $\epsilon$-nefness, which is a replacement of the Shokurov-type polytope for generalized foliated quadruples.

\medskip

\noindent\textbf{Acknowledgements.} The work is supported by the National Key R\&D Program of China \#2024YFA1014400. The second author is partially supported by an AMS-Simons Travel Grant. The third author was partially supported by NSF research grants no. DMS-1801851 and DMS-1952522, as well as a grant from the Simons Foundation (Award Number: 256202). The third author is supported by another grant from the Simons Foundation. The authors would like to thank Caucher Birkar, Paolo Cascini, Guodu Chen, Christopher D. Hacon, Jingjun Han, Junpeng Jiao, Yuchen Liu, Vyacheslav V. Shokurov, Calum Spicer, Roberto Svaldi, Chenyang Xu, Qingyuan Xue, Ziwen Zhu and Ziquan Zhuang for useful discussions. We would also like to thank the referees for their work and helpful comments.

\section{Preliminaries}\label{sec: preliminaries} 

We will adopt the standard notations and definitions on MMP in \cite{KM98,BCHM10} and use them freely. For foliations and foliated triples, we adopt the notations and definitions in \cite{CHLX23} which generally align with \cite{CS20, ACSS21, CS21} with possible minor differences. 

\subsection{Special notation}

\begin{defn}
    A \emph{contraction} is a projective morphism between varieties such that $f_*\mathcal{O}_X=\mathcal{O}_Y$.
\end{defn}

\begin{nota}
    Let $f: X\dashrightarrow X'$ be a birational map between normal varieties. We denote by $\Exc(f)$ the reduced divisor supported on the codimension one part of the exceptional locus of $f$.
\end{nota}

\begin{nota}
    Let $X$ be a normal variety and $D,D'$ two $\Rr$-divisors on $X$. We define 
    $D\wedge D':=\sum_P\min\{\mult_PD,\mult_PD'\}P$ where the sum runs through all the prime divisors $P$ on $X$. We denote by $\Supp D$ the reduced divisor supported on $D$.
\end{nota}

\begin{defn}
Let $m$ be a positive integer and $\bm{v}\in\mathbb R^m$. The \emph{rational envelope} of $\bm{v}$ is the minimal rational affine subspace of $\mathbb R^m$ which contains $\bm{v}$. For example, if $m=2$ and $\bm{v}=\left(\frac{\sqrt{2}}{2},1-\frac{\sqrt{2}}{2}\right)$, then the rational envelope of $\bm{v}$ is $(x_1+x_2=1)\subset\mathbb R^2_{x_1x_2}$.
\end{defn}

\begin{nota}
    A \emph{general choice} of a real number $a$ is a choice of a real number such that $a\not\in\mathbb Q(\Ii_0)$ for a finite set $\Ii_0\subset\mathbb R$. Here $\mathbb Q(\Ii_0)$ is the field extension of $\mathbb Q$ by elements in $\Ii_0$. We also say that $a$ is \emph{general} in $\mathbb R/\mathbb Q$.
\end{nota}

\subsection{Foliations}

\begin{defn}[Foliations, {cf. \cite{ACSS21,CS21}}]\label{defn: foliation}
Let $X$ be a normal variety. A \emph{foliation} on $X$ is a coherent sheaf $\Ff\subset T_X$ such that
\begin{enumerate}
    \item $\Ff$ is saturated in $T_X$, i.e. $T_X/\Ff$ is torsion free, and
    \item $\Ff$ is closed under the Lie bracket.
\end{enumerate}
The \emph{rank} of the foliation $\Ff$ is the rank of $\Ff$ as a sheaf and is denoted by $\rk\Ff$. The \emph{co-rank} of $\Ff$ is $\dim X-\rk\Ff$. The \emph{canonical divisor} of $\Ff$ is a divisor $K_\Ff$ such that $\mathcal{O}_X(-K_{\mathcal{F}})\cong\mathrm{det}(\Ff)$. If $\Ff=0$, then we say that $\Ff$ is a \emph{foliation by points}.

Given any dominant map $h: Y\dashrightarrow X$, we denote by $h^{-1}\Ff$ the \emph{pullback} of $\Ff$ on $Y$ as constructed in \cite[3.2]{Dru21} and say that $h^{-1}\Ff$ is \emph{induced by} $\Ff$. Given any birational map $g: X\dashrightarrow X'$, we denote by $g_*\Ff:=(g^{-1})^{-1}\Ff$ the \emph{pushforward} of $\Ff$ on $X'$ and also say that $g_*\Ff$ is \emph{induced by} $\Ff$. We say that $\Ff$ is an \emph{algebraically integrable foliation} if there exists a dominant map $f: X\dashrightarrow Z$ such that $\Ff=f^{-1}\Ff_Z$, where $\Ff_Z$ is the foliation by points on $Z$, and we say that $\Ff$ is \emph{induced by} $f$.

A subvariety $S\subset X$ is called \emph{$\Ff$-invariant} if for any open subset $U\subset X$ and any section $\partial\in H^0(U,\Ff)$, we have $\partial(\mathcal{I}_{S\cap U})\subset \mathcal{I}_{S\cap U}$, 
where $\mathcal{I}_{S\cap U}$ is the ideal sheaf of $S\cap U$.  For any prime divisor $P$ on $X$, we denote $\epsilon_{\Ff}(P):=1$ if $P$ is not $\Ff$-invariant and $\epsilon_{\Ff}(P):=0$ if $P$ is $\Ff$-invariant. For any prime divisor $E$ over $X$, we define $\epsilon_{\Ff}(E):=\epsilon_{\Ff_Y}(E)$ where $h: Y\dashrightarrow X$ is a birational map such that $E$ is on $Y$ and $\Ff_Y:=h^{-1}\Ff$. For any $\Rr$-divisor $D$ on $X$, we denote by $D^{\Ff-\mathrm{ninv}}$ the reduced divisor supported on the union of non-$\Ff$-invariant components of $D$.
\end{defn}

\begin{defn}[Tangent, {cf. \cite[Section 3.4]{ACSS21}}]\label{defn: tangent to foliation}
 Let $X$ be a normal variety, $\Ff$ a foliation on $X$, and $V\subset X$ a subvariety. Suppose that $\Ff$ is a foliation induced by a dominant rational map $X\dashrightarrow Z$. We say that $V$ is \emph{tangent} to $\Ff$ if there exists a birational morphism $\mu: X'\rightarrow X$, an equidimensional contraction $f': X'\rightarrow Z$, and a subvariety $V'\subset X'$, such that
    \begin{enumerate}
    \item $\mu^{-1}\Ff$ is induced by $f'$, and
        \item $V'$ is contained in a fiber of $f'$ and $\mu(V')=V$.
    \end{enumerate}
\end{defn}

\begin{defn}[Foliated triples]\label{defn: foliated triple}
 A \emph{foliated triple} $(X,\Ff,B)/U$ consists of a normal quasi-projective variety $X$, a foliation $\Ff$ on $X$, an $\Rr$-divisor $B\geq 0$ on $X$, and a projective morphism $X\rightarrow U$, such that $K_{\Ff}+B$ is $\mathbb R$-Cartier.

  If $\Ff=T_X$, then we may drop $\Ff$ and say that $(X,B)/U$ is a \emph{pair}. If $U$ is not important, then we may drop $U$. If $\Ff$ is algebraically integrable, then we say that $(X,\Ff,B)$ is algebraically integrable. If $X$ is $\Qq$-factorial, then we say that $(X,\Ff,B)$ is $\Qq$-factorial. If we allow $B$ to have negative coefficients, then we shall add the prefix ``sub-". If $B$ is a $\Qq$-divisor then we may add the prefix ``$\mathbb Q$-".
\end{defn}

\begin{defn}[Singularities]\label{defn: foliation singularity}
Let $(X,\Ff,B)$ be a foliated triple. For any prime divisor $E$ over $X$, let $f: Y\rightarrow X$ be a birational morphism such that $E$ is on $Y$, and suppose that
$$K_{\Ff_Y}+B_Y:=f^*(K_\Ff+B)$$
where $\Ff_Y:=f^{-1}\Ff$. We define $a(E,\Ff,B):=-\mult_EB_Y$ to be the \emph{discrepancy} of $E$ with respect to $(X,\Ff,B)$. If $\Ff=T_X$, then we define $a(E,X,B):=a(E,\Ff,B)$ which is the usual discrepancy for pairs.

We say that $(X,\Ff,B)$ is \emph{lc} (resp. \emph{klt}) if $a(E,\Ff,B)\geq -1$ (resp. $>-1$) for any prime divisor $E$ over $X$. For foliated sub-triples, we define singularities in the same way and we shall add the prefix ``sub-" for the descriptions of singularities.

An \emph{lc place} of $(X,\Ff,B)$ is a prime divisor $E$ over $X$ such that $a(E,\Ff,B)=-\epsilon_{\Ff}(E)$. An \emph{lc center} of $(X,\Ff,B)$ is the center of an lc place of $(X,\Ff,B)$ on $X$.

We remark that our definition of lc and klt singularities has some differences compared with the classical definition \cite{McQ08,CS20,ACSS21,CS21,CHLX23}, where the $-1$ in the inequality is replaced with $-\epsilon_{\Ff}(E)$. The next lemma shows that the two definitions on lc coincide so we are free to use results on ``lc foliations" in literature. Moreover, there are good reasons why we refine the definition of klt singularities. We refer the readers to Remark \ref{rem: new definition klt} for details.
\end{defn}

\begin{lem}\label{lem: equivalence definition lc}
    Let $(X,\Ff,B)$ be a foliated sub-triple. The following two conditions are equivalent:
    \begin{enumerate}
        \item $(X,\Ff,B)$ is sub-lc.
        \item $a(E,\Ff,B)\geq -\epsilon_{\Ff}(E)$ for any prime divisor $E$ over $X$.
    \end{enumerate}
\end{lem}
\begin{proof}
    It is clear that (2) implies (1) so we only need to show that (1) implies (2). Suppose the lemma does not hold. Then there exists a prime divisor $E$ over $X$ such that $E$ is $\Ff$-invariant and $a(E,\Ff,B)<0$. Possibly replacing $X$ by a high model, we may assume that $E$ is on $X$ and $X$ is smooth. Thus $E$ is a component of $B$ and $\mult_EB>0$. This contradicts \cite[Remark 2.3]{CS21}.
\end{proof}

\begin{defn}[Potentially klt]\label{defn: potentially klt}
Let $X$ be a normal quasi-projective variety. We say that $X$ is \emph{potentially klt} if $(X,\Delta)$ is klt for some $\Rr$-divisor $\Delta\geq 0$. 
\end{defn}

\begin{lem}\label{lem: gklt is klt}
    Let $(X,B)/U$ be an lc pair such that $X$ is potentially klt and $A$ an ample$/U$ $\Rr$-divisor. Then there exists a klt pair $(X,\Delta)$ such that $\Delta\sim_{\mathbb R,U}B+A$.
\end{lem}
\begin{proof}
    There exist a klt pair $(X,\Delta_0)$ and a real number $\epsilon>0$ sufficiently small such that $H_0:=A+\epsilon(B-\Delta_0)$ is ample$/U$. Let $H$ be a general member in $|H_0|_{\mathbb R/U}$. Then $\Delta:=(1-\epsilon)B+\epsilon\Delta_0+H$ satisfies our requirements. 
\end{proof}

\subsection{Special algebraically integrable foliations}

\begin{defn}[Foliated log resolutions]\label{defn: log resolution}
We refer the readers to \cite[Definition 6.2.1]{CHLX23} or \cite[3.2. Log canonical foliated pairs]{ACSS21} for the definition of being foliated log smooth. 

Let $X$ be a normal quasi-projective variety, $B$ an $\Rr$-divisor on $X$, and $\Ff$ an algebraically integrable foliation on $X$. A \emph{foliated log resolution} of $(X,\Ff,B)$ is a birational morphism $h: X'\rightarrow X$ such that 
$$(X',\Ff':=h^{-1}\Ff,B':=h^{-1}_*B+\Exc(h))$$ 
is foliated log smooth. The existence of foliated log resolutions for any such $(X,\Ff,B)$ is guaranteed by \cite[Lemma 6.2.4]{CHLX23}.
\end{defn}

\begin{defn}[Property $(*)$ foliations, {\cite[Definition 3.8]{ACSS21}, \cite[Definition 7.2.2]{CHLX23}}]\label{defn: foliation property *}
Let $(X,\Ff,B)$ be a foliated triple. Let $G\geq 0$ be a reduced divisor on $X$ and $f: X\rightarrow Z$ a contraction. We say that $(X,\Ff,B;G)/Z$ satisfies \emph{Property $(*)$} if the following conditions hold.
\begin{enumerate}
\item $\Ff$ is induced by $f$ and $G$ is an $\Ff$-invariant divisor.
\item $f(G)$ is of pure codimension $1$, $(Z,f(G))$ is log smooth, and $G=f^{-1}(f(G))$.
\item For any closed point $z\in Z$ and any reduced divisor  $\Sigma\ge f(G)$ on $Z$ such that  $(Z,\Sigma)$ is log smooth near $z$, $(X,B+G+f^*(\Sigma-f(G)))$ is lc over a neighborhood of $z$.
\end{enumerate}
We say that $f$, $Z$, and $G$ are \emph{associated} with $(X,\Ff,B)$.
\end{defn}

\begin{prop}[cf. {\cite[Proposition 3.6]{ACSS21}, \cite[Proposition 7.3.6]{CHLX23}}]\label{prop: weak cbf gfq}
Let $(X,\Ff,B)$ be a foliated triple. Let  $G\geq 0$ be a reduced divisor on $X$ and $f: X\rightarrow Z$ an equidimensional contraction, such that $(X,\Ff,B;G)/Z$ satisfies Property $(*)$ and $B$ is horizontal$/Z$. Then
$$K_{\Ff}+B\sim_{Z}K_X+B+G.$$
\end{prop}

\begin{defn}
Let $f: X\rightarrow Z$ be a projective morphism between normal quasi-projective varieties and $G\geq 0$ an $\Rr$-divisor on $X$. We say that $G$ is \emph{super$/Z$} if there exist ample Cartier divisors $H_1,\dots,H_{m}$ on $Z$ such that $G\geq\sum_{i=1}^{m}f^*H_i$, where $m:=2\dim X+1$.
\end{defn}

\begin{defn}[ACSS, {cf. \cite[Definitions 5.4.2, 7.2.2, 7.2.3]{CHLX23}}]\label{defn: ACSS f-triple}
Let $(X,\Ff,B)$ be an lc foliated triple, $G\geq 0$ a reduced divisor on $X$, and $f: X\rightarrow Z$ a contraction. We say that $(X,\Ff,B;G)/Z$ is \emph{ACSS} if the following conditions hold:
\begin{enumerate}    
\item $(X,\Ff,B;G)/Z$ satisfies Property $(*)$.
\item $f$ is equidimensional.
\item There exists an $\Rr$-Cartier $\Rr$-divisor $D\geq 0$ on $X$, such that  $\Supp\{B\}\subset\Supp D$, and for any reduced divisor $\Sigma\geq f(G)$ such that $(Z,\Sigma)$ is log smooth, $$(X,B+D+G+f^*(\Sigma-f(G)))$$ 
      is qdlt (cf. \cite[Definition 35]{dFKX17}).
\item For any lc center of $(X,\Ff,B)$ with generic point $\eta$, over a neighborhood of $\eta$,
    \begin{enumerate}
      \item $\eta$ is the generic point of an lc center of $(X,\Ff,\lfloor B\rfloor)$, and
       \item $f: (X,B+G)\rightarrow (Z,f(G))$ is a toroidal morphism,
    \end{enumerate}
\end{enumerate}
If  $(X,\Ff,B;G)/Z$ is ACSS and $G$ is super$/Z$, then we say that $(X,\Ff,B;G)/Z$ is super ACSS. If $(X,\Ff,B;G)/Z$ is (super) ACSS, then we say that $(X,\Ff,B)/Z$ and $(X,\Ff,B)$ are (super) ACSS.
\end{defn}

\subsection{Birational maps in MMP}

\begin{defn}
    Let $X\rightarrow U$ be a projective morphism from a normal quasi-projective variety to a variety.  Let $D$ be an $\Rr$-Cartier $\Rr$-divisor on $X$ and $\phi: X\dashrightarrow X'$ a birational map$/U$. Then we say that $X'$ is a \emph{birational model} of $X$. We say that $\phi$ is $D$-non-positive (resp. $D$-negative, $D$-trivial, $D$-non-negative, $D$-positive) if the following conditions hold:
    \begin{enumerate}
    \item $\phi$ does not extract any divisor.
    \item $D':=\phi_*D$ is $\Rr$-Cartier.
    \item There exists a resolution of indeterminacy $p: W\rightarrow X$ and $q: W\rightarrow X'$, such that
    $$p^*D=q^*D'+F$$
    where $F\geq 0$ (resp. $F\geq 0$ and $\Supp p_*F=\Exc(\phi)$, $F=0$, $0\geq F$, $0\geq F$ and $\Supp p_*F=\Exc(\phi)$).
    \end{enumerate}
\end{defn}

\begin{defn}
    Let $X\rightarrow U$ be a projective morphism from a normal quasi-projective variety to a variety. Let $D$ be an $\Rr$-Cartier $\Rr$-divisor on $X$ and $f: X\rightarrow Z$ a contraction$/U$.  We say that $f$ is a \emph{$D$-Mori fiber space$/U$} if $f$ is a contraction of a $D$-negative extremal ray$/U$ and $\dim X>\dim Z$. If $f: X\rightarrow Z$ is a \emph{$D$-Mori fiber space$/U$} for some $\Rr$-Cartier $\Rr$-divisor $D$, then we say that $f: X\rightarrow Z$ is a \emph{Mori fiber space$/U$}. If $f$ is obviously a contraction$/U$ or the ``$/U$" property is not important, then we may drop the ``$/U$" in the definitions.
\end{defn}

\begin{defn} Let $X\rightarrow U$ be a projective morphism from a normal quasi-projective variety to a variety. Let $D$ be an $\Rr$-Cartier $\Rr$-divisor on $X$, $\phi: X\dashrightarrow X'$ a $D$-negative map$/U$ and $D':=\phi_*D$.
\begin{enumerate}
\item We say that $X'$ is a \emph{minimal model}$/U$ of $D$ if $D'$ is nef$/U$.
\item We say that $X'$ is a \emph{good minimal model}$/U$ of $D$ if $D'$ is semi-ample$/U$.
\item A contraction$/U$ $f: X'\rightarrow Z$ is called a \emph{Mori fiber space}$/U$ of $D$ if $f$ is a $D'$-Mori fiber space$/U$.
\end{enumerate}
\end{defn}

\begin{lem}\label{lem: two supporting function trivial}
Let $X\rightarrow U$ be a projective morphism from a normal quasi-projective variety to a variety and $F$ an extremal face in $\overline{NE}(X/U)$. Let $H_1,H_2$ be two supporting functions$/U$ of $F$ and $\phi: X\dashrightarrow X'$ an $H_1$-trivial birational map$/U$. Then $\phi$ is $H_2$-trivial.
\end{lem}
\begin{proof}
    Let $p: W\rightarrow X$ and $q: W\rightarrow X'$ be a resolution of indeterminacy. Then there exists a unique extremal face $F_W$ of $\overline{NE}(W/U)$ such that $p_*F_W=F$, $\overline{NE}(W/X)\subset F_W$, and $p^*H_1,p^*H_2$ are both supporting functions of $F_W$. 

    Since $\phi$ is $H_1$-trivial, $q$ is $p^*H_1$-trivial. Therefore, $q$ only contracts $p^*H_1$-trivial extremal rays in $\overline{NE}(W/U)$, so $q$ only contracts $p^*H_2$-trivial extremal rays in $\overline{NE}(W/U)$. Thus $q$ is $p^*H_2$-trivial, so $\phi$ is $H_2$-trivial.
\end{proof}

\begin{lem}\label{lem: general real coefficient b-trivial}
Let $X\rightarrow U$ be a projective morphism from a normal quasi-projective variety to a variety. Let $A,B$ be two $\Rr$-divisors on $X$ and let $t$ be a real number such that $t$ is general in $\mathbb R/\mathbb Q$ and $A+tB$ is $\Rr$-Cartier. Then $A,B$ are $\Rr$-Cartier, and any $(A+tB)$-trivial map $\phi: X\dashrightarrow X'$ is $A$-trivial and $B$-trivial.
\end{lem}
\begin{proof}
Let $A'$ and $B'$ be the images of $A,B$ on $X'$ respectively. Then $A'+tB'$ is $\Rr$-Cartier.     By \cite[Lemma 5.3]{HLS24},  $A,B,A',B'$ are $\Rr$-Cartier. Let $p: W\rightarrow X$ and $q: W\rightarrow X'$ be a resolution of indeterminacy, then
$$p^*A+tp^*B=p^*(A+tB)=q^*(A'+tB')=q^*A'+tq^*B'.$$
 By \cite[Lemma 5.3]{HLS24}, $p^*A=q^*A'$ and $p^*B=q^*B'$. The lemma follows.
\end{proof}

\subsection{Relative Nakayama-Zariski decompositions}

\begin{defn}
    Let $\pi: X\rightarrow U$ be a projective morphism from a normal variety to a variety, $D$ a pseudo-effective$/U$ $\Rr$-Cartier $\Rr$-divisor on $X$, and $P$ a prime divisor on $X$. We define $\sigma_{P}(X/U,D)$ as in \cite[Definition 3.1]{LX23a} by considering $\sigma_{P}(X/U,D)$ as a number in  $[0,+\infty)\cup\{+\infty\}$. We define $N_{\sigma}(X/U,D)=\sum_Q\sigma_Q(X/U,D)Q$
    where the sum runs through all prime divisors on $X$ and consider it as a formal sum of divisors with coefficients in $[0,+\infty)\cup\{+\infty\}$.
\end{defn}

\begin{lem}[{\cite[Lemma 3.4(2)(3), Lemma 3.7(4)]{LX23a}}]\label{lem: nz keep under pullback}
Let $\pi: X\rightarrow U$ be a projective morphism from a normal variety to a variety and $D$ a pseudo-effective$/U$ $\Rr$-Cartier $\Rr$-divisor on $X$. Let $f: Y\rightarrow X$ be a projective birational morphism. Then:
\begin{enumerate}
    \item For any exceptional$/X$ $\Rr$-Cartier $\Rr$-divisor $E\ge0$ and any prime divisor $P$ on $Y$, we have $$\sigma_P(Y/U,f^*D+E)=\sigma_P(Y/U,f^*D)+\mult_PE.$$
    \item For any exceptional$/X$ $\Rr$-Cartier $\Rr$-divisor $E\ge0$ on $Y$, we have 
    $$N_{\sigma}(X/U,D)=f_*N_{\sigma}(Y/U,f^*D+E).$$
    \item $\Supp N_{\sigma}(X/U,D)$ coincides with the divisorial part of $\Bb_-(X/U,D)$.
\end{enumerate}
\end{lem}

\begin{lem}\label{lem: nz for lc divisor}
Let $X\rightarrow U$ be a projective morphism from a normal variety to a variety and $\phi: X\dashrightarrow X'$ a birational map$/U$. Let $D$ be an $\Rr$-Cartier $\Rr$-divisor on $X$ such that $\phi$ is $D$-negative and $D':=\phi_*D$. Then:
\begin{enumerate}
    \item The divisors contracted by $\phi$ are contained in $\Supp N_{\sigma}(X/U,D)$.
    \item If $D'$ is movable$/U$, then $\Supp N_{\sigma}(X/U,D)$ is the set of all $\phi$-exceptional divisors.
    \end{enumerate}
 \end{lem}
\begin{proof}
Let $p: W\rightarrow X$ and $q: W\rightarrow X'$ be a resolution of indeterminacy. Then
$$p^*D=q^*D'+E$$
for some $E\geq 0$ that is exceptional$/X'$ and $\Supp E$ contains the strict transforms of all $\phi$-exceptional divisors on $W$. By Lemma \ref{lem: nz keep under pullback}(1),
$$
\Supp E\subset\Supp N_\sigma(W/U,q^*D'+E)=\Supp N_\sigma(W/U,p^*D),
$$
By Lemma \ref{lem: nz keep under pullback}(2), $\Supp p_*E\subset\Supp N_\sigma(X/U,D)$. Therefore, any $\phi$-exceptional divisor is contained in $\Supp N_\sigma(X/U,D)$.

If $D'$ is movable$/U$, then by Lemma \ref{lem: nz keep under pullback}(2), $q_*N_\sigma(W/U,q^*D'+E)=0$. Thus $$\Supp N_\sigma(W/U,q^*D'+E)$$ is $q$-exceptional. By Lemma \ref{lem: nz keep under pullback}(2) again, we have $\Supp N_\sigma(X/U,D)=\Supp p_*N_\sigma(W/U,q^*D'+E)$, whose components are all $\phi$-exceptional. (2) follows from (1).
\end{proof}

\subsection{Generalized pairs and generalized foliated quadruples}

\begin{rem}
Generalized pairs (\cite[Definition 1.4]{BZ16}) and generalized foliated quadruples (\cite[Definition 1.2]{LLM23},\cite[Definition 3.4.3]{CHLX23}) will be inevitably used in this paper. They are crucial for our proofs since Bertini-type theorems fail for foliations (cf. \cite[Example 3.4]{DLM23}). We need this notion to discuss the structures induced by MMP in a more accurate manner.

For $\bb$-divisors and generalized pairs, we will follow the notations and definitions in \cite{BZ16,HL23}. For generalized foliated quadruples, we shall follow \cite{CHLX23}. 

Generalized pairs and generalized foliated quadruples are very technical concepts. To make the statements in this paper more concise, for most results whose proofs for generalized foliated quadruples are similar to the proofs for foliated triples, we will only prove the foliated triple version and will not prove the generalized foliated quadruple version. We shall state the corresponding generalized foliated quadruple version in Appendix \ref{sec: gfq case}. We will freely use results in \cite{CHLX23} on generalized foliated quadruples.
\end{rem}

We need some results on NQC $\Rr$-divisors which are related to generalized pairs and generalized foliated quadruples.

\begin{defn}[NQC]
    Let $X\rightarrow U$ be a projective morphism from a normal quasi-projective variety to a variety. Let $D$ be a nef $\Rr$-Cartier $\Rr$-divisor on $X$ and $\Mm$ a nef $\bb$-divisor on $X$.

    We say that $D$ is NQC$/U$ if $D=\sum d_iD_i$, where each $d_i\geq 0$ and each $D_i$ is a nef$/U$ Cartier divisor. We say that $\Mm$ is NQC$/U$ if $\Mm=\sum \mu_i\Mm_i$, where each $\mu_i\geq 0$ and each $\Mm_i$ is a nef$/U$ Cartier $\bb$-divisor.
\end{defn}

\begin{lem}[{cf. \cite[Lemma 4.4(3)]{BZ16}}]\label{lem: bz16 4.4(3)}
    Let $(X,B)/U$ be a $\Qq$-factorial lc pair and $L$ an NQC$/U$ $\Rr$-divisor on $X$. Assume that $X$ is klt. Then there exists a positive real number $l_0$ such that any sequence of steps of a $(K_X+B+lL)$-MMP$/U$ is $L$-trivial for any $l>l_0$.
\end{lem}

\begin{lem}\label{lem: lc+ample nef=nqc}
    Let $(X,\Ff,B)/U$ be an lc algebraically integrable foliated triple and $D$ a nef$/U$ $\Rr$-divisor on $X$ such that $D-(K_{\Ff}+B)$ is ample$/U$. Then $D$ is NQC$/U$.
\end{lem}
\begin{proof}
We write $D=\sum_{i=1}^cr_iD_i$, where $r_1,\dots,r_c$ are linearly independent over $\mathbb Q$ and each $D_i$ is a Cartier divisor. Let $\bm{r}=(r_1,\dots,r_c)$. We define $D(\bm{v}):=\sum_{i=1}^cv_iD_i$ for any $\bm{v}=(v_1,\dots,v_c)\in\mathbb R^c$. By \cite[Lemma 5.3]{HLS24}, $D(\bm{v})$ is $\Rr$-Cartier for any $\bm{v}\in\mathbb R^c$. 

Let $L:=D-(K_\Ff+B)$. Since ample$/U$ is an open condition, there exists an open set $V\ni\bm{r}$ in $\mathbb R^c$, such that 
$\frac{1}{2}L+D(\bm{v})-D$ is ample$/U$ for any $\bm{v}\in V$.

By \cite[Theorem 2.3.1]{CHLX23}, there exist finitely many $\left(K_{\Ff}+B+\frac{1}{2}L\right)$-negative extremal rays$/U$ $R_1,\dots,R_l$, and $R_j=\mathbb R_+[C_j]$ for some curve $C_j$. Since $D$ is nef, $D\cdot C_j\geq 0$ for each $j$. Thus possibly shrinking $V$, we may assume that for any $\bm{v}\in V$, we have that $D(\bm{v})\cdot C_j>0$ for any $j$ such that $D\cdot C_j>0$. Since $r_1,\dots,r_c$ are linearly independent over $\mathbb Q$, for any $j$ such that $D\cdot C_j=0$, we have $D(\bm{v})\cdot C_j=0$ for any $\bm{v}\in\mathbb R^c$. Therefore,  $D(\bm{v})\cdot C_j\geq 0$ for any $j$ and any $\bm{v}\in V$.

By the cone theorem \cite[Theorem 2.3.1]{CHLX23}, for any curve $C$ on $X$, we may write
$[C]=\eta+\sum_{i=1}^l a_i[C_i]$
where $a_1,\dots,a_l\geq 0$ and $\eta\in\overline{NE}(X/U)_{K_{\Ff}+B+\frac{1}{2}L\geq 0}$. For any $\bm{v}\in V$, since
$$D(\bm{v})\cdot \eta=\left(K_\Ff+B+\frac{1}{2}L\right)\cdot \eta+\left(\frac{1}{2}L+D(\bm{v})-D\right)\cdot \eta\geq 0,$$
$D(\bm{v})\cdot C\geq 0$. Therefore, $D(\bm{v})$ is nef$/U$ for any $\bm{v}\in V$. We let $\bm{v}_1,\dots,\bm{v}_{c+1}\in V\cap\mathbb Q^c$ be rational points such that $\bm{r}$ is in the interior of the convex hull of $\bm{v}_1,\dots,\bm{v}_{c+1}$. Then there exist positive real numbers $a_1,\dots,a_{c+1}$ such that $\sum_{i=1}^{c+1}a_i=1$ and $\sum_{i=1}^{c+1}a_i\bm{v}_i=\bm{r}$. Since $D(\bm{v}_i)$ is a nef$/U$ $\Qq$-divisor for each $i$ and $D=\sum_{i=1}^{c+1}a_iD(\bm{v}_i)$, $D$ is NQC$/U$.
\end{proof}

\begin{lem}\label{lem: trivial ray when -delta triple}
Let $(X,\Ff,B)/U$ be an lc algebraically integrable foliated triple and $D$ an $\Rr$-divisor on $X$ such that  $K_{\Ff}+B+D$ is NQC$/U$. Then there exists $\delta_0\in (0,1)$ such that for any $\delta\in (0,\delta_0)$, any $(K_{\Ff}+B+(1-\delta)D)$-non-positive extremal ray$/U$ is $(K_{\Ff}+B+D)$-trivial.
\end{lem}
\begin{proof}
Let $\pi: X\rightarrow U$ be the induced morphism. Since $K_{\Ff}+B+D$ is NQC$/U$, there exists a positive real number $\epsilon$ such that $(K_{\Ff}+B+D)\cdot C\geq\epsilon$ for any curve $C$ such that $\pi(C)=\{pt\}$ and $(K_{\Ff}+B+D)\cdot C>0$. 

Let $d:=\dim X$. We show that $\delta_0:=\frac{\epsilon}{2d+\epsilon}$ satisfies our requirements. Let $R$ be a $(K_{\Ff}+B+(1-\delta)D)$-non-positive extremal ray$/U$. If $R$ is not $(K_{\Ff}+B+D)$-trivial, then $R$ is $(K_{\Ff}+B+D)$-positive, hence $(K_{\Ff}+B)$-negative. By the length of extremal rays \cite[Theorem 2.3.1]{CHLX23}, $R$ is spanned by a curve $C$ such that $\pi(C)=\{pt\}$ and $0<-(K_{\Ff}+B)\cdot C\leq 2d$. Thus for any $\delta\in (0,\delta_0)$,
\begin{align*}
   0&\ge(K_{\Ff}+B+(1-\delta)D)\cdot C=(1-\delta)(K_{\Ff}+B+D)\cdot C+\delta(K_{\Ff}+B)\cdot C\\
   &\geq (1-\delta)\epsilon-2d\delta>\epsilon-(2d+\epsilon)\delta_0=0,
\end{align*}
which is not possible.
\end{proof}

\section{Core models}\label{sec: core model}

The goal of this section is to introduce two new types of birational models for algebraically integrable foliations, namely \emph{core models} and \emph{simple models}. We shall also recall the definition of \emph{ACSS models} defined in \cite{CHLX23,DLM23}.

Roughly speaking, \emph{simple models} are birational models of algebraically integrable foliations that are weaker than ``$(*)$-models but still have potentially lc pair structures, while \emph{core models} are unique simple models which satisfy certain universal property. The use of core models is crucial for the proof of our main theorems. 

\subsection{Core models for two contractions}

\begin{deflem}\label{deflem: relative ample model}
    Let $X',X,Z$ be normal quasi-projective varieties and $h: X'\rightarrow X$, $f: X'\rightarrow Z$ contractions. Then there exists a unique normal quasi-projective variety $\bar X$ up to isomorphisms with two contractions $\bar h: \bar X\rightarrow X$ and $\bar f: \bar X\rightarrow Z$ satisfying the following.
    \begin{enumerate}
        \item  For any ample  $\Rr$-divisor $A$ on $X$, $\bar h^*A$ is ample$/Z$.
        \item  There exists a contraction $g: X'\rightarrow \bar X$ such that $\bar h\circ g=h$ and $\bar f\circ g=f$.
    \end{enumerate}
The variety $\bar X$ is called the \emph{core model} of $(h,f)$ associated with $(\bar h,\bar f)$.

Moreover, for any dominant map $\phi: X'\dashrightarrow X''$ such that $K(X'')$ is algebraically closed in $K(X')$ (e.g. $\phi$ is a birational map or a contraction), and any contractions $h'': X''\rightarrow X$ and $f'': X''\rightarrow Z$ such that $h''\circ\phi=h$ and $f''\circ\phi=f$, $\bar X$ is the core model of $(h'',f'')$ associated with $(\bar h,\bar f)$.
\end{deflem}

\begin{proof}

\noindent\textbf{Step 1}. In this step we construct $\bar h,\bar f, g$ which satisfy (2) such that (1) holds for a fixed ample $\Rr$-divisor $A$, and $\bar X$ is unique. Let $A$ be a fixed ample $\Rr$-divisor on $X$. Then $h^*A$ is semi-ample, hence semi-ample$/Z$. Let $g: X'\rightarrow\bar X$ be the ample model$/Z$ of $h^*A$. Then there exists an induced contraction $\bar f: \bar X\rightarrow Z$. Since the ample model of $h^*A$ is $h: X'\rightarrow X$, we have an induced contraction $\bar h: \bar X\rightarrow X$. We denote it by $\bar h$. By the uniqueness of ample models, $\bar X$ is unique.

\medskip

\noindent\textbf{Step 2}. In this step we show that (1) holds for any ample $\Rr$-divisor. Suppose that there exists an ample $\Rr$-divisor $H$ on $X$ such that $\bar h^*H$ is not ample$/Z$. Then by applying \textbf{Step 1} to $H,\bar h,\bar f$, there exist contractions $g': \bar X\rightarrow\bar X'$, $\bar h': \bar X'\rightarrow X$, and $\bar f': \bar X'\rightarrow Z$ such that $\bar h'^*H$ is ample$/Z$. Since $\bar h^*H$ is not ample$/Z$, $g'$ is not an isomorphism. Since $g'$ is a contraction$/Z$ and $\bar h^*A=g'^*\bar h'^*A$ is ample$/Z$, $g'$ is finite and thus an isomorphism, which is a contradiction.

\medskip

\noindent\textbf{Step 3}. In this step we prove the moreover part. There exist a birational morphism $p: W\rightarrow X'$ from a normal quasi-projective variety $W$ and a projective surjective morphism $q: W\rightarrow X''$ such that $q=\phi\circ p$. Since $K(X'')$ is algebraically closed in $K(X')$, $q$ is a contraction. Let $\bar X''$ be the core model of $(h'',f'')$ associated with $(\bar h'',\bar f'')$ and $g'': X''\rightarrow\bar X''$ the induced contraction. Since both $g\circ p$ and $g''\circ q$ are ample models$/Z$ of  $(h\circ p)^*A=(h''\circ q)^*A$, the ample model$/Z$ of $h^*A$ is isomorphic to the ample model$/Z$ of $h''^*A$. Thus $\bar X$ is the core model of $(h'',f'')$ associated with $(\bar h,\bar f)$.
\end{proof}

\subsection{Core models for algebraically integrable foliations}

\begin{defn}[Simple modifications]\label{defn: simple model}
    Let $(X,\Ff,B)$ and $(X',\Ff',B')$ be two algebraically integrable foliated triples and $h: X'\rightarrow X$ a birational morphism. Let $f: X'\rightarrow Z$ be a contraction and $G$ a reduced divisor on $X'$. We say that $h: (X',\Ff',B';G)/Z\rightarrow (X,\Ff,B)$ is a \emph{simple modification} if the following conditions hold.
    \begin{enumerate}
        \item $\Ff':=h^{-1}\Ff$ and $B':=h^{-1}_*(B^{\Ff-\mathrm{ninv}}\wedge\Supp B)+\Exc(h)^{\Ff'-\mathrm{ninv}}$. 
        \item $(X',\Ff',B')$ is lc.
        \item  $a(E,\Ff,B)\leq-\epsilon_{\Ff}(E)$ for any $h$-exceptional prime divisor $E$. 
        \item $K_{\Ff'}+B'\sim_ZK_{X'}+B'+G.$
        \item $(X',\Ff',B';G)/Z$ satisfies Property $(*)$.
    \end{enumerate}
    
    We say that $h: (X',\Ff',B',G)/Z\rightarrow (X,\Ff,B)$ is an \emph{ACSS modification} if it is a simple modification and  $(X',\Ff',B';G)/Z$ is ACSS.

    We say that $h: (X',\Ff',B',G)/Z\rightarrow (X,\Ff,B)$ is a \emph{core modification} if it is a simple modification and $h^*A$ is ample$/Z$ for any ample $\Rr$-divisor $A$ on $X$. 

    If $h: (X',\Ff',B',G)/Z\rightarrow (X,\Ff,B)$ is a simple (resp. ACSS, core) modification, then we also say that $h$ is a simple (resp. ACSS, core) modification of $(X,\Ff,B)$.
\end{defn}

\begin{defn}[Core models and ACSS models]\label{defn: core model}
Let $(X,\Ff,B)$ be an algebraically integrable foliated triple and $h: (X',\Ff',B',G)/Z\rightarrow (X,\Ff,B)$ a simple (resp. ACSS, core) modification. 

We say that $(X',\Ff',B';G)/Z$, $(X',\Ff',B')/Z$, and $(X',\Ff',B')$ are \emph{simple} (resp. \emph{ACSS, core}) \emph{models} of $(X,\Ff,B)$. Moreover, we say that $h$ is
\begin{enumerate}
\item \emph{$\Qq$-factorial} if $X'$ is $\Qq$-factorial,
\item \emph{strict} if the $\Ff$-invariant part of $\Supp\Exc(h)$ is contained in $\Supp G$,
\item \emph{super} if $G$ is super$/Z$.
\end{enumerate}
\end{defn}

We will frequently use the following result on the existence of ACSS models:

\begin{thm}[{\cite[Theorem 2.5.1]{CHLX23},\cite[Theorem 3.10]{ACSS21}}]\label{thm: eo acss model}
    Let $(X,\Ff,B)$ be an lc algebraically integrable foliated triple. Then $(X,\Ff,B)$ has an ACSS model $h: (X',\Ff',B';G)/Z\rightarrow (X,\Ff,B)$ that is $\Qq$-factorial, strict and super.
\end{thm}

\begin{lem}\label{lem: existence of core model}
    Let $(X,\Ff,B)$ be an lc algebraically integrable foliated triple and let $h: (X',\Ff',B';G)/Z\rightarrow (X,\Ff,B)$ be a simple model. Let $f: X'\rightarrow Z$ be the associated contraction and let $\bar X$ be the core model of $(h,f)$ associated with $(\bar h,\bar f)$. Let $g: X'\rightarrow\bar X$ be the induced birational morphism, $\bar\Ff:=g_*\Ff',\bar B:=g_*B'$, and $\bar G:=g_*G$. 
    
    Assume that $f$ is equidimensional. Then:
    \begin{enumerate}
    \item $K_{\Ff'}+B'=g^*(K_{\bar\Ff}+\bar B)$.
        \item $K_{X'}+B'+G=g^*(K_{\bar X}+\bar B+\bar G)$.
        \item $\bar h: (\bar X,\bar\Ff,\bar B;\bar G)/Z\rightarrow (X,\Ff,B)$
        is a core model.
        \item If $h: (X',\Ff',B';G)/Z\rightarrow (X,\Ff,B)$ is strict (resp. super), then $\bar h: (\bar X,\bar\Ff,\bar B;\bar G)/Z\rightarrow (X,\Ff,B)$ is strict (resp. super).
    \end{enumerate}
\end{lem}
\begin{proof}
(1) Let
    $K_{\bar\Ff}+\bar B':=\bar h^*(K_{\Ff}+B).$
    Then
    $$K_{\bar \Ff}+\bar B=g_*(K_{\Ff'}+B')=g_*h^*(K_{\Ff}+B)=g_*g^*(K_{\bar\Ff}+\bar B')=K_{\bar\Ff}+\bar B'.$$
    Therefore, $\bar B'=\bar B$, so $K_{\Ff'}+B'=g^*(K_{\bar\Ff}+\bar B)$.

(2) By Proposition \ref{prop: weak cbf gfq}, we have
$$K_{X'}+B'+G\sim_ZK_{\Ff'}+B'=g^*(K_{\bar\Ff}+\bar B)\sim_{\mathbb R,\bar X}0.$$
Since $g$ is a birational morphism$/Z$, 
$$K_{X'}+B'+G\sim_{\mathbb R,\bar X}0.$$
Then we have
$$K_{X'}+B'+G=g^*g_*(K_{X'}+B'+G)=g^*(K_{\bar X}+\bar B+\bar G).$$

(3) Since $\bar X$ is the core model of $(h,f)$, we only need to show that $\bar h: (\bar X,\bar\Ff,\bar B;\bar G)/Z\rightarrow (X,\Ff,B)$ is a simple model by checking each condition of Definition \ref{defn: simple model}.

Definition \ref{defn: simple model}(1):  By our construction, $\bar\Ff=\bar h^{-1}\Ff$ and $\bar B=\bar h^{-1}_*B+\Exc(\bar h)^{\bar\Ff-\mathrm{ninv}}$.

Definition \ref{defn: simple model}(2): By (1), $K_{\bar\Ff}+\bar B=\bar h^*(K_{\Ff}+B)$. Since $(X,\Ff,B)$ is lc, $(\bar X,\bar\Ff,\bar B)$ is lc.

Definition \ref{defn: simple model}(3): Since any $\bar h$-exceptional divisor is also $h$-exceptional, it follows from our construction. 

Definition \ref{defn: simple model}(4): It follows from (1) and (2).

Definition \ref{defn: simple model}(5): We only need to check Definition \ref{defn: foliation property *} for $(\bar X,\bar\Ff,\bar B,\bar G)/Z$. Definition \ref{defn: foliation property *}(1): Since $\Ff'$ is induced by $f$, $\bar\Ff$ is induced by $\bar f$. Since $G$ is $\Ff'$-invariant, $\bar G$ is $\bar\Ff$-invariant. Definition \ref{defn: foliation property *}(2): $(Z,f(G)=\bar f(\bar G))$ is log smooth by assumption. Since $G=f^{-1}(f(G))$, $\bar G=\bar f^{-1}(\bar f(G))$. Definition \ref{defn: foliation property *}(3):  For any closed point $z\in Z$ and any reduced divisor $\Sigma\geq\bar f(\bar G)$ on $Z$ such that $(Z,\Sigma)$ is log smooth near $z$, $(X',B'+f^*(\Sigma-\bar f(\bar G)))$ is lc over a neighborhood of $z$. By (2),
$$K_{X'}+B'+f^*(\Sigma-f(G))=g^*(K_{\bar X}+\bar B+\bar f^*(\Sigma-\bar f(\bar G))),$$
so $(\bar X,\bar B+\bar f^*(\Sigma-\bar f(\bar G)))$ is lc over a neighborhood of $z$.

(4) It follows from the definitions of being strict or super.
\end{proof}

\begin{lem}\label{lem: existence of klt pair on strict simple models}
       Let $(X,\Ff,B)$ be an lc algebraically integrable foliated triple and let $h: (X',\Ff',B';G)/Z\rightarrow (X,\Ff,B)$ be a strict simple model. If $X$ is potentially klt, then $X'$ is potentially klt.
\end{lem}
\begin{proof}
Let $(X,\Delta)$ be a klt pair and let $K_{X'}+\tilde\Delta':=h^*(K_X+\Delta)$. Then $(X',\tilde\Delta')$ is sub-klt. Since $(X',B'+G)$ is lc, $(X',\delta\tilde\Delta'+(1-\delta)(B'+G))$ is sub-klt for any $\delta\in (0,1)$. Since $G$ contains all $h$-exceptional prime divisors that are $\Ff'$-invariant, and since $\Supp B'$ contains all $h$-exceptional prime divisors that are not $\Ff'$-invariant, we have $\delta\tilde\Delta'+(1-\delta)(B'+G)\geq 0$ for any $0<\delta\ll 1$. Therefore, $(X',\delta\tilde\Delta'+(1-\delta)(B'+G))$ is klt for any $0<\delta\ll 1$. In particular, $X'$ is potentially klt.
\end{proof}

\begin{rem}
    \cite{ACSS21,CS25a} define ``$(*)$-models" and \cite{CHLX23} defines ``great ACSS models". We do not need these models in this paper. Nevertheless, we provide the readers with the following table on the properties of different types of models. Note that ``$(*)$-models" in \cite{ACSS21} and \cite{CS25a} are defined differently.
\begin{table}[ht]
   \caption{Different types of simple models}
    \label{tab:log structures}
\begin{center}
        \begin{tabular}{|c|c|c|c|c|c|c|}
\hline
& $\Qq$-factorial & Strict & Super & $X'\rightarrow Z$ equidim & $(X',B'+G)$ & $X'$ \\
\hline
Simple models & N &  N & N & N  & lc & lc$^\bullet$\\
\hline
Core models & N & N & N & N & lc & lc$^\bullet$\\
\hline
$(*)$-models (\cite{ACSS21}) & N & N & N & Y & lc & klt \\
\hline
ACSS models & N & N & N & Y & qdlt& klt$^\bullet$\\
\hline
$(*)$-models (\cite{CS25a}) & Y & Y & N & Y & lc & klt \\
\hline
Great ACSS models & N & Y & Y & Y & qdlt & klt$^\bullet$ \\
\hline
\end{tabular}
\end{center}
Y: Yes. N: Not necessarily true. $^\bullet$: holds when $X'$ is $\Qq$-factorial.
\end{table}
\end{rem}

\section{Models for foliations}\label{sec: models}

The goal of this section is to introduce and study the basic behaviors of different types of models for foliated triples: weak lc models, minimal models, good minimal models, etc. We will also introduce minimal models in the sense of Birkar-Shokurov and log minimal models for foliations. Results in this section are similar to results in \cite[Section 2]{Bir12} and \cite[Section 3]{HL23} with some differences as we need to take invariant lc centers into consideration.

\subsection{Definitions of minimal models and Mori fiber spaces}

\begin{rem}
    In the classical definition of models, ``log minimal model",``good minimal model", or ``log terminal model" (cf. \cite{BCHM10,Bir12}) usually requires that the model is $\Qq$-factorial dlt. This is because the initial structure on which we start running the MMP is usually $\Qq$-factorial dlt. For foliations this is replaced by the condition ``$\Qq$-factorial ACSS" \cite{CHLX23}. However, since the singularities we are going to come up with in this paper is usually worse than $\Qq$-factorial ACSS, we have to change these definitions a little bit in order to deal with worse singularities. On the other hand, we still want to consider models of the objects we study that have nice singularities after possibly extracting some lc places. Considering all these issues, we will slightly change the notations in previous literature and define different models in the following way:
    \begin{itemize}
        \item For models requiring good singularities (e.g. $\Qq$-factorial ACSS), we always keep the word ``log". We always allow extraction of lc centers when considering these models.
        \item For models without these strict singularity conditions (e.g. only requiring lc), we shall not use the word ``log". Moreover, if we allow extraction of lc centers, then we shall add the prefix ``bs-" or write ``in the sense of Birkar-Shokurov".
    \end{itemize}
\end{rem}

\begin{defn}[Log birational model]\label{defn: log birational model}
  Let $(X,\Ff,B)/U$ be a foliated triple, $\phi: X\dashrightarrow X'$ a birational map over $U$, $E:=\Exc(\phi^{-1})$ the reduced $\phi^{-1}$-exceptional divisor, and $\Ff':=\phi_*\Ff$. Assume that $a(D,\Ff,B)\leq-\epsilon_{\Ff}(D)$ for any component $D$ of $E$. We let
  $$B':=\phi_*B+\sum_D(-a(D,\Ff,B))E$$
and say that $(X',\Ff',B')/U$ is a \emph{log birational model} of $(X,\Ff,B)/U$, where the sum runs through all components of $E$.
\end{defn}

\begin{defn}[Minimal models]\label{defn: minimal model}
    Let $(X,\Ff,B)/U$ be a foliated triple and $(X',\Ff',B')/U$ a log birational model of $(X,\Ff,B)/U$ such that $K_{\Ff'}+B'$ is nef$/U$. 
    \begin{enumerate}
        \item We say that $(X',\Ff',B')/U$ is a \emph{bs-weak lc model} or \emph{weak lc model in the sense of Birkar-Shokurov} of $(X,\Ff,B)/U$, if for any prime divisor $D$ on $X$ which is exceptional over $X'$, $$a(D,\Ff,B)\leq a(D,\Ff',B').$$
        \item We say that $(X',\Ff',B')/U$ is a \emph{bs-minimal model} or \emph{minimal model in the sense of Birkar-Shokurov} of $(X,\Ff,B)/U$, if for any prime divisor $D$ on $X$ which is exceptional over $X'$, $$a(D,\Ff,B)<a(D,\Ff',B').$$
          \item We say that $(X',\Ff',B')/U$ is a \emph{bs-semi-ample model} or \emph{semi-ample model in the sense of Birkar-Shokurov} of $(X,\Ff,B)/U$ if it is a bs-weak lc model of $(X,\Ff,B)/U$ and $K_{\Ff'}+B'$ is semi-ample$/U$. 
        \item We say that $(X',\Ff',B')/U$ is a \emph{bs-good minimal model} or \emph{good minimal model in the sense of Birkar-Shokurov} of $(X,\Ff,B)/U$ if it is a bs-minimal model of $(X,\Ff,B)/U$ and $K_{\Ff'}+B'$ is semi-ample$/U$. 
        \end{enumerate}
If, in addition, the induced birational map $X\dashrightarrow X'$ does not extract any divisor, then we say remove the initial ``bs-" or the phrase ``in the sense of Birkar-Shokurov" in the previous definitions. 
\begin{enumerate}
        \item[(5)] We say that $(X',\Ff',B')/U$ is a \emph{log minimal model} of $(X,\Ff,B)/U$ if it is a bs-minimal model of $(X,\Ff,B)$ and $(X',\Ff',B')$ is $\Qq$-factorial ACSS.
        \item[(6)] We say that $(X',\Ff',B')/U$ is a \emph{good log minimal model} of $(X,\Ff,B)/U$ if it is a log minimal model of $(X,\Ff,B)$ and $K_{\Ff'}+B'$ is semi-ample$/U$. 
\end{enumerate}

We remark that, similar to \cite{CHLX23}, the definition of ``log minimal model" in our paper does not coincide with the classical definition with $\Ff=T_X$ as $\Qq$-factorial ACSS is equivalent to $\Qq$-factorial qdlt instead of $\Qq$-factorial dlt when $\Ff=T_X$. 
\end{defn}

\begin{defn}[Mori fiber space]\label{defn: mfs}
Let $(X,\Ff,B)/U$ be a foliated triple and let $(X',\Ff',B')/U$ be a log birational model of $(X,\Ff,B)/U$. Let $f: X'\rightarrow Z$ be a $(K_{\Ff'}+B')$-Mori fiber space$/U$.
\begin{enumerate}
\item We say that $(X',\Ff',B')\rightarrow Z$ is a \emph{bs-Mori fiber space}, or a \emph{Mori fiber space in the sense of Birkar-Shokurov} of $(X,\Ff,B)/U$, if for any prime divisor $D$ on $X$ which is exceptional over $X'$, 
$$a(D,\Ff,B)<a(D,\Ff',B').$$
\item We say that $(X',\Ff',B')\rightarrow Z$ is a \emph{Mori fiber space} of $(X,\Ff,B)/U$ if $(X',\Ff',B')\rightarrow Z$  is a bs-Mori fiber space of $(X,\Ff,B)/U$ and the induced birational map $X\dashrightarrow X'$ does not extract any divisor.
\item We say that  $(X',\Ff',B')\rightarrow Z$ is a \emph{log Mori fiber space} of $(X,\Ff,B)/U$ if it is a bs-Mori fiber space of $(X,\Ff,B)/U$ and $(X',\Ff',B')$ is $\Qq$-factorial ACSS.
\end{enumerate}
\end{defn}

\begin{rem}
The condition ``$\Qq$-factorial ACSS" is a condition only for algebraically integrable foliations. Therefore, ``log minimal model", ``good log minimal model", and ``log Mori fiber space" are only well-defined for algebraically integrable foliations. However, Definition \ref{defn: minimal model}(1-4) and Definitions \ref{defn: mfs}(1-2) are well-defined for arbitrary foliations. Therefore, many results in this section also hold for arbitrary foliations.

We also remark that we do not have any requirement on the singularities of $(X,\Ff,B)$ and $(X',\Ff',B')$ in  Definition \ref{defn: minimal model}(1-4) and Definitions \ref{defn: mfs}(1-2). This is because in many cases, we want to consider a generalized foliated quadruple polarized by an ample divisor $A$. Due to the failure of Bertini-type theorems for foliations, usually the only thing we can do is to consider a generalized foliated quadruple structure $(X,\Ff,B,\bar A)$, i.e. we let $A$ be the nef part. However, this is inconvenient when the foliation is associated with some other pair structure, as many theorems on pairs consider structures of the form $(X,B+A)$ instead. Therefore, if we do not have any singularity restrictions on the models, then using $(X,\Ff,B+A)$ will bring us more flexibility when applying results of usual pairs.
\end{rem}

\subsection{Basic properties of models}

In this subsection we prove several basics properties on models of foliated triples. We remark that results in this section works for any foliated triples without any requirement on algebraic integrability nor singularities, so we expect results in this subsection to be useful for further applications, particularly to non-algebraically integrable foliations.

\begin{lem}[cf.{ \cite[Remark 2.6]{Bir12}, \cite[Lemma 3.4]{HL23}}]\label{lem: hl23 2.6}
Let $(X,\Ff,B)/U$ be a foliated triple and let $(X',\Ff',B')/U$ a bs-weak lc model of $(X,\Ff,B)/U$ associated with the birational map $\phi: X\dashrightarrow X'$. Let $p: W\rightarrow X$ and $q: W\rightarrow X'$ be birational morphisms such that $q=\phi\circ p$. Assume that
$$p^*(K_\Ff+B)=q^*(K_{\Ff'}+B')+E,$$
then $E\geq 0$ and is exceptional$/X'$.
\end{lem}
\begin{proof}
For any prime divisor $D$ that is an irreducible component of $E$, $$\mult_DE=a(D,\Ff',B')-a(D,\Ff,B).$$ 
Therefore, if $D$ is not exceptional$/X$, then:
\begin{itemize}
    \item If $D$ is not exceptional$/X'$, then $\mult_DE=0$ by Definition \ref{defn: log birational model}.
    \item If $D$ is exceptional$/X'$, then $\mult_DE\geq 0$ by Definition \ref{defn: minimal model}(1).
\end{itemize}
Therefore, $p_*E\geq 0$. Since $K_{\Ff'}+B'$ is nef$/U$, $q^*(K_{\Ff'}+B')$ is nef$/X$, hence $E$ is anti-nef$/X$. By the negativity lemma, $E\geq 0$.

If $E$ is not exceptional$/X'$, then there exists a component $D$ of $E$ that is not exceptional$/X'$. If $D$ is not exceptional$/X$, then $\mult_DE=0$ by Definition \ref{defn: log birational model}, a contradiction. Thus $D$ is exceptional over $X$. In particular, $\phi$ extracts $D$. Since $(X',\Ff',B')/U$ is a log birational model of $(X,\Ff,B)$, $$a(D,\Ff',B')=a(D,\Ff,B),$$ 
which implies that $\mult_DE=0$, a contradiction.
\end{proof}

\begin{lem}[cf.{ \cite[Remark 2.7]{Bir12}, \cite[Lemma 3.5]{HL23}}]\label{lem: g-pair version bir12 2.7}
Let $(X,\Ff,B)/U$ be a foliated triple. Let $(X_1,\Ff_1,B_1)/U$ and $(X_2,\Ff_2,B_2)/U$ be two bs-weak lc models of $(X,\Ff,B)/U$ with induced birational maps $\phi: X_1\dashrightarrow X_2$. Let $h_1: W\rightarrow X_1$ and $h_2: W\rightarrow X_2$ be two birational morphisms such that $\phi\circ h_1=h_2$. Then:
\begin{enumerate}
    \item $$h_1^*(K_{\Ff_1}+B_1)=h_2^*(K_{\Ff_2}+B_2).$$
    \item If $K_{\Ff_2}+B_2$ is semi-ample$/U$, then $K_{\Ff_1}+B_1$ is semi-ample$/U$.
    \item If $K_{\Ff_2}+B_2$ is ample$/U$, then $\phi$ is a morphism.
\end{enumerate}
\end{lem}
\begin{proof}
Let $\phi_1: X\dashrightarrow X_1$ and $\phi_2: X\dashrightarrow X_2$ be the induced birational maps. Possibly replacing $W$ with a higher model, we may assume that the induced birational map $h: W\rightarrow X$ is a morphism. Let $$E_i:=h^*(K_X+B)-h_i^*(K_{X_i}+B_i)$$
for $i\in\{1,2\}$. By Lemma \ref{lem: hl23 2.6}, $E_i\geq 0$ and is exceptional over $X_i$ for $i\in\{1,2\}$. Thus $h_{1,*}(E_2-E_1)\geq 0$ and $E_1-E_2$ is nef$/X_1$, and $h_{2,*}(E_1-E_2)\geq 0$ and $E_2-E_1$ is nef$/X_2$. By the negativity lemma, $E_2-E_1\geq 0$ and $E_1-E_2\geq 0$. Thus $E_1=E_2$, which implies (1). (2) immediately follows from (1). By (1), if $K_{\Ff_2}+B_2$ is ample$/U$, then $h_2: W\rightarrow X_2$ is the ample model$/U$ of $h^*(K_{\Ff_1}+B_1)$, hence $\phi$ is the ample model$/U$ of $K_{\Ff_1}+B_1$. Since $K_{\Ff_1}+B_1$ is semi-ample$/U$, $\phi$ is a morphism. This implies (3).
\end{proof}

\begin{lem}\label{lem: numerical equivalence model}
    Let $r$ be a positive real number. Let $(X,\Ff_1,B_1)/U$ and $(X,\Ff_2,B_2)/U$ be two foliated triples such that
    $$K_{\Ff_2}+B_2\equiv_U r(K_{\Ff_1}+B_1)$$
    Let $(X',\Ff'_1,B_1')/U$ be a weak lc model (resp. minimal model) of $(X,\Ff_1,B_1)/U$ with induced birational map $\phi: X\dashrightarrow X'$. Let $\Ff_2':=\phi_*\Ff$ and $B_2':=\phi_*B_2$. Then $(X',\Ff_2',B_2')/U$ is a weak lc model (resp. minimal model) of $(X,\Ff_2,B_2)/U$.

    If  $(X',\Ff'_1,B_1')/U$ is a semi-ample model (resp. good minimal model) of $(X,\Ff_1,B_1)/U$ and
    $$K_{\Ff_2}+B_2\sim_{\mathbb R,U} r(K_{\Ff_1}+B_1),$$
        $(X',\Ff'_2,B_2')/U$ is a semi-ample model (resp. good minimal model) of $(X,\Ff_2,B_2)/U$.
\end{lem}
\begin{proof}
Let $p: W\rightarrow X$ and $q: W\rightarrow X'$ be a resolution of indeterminacy. By Lemma \ref{lem: hl23 2.6},
$$p^*(K_{\Ff_1}+B_1)=q^*(K_{\Ff_1'}+B_1')+E$$
for some $\Rr$-divisor $E\geq 0$ that is exceptional$/X'$. Then
$$p^*(K_{\Ff_2}+B_2)\equiv_U rq^*(K_{\Ff_1'}+B_1')+rE,$$
so
$$K_{\Ff_2'}+B_2'=q_*p^*(K_{\Ff_1}+B_1)\equiv q_*(rq^*(K_{\Ff_1'}+B_1')+rE)=r(K_{\Ff_1'}+B_1')$$
is nef$/U$. Moreover, if $K_{\Ff_1'}+B_1'$ is semi-ample$/U$ and $K_{\Ff_2}+B_2\sim_{\mathbb R,U} r(K_{\Ff_1}+B_1)$, then 
$$K_{\Ff_2'}+B_2'\sim_{\mathbb R,U}r(K_{\Ff_1'}+B_1')$$
is semi-ample$/U$. 

We have
$$p^*(K_{\Ff_2}+B_2)\equiv_U q^*(K_{\Ff_2'}+B_2')+rE.$$
Therefore, for any prime divisor $D$ on $X$ which is exceptional over $X'$,
$$a(D,\Ff_2',B_2')-a(D,\Ff_2,B_2)=-\mult_Dp_*(rE)=r(a(D,\Ff_1',B_1')-a(D,\Ff_1,B_1)).$$
Therefore, 
$a(D,\Ff_2,B_2)\leq\text{(resp. }<\text{) }a(D,\Ff_2',B_2')$ if and only if $a(D,\Ff_1,B_1)\leq\text{(resp. }<\text{) }a(D,\Ff_1',B_1')$. The lemma follows immediately from the definitions.
\end{proof}

\subsection{Models under foliated log resolutions}

From now on, we shall focus on different models of foliations that are algebraically integrable and lc. We first study the relationship between different types of models and foliated log resolutions. Of course, we expect results in this section to hold in greater generalities provided that there is a proper definition of ``foliated log resolution" for non-algebraically integrable foliations.

We first recall the following result.

\begin{thm}[{\cite[Theorem 9.4.1]{CHLX23}}]\label{thm: chlx23 9.4.1}
    Let $(X,\Ff,B)/U$ be a $\Qq$-factorial ACSS algebraically integrable foliated triple such that $K_{\Ff}+B\sim_{\mathbb R,U}E\geq 0$ and $E$ is very exceptional$/U$. Then we may run a $(K_{\Ff}+B)$-MMP$/U$ with scaling of an ample$/U$ $\Rr$-divisor $A$ and any such MMP terminates with a good log minimal model $(X',\Ff',B')/U$ such that $K_{\Ff'}+B'\sim_{\mathbb R,U}0$.
\end{thm}

See \cite[Definition 3.1]{Bir12} for the definition of very exceptional divisors. In particular, exceptional divisors coincide with very exceptional divisors in the context of birational morphisms.

\begin{defn}[Foliated log smooth model]\label{defn: log smooth models}
Let $(X,\Ff,B)$ be an lc algebraically integrable foliated triple and $h: X'\rightarrow X$ a foliated log resolution of $(X,\Ff,B)$ (cf. Definition \ref{defn: log resolution}). Let $\Ff':=h^{-1}\Ff$, and let $B'\geq 0$ and $E\geq 0$ be two $\Rr$-divisors on $X'$ satisfying the following.
\begin{enumerate}
    \item $K_{\Ff'}+B'=h^*(K_\Ff+B)+E$,
    \item $(X',\Ff',B')$ is foliated log smooth and lc.
    \item $E$ is $h$-exceptional.
    \item For any $h$-exceptional prime divisor $D$ such that $$a(D,X,B)>-\epsilon_{\Ff}(D),$$ $D$ is a component of $E$.
\end{enumerate}
We say that $(X',\Ff',B')$ is a \emph{foliated log smooth model} of $(X,\Ff,B)$. 
\end{defn}

\begin{lem}\label{lem: g-pair version bir12 2.8}
Let $(X,\Ff,B)/U$ be an lc algebraically integrable foliated triple. Let $(W,\Ff_W,B_W)$ be a foliated log smooth model of $(X,\Ff,B)$. 

Then any bs-weak lc model (resp. bs-minimal model, bs-semi-ample model, bs-good minimal model, log minimal model, good log minimal model) of $(W,\Ff_W,B_W)/U$ is a bs-weak lc model (resp. bs-minimal model, bs-semi-ample model, bs-good minimal model, log minimal model, good log minimal model) of $(X,\Ff,B)/U$. 
\end{lem}

\begin{proof}
We let $h: W\rightarrow X$ be the induced birational morphism. We may write
$$K_{\Ff_W}+B_W=h^*(K_{\Ff}+B)+E$$
for some $E\geq 0$ that is $h$-exceptional, and $D\subset\Supp E$ for any $h$-exceptional prime divisor $D$ such that $a(D,X,B)>-\epsilon_{\Ff}(E)$.

\begin{claim}\label{claim: log smooth model log discrepancy compare}
Let $(X',\Ff',B')/U$ be a bs-weak lc model of $(W,\Ff_W,B_W)/U$. Then $$a(D,\Ff,B)\leq a(D,\Ff',B')$$ for any prime divisor $D$ over $X$.
\end{claim}
\begin{proof}
Let $\phi_W: W\dashrightarrow X'$ be the induced birational map, and let $p: V\rightarrow W$ and $q: V\rightarrow X'$ be a common resolution such that $q=\phi_W\circ p$. By Lemma \ref{lem: hl23 2.6},
$$p^*(K_{\Ff_W}+B_W)=q^*(K_{\Ff'}+B')+F$$
for some $F\geq 0$ that is exceptional over $X'$. Then we have
$$p^*h^*(K_{\Ff}+B)=q^*(K_{\Ff'}+B')+F-p^*E,$$
so
$$p^*E-F\sim_{\Rr,X}q^*(K_{\Ff'}+B')$$
is nef$/X$. Since $h_*p_*(F-p^*E)=h_*p_*F\geq 0$, by the negativity lemma, $F\geq p^*E$. Thus $a(D,\Ff,B)\leq a(D,\Ff',B')$ for any prime divisor $D$ over $X$.
\end{proof}

\noindent\textit{Proof of Lemma \ref{lem: g-pair version bir12 2.8} continued}. First we prove the bs-weak lc model case. Let $(X',\Ff',B')/U$ be a bs-weak lc model of $(W,\Ff_W,B_W)/U$ with induced birational map $\phi_W: W\dashrightarrow X'$.  By Claim \ref{claim: log smooth model log discrepancy compare}, we only need to show that $(X',\Ff',B')/U$ is a log birational model of $(X,\Ff,B)/U$. 

Let $\phi: X\dashrightarrow X'$ be the induced morphism and 
$$\tilde B':=\phi_*B+\Exc(\phi^{-1})^{\Ff'-\mathrm{ninv}},$$ then we only need to show that $B'=\tilde B'$. Since  $(X',\Ff',B')/U$ is a bs-weak lc model of $(W,\Ff_W,B_W)/U$, we have
$$B'=(\phi_W)_*B_W+\Exc(\phi_W^{-1})^{\Ff'-\mathrm{ninv}}.$$ 
Let $D$ be a prime divisor on $X'$. There are three cases:

\medskip

\noindent\textbf{Case 1}. $D$ is not exceptional over $X$. In this case,
    $$-\mult_D\tilde B'=a(D,\Ff',\tilde B')=a(D,\Ff,B)=a(D,\Ff_W,B_W)=a(D,\Ff',B')=-\mult_DB',$$
so $\mult_DB'=\mult_D\tilde B'$.

\medskip

\noindent\textbf{Case 2}. $D$ is exceptional over $W$. In this case, $D$ is a component of $\Exc(\phi_W^{-1})$ and a component of $\Exc(\phi^{-1})$, hence
$$\mult_DB'=\epsilon_{\Ff}(D)=\mult_DB''.$$

\medskip

\noindent\textbf{Case 3}. $D$ is exceptional over $X$ but not exceptional over $W$. In this case,
$$-\mult_DB'=a(D,\Ff',B')=a(D,\Ff_W,B_W).$$
Since $E\geq 0$,
$a(D,\Ff_W,B_W)\leq a(D,\Ff,B).$
By Claim \ref{claim: log smooth model log discrepancy compare}, 
$a(D,\Ff,B)\leq a(D,\Ff',B').$
Thus
$$-\mult_DB'=a(D,\Ff,B)=a(D,\Ff',B')=a(D,\Ff_W,B_W).$$
By Definition \ref{defn: log smooth models}(4), $a(D,\Ff,B)=-\epsilon_{\Ff}(D),$
which implies that
$$\mult_DB'=\epsilon_{\Ff}(D)=\mult_D\Exc(\phi^{-1})^{\Ff'-\mathrm{ninv}}=\mult_D\tilde B'.$$
Thus $B'=\tilde B'$, so $(X',\Ff',B')/U$ is a log birational model of $(X,\Ff,B)/U$, and we are done for the bs-weak lc model case.

Next we prove the bs-minimal model case. Suppose that $(X',\Ff',B')/U$ be a bs-minimal model of $(W,\Ff_W,B_W)/U$. For any prime divisor $D$ on $X$ which is exceptional over $X'$, $h^{-1}_*D$ is a prime divisor on $W$ which is exceptional over $X'$. Thus
$$a(D,\Ff,B)=a(D,\Ff_W,B_W)<a(D,\Ff',B').$$
The bs-minimal model case immediately follows from the bs-weak lc model case. 

The bs-semi-ample model, bs-good minimal model, log minimal model, and good log minimal model cases follow immediately from the bs-weak lc model and the bs-minimal model cases.
\end{proof}

\subsection{Models under pullbacks}

\begin{lem}\label{lem: foliation lsm has lmm}
Let $(X,\Ff,B)/U$ be an lc algebraically integrable foliated triple and $(X',\Ff',B')/U$ a bs-weak lc model of $(X,\Ff,B)/U$. Let $(W,\Ff_W,B_W)$ be a foliated log smooth model of $(X,\Ff,B)$ such that the induced birational map $\phi_W: W\dashrightarrow X'$ is a morphism. 

Then we may run a $(K_{\Ff_W}+B_W)$-MMP$/X'$ with scaling of an ample$/X'$ $\Rr$-divisor which terminates with a good minimal model $(Y,\Ff_Y,B_Y)/X'$ of $(W,\Ff_W,B_W)/X'$ such that $$K_{\Ff_Y}+B_Y=q^*(K_{\Ff'}+B').$$
where $q: Y\rightarrow X'$ is the induced morphism. In particular, $(Y,\Ff_Y,B_Y)/U$ is a log minimal model of $(W,\Ff_W,B_W)/U$.
\end{lem}
\begin{proof}
    Let $h: W\rightarrow X$ be the induced birational morphism. We have
    $$K_{\Ff_W}+B_W=h^*(K_\Ff+B)+E$$
for some $E\geq 0$ that is exceptional$/X$. By Lemma \ref{lem: hl23 2.6}, we have
$$h^*(K_{\Ff}+B)=\phi_W^*(K_{\Ff'}+B')+F$$
where $F\geq 0$ is exceptional$/X'$. Thus
$$K_{\Ff_W}+B_W\sim_{\mathbb R,X'}F+E.$$
\begin{claim}\label{claim: wglc to lmm E exceptional}
$E$ is exceptional$/X'$.
\end{claim}
\begin{proof}
Let $D$ be a component of $E$. By Definition \ref{defn: log smooth models}(4), $a(D,\Ff,B)>-\epsilon_{\Ff}(E)$ and $D$ is exceptional$/X$. 

Assume that $D$ is not exceptional over $X'$. Since $(X',\Ff',B')/U$ is a log birational model of $(X,\Ff,B)/U$ and $(X,\Ff,B)$ is lc, $a(D,\Ff',B')=-\epsilon_{\Ff}(E)$. Since $F\geq 0$, $a(D,\Ff,B)\leq a(D,\Ff',B')$. Thus $a(D,\Ff,B)=-\epsilon_{\Ff}(E)$, hence $D$ is not a component of $E$, a contradiction.
\end{proof}
\noindent\textit{Proof of Lemma \ref{lem: foliation lsm has lmm} continued}. By Claim \ref{claim: wglc to lmm E exceptional}, $F+E$ is exceptional over $X'$. By Theorem \ref{thm: chlx23 9.4.1}, we may run a $(K_{\Ff_W}+B_W)$-MMP$/X'$ with scaling of an ample$/X'$ divisor, which terminates with a good minimal model $(Y,\Ff_Y,B_Y)/X'$ of $(W,\Ff_W,B_W)/X'$ such that $K_{\Ff_Y}+B_Y\sim_{\Rr,X'}0$. In particular, $(Y,\Ff_Y,B_Y)$ is $\Qq$-factorial ACSS, and $a(D,\Ff_W,B_W)<a(D,\Ff_Y,B_Y)$ for any prime divisor $D$ on $W$ that is exceptional$/Y$. By the negativity lemma, 
$$K_{\Ff_Y}+B_Y=q^*(K_{\Ff'}+B').$$ 
The lemma follows.
\end{proof}

\begin{lem}\label{lem: g-pair weak glc imply lmm}
Let $(X,\Ff,B)/U$ be an lc algebraically integrable foliated triple. If $(X,\Ff,B)/U$ has a bs-weak lc model (resp. bs-semi-ample model), then $(X,\Ff,B)/U$ has a log minimal model (resp. good log minimal model).
\end{lem}
\begin{proof}
By Lemma \ref{lem: g-pair version bir12 2.7} we only need to prove the bs-weak lc model case. The lemma follows immediately from Lemmas \ref{lem: g-pair version bir12 2.8} and \ref{lem: foliation lsm has lmm}.
\end{proof}

\begin{lem}\label{lem: same weak glc model under pullback}
Let $(X,\Ff,B)/U$ and $(Y,\Ff_Y,B_Y)/U$ be two lc algebraically integrable foliated triples, and let $f: Y\rightarrow X$ be a birational morphism such that
$$K_{\Ff_Y}+B_Y=f^*(K_\Ff+B)+E$$
for some $E\geq 0$ that is exceptional$/X$ and $f_*\Ff_Y=\Ff$. Then:
\begin{enumerate}
    \item Any bs-weak lc model of $(X,\Ff,B)/U$ is a bs-weak lc model of $(Y,\Ff_Y,B_Y)/U$.
    \item If $(X,\Ff,B)/U$ has a bs-weak lc model (resp. bs-semi-ample model), then $(Y,\Ff_Y,B_Y)/U$ has a log minimal model (resp. good log minimal model).
\end{enumerate}
\end{lem}
\begin{proof}
(1) Let $(X',\Ff',B')/U$ be a bs-weak lc model of $(X,\Ff,B)/U$, $\phi: X\dashrightarrow X'$ the induced birational map, and $\phi_Y:=\phi\circ f$. Let $p: W\rightarrow Y$ and $q: W\rightarrow X'$ be a resolution of indeterminacy, and let $h:=f\circ p$. By Lemma \ref{lem: hl23 2.6},
$$h^*(K_\Ff+B)=q^*(K_{\Ff'}+B')+F$$
for some $F\geq 0$ that is exceptional over $X'$. Thus 
$$p^*(K_{\Ff_Y}+B_Y)=q^*(K_{\Ff'}+B')+p^*E+F.$$
Thus $a(D,\Ff_Y,B_Y)\leq a(D,\Ff',B')$ for any prime divisor $D$ over $X'$. In particular, if $a(D,\Ff',B')=-\epsilon_{\Ff}(D)$, then $a(D,\Ff_Y,B_Y)=-\epsilon_{\Ff}(D)$.

Since $(X',\Ff',B')/U$ is a log birational model of $(X,\Ff,B)/U$ and $(X,\Ff,B)$ is lc, $$B'=\phi_*B+\Exc(\phi^{-1})^{\Ff'-\mathrm{ninv}}.$$ 
Let 
$$\tilde B':=(\phi_Y)_*B_Y+\Exc(\phi_Y^{-1})^{\Ff'-\mathrm{ninv}}.$$ For any prime divisor $D$ on $X'$, there are two cases:

\medskip

\noindent\textbf{Case 1}. $D$ is not exceptional over $X$. In this case,
$$\mult_DB'=a(D,\Ff',B')=a(D,\Ff,B)=a(D,\Ff_Y,B_Y)=a(D,\Ff',\tilde B')=-\mult_D\tilde B',$$
so $\mult_DB'=\mult_D\tilde B'$.

\medskip

\noindent\textbf{Case 2}. $D$ is exceptional over $X$. In this case, 
$$a(D,\Ff',B')=-\mult_DB'=-\epsilon_{\Ff}(D).$$
Since $a(D,\Ff_Y,B_Y)\leq a(D,\Ff',B')$ and $(Y,\Ff_Y,B_Y)$ is lc, $a(D,\Ff_Y,B_Y)=-\epsilon_{\Ff}(D)$. Therefore, if $D$ is not exceptional over $Y$, then
$$\mult_D\tilde B'=\mult_DB_Y=-a(D,\Ff_Y,B_Y)=\epsilon_{\Ff}(D)=\mult_DB',$$
and if $D$ is exceptional over $Y$, then
$$\mult_D\tilde B'=\mult_D\Exc(\phi_Y^{-1})^{\Ff'-\mathrm{ninv}}=\epsilon_{\Ff}(D)=\mult_DB'.$$
Thus $B'=B''$, hence $(X',\Ff',B')/U$ is a log birational model of $(Y,\Ff_Y,B_Y)/U$. Since $K_{\Ff'}+B'$ is nef$/U$, and
$a(D,\Ff_Y,B_Y)\leq a(D,\Ff',B')$ for any prime divisor $D$ over $X'$, $(X',\Ff',B')/U$ is a bs-weak lc model of $(Y,\Ff_Y,B_Y)/U$, and we get (1).

(2) follows from (1), Lemma \ref{lem: g-pair weak glc imply lmm}, and Lemma \ref{lem: g-pair version bir12 2.7}.
\end{proof}

\subsection{Minimal models and core models}

In this subsection, we shall use core models to study how the (bs-)minimal models of $(X,\Ff,B)/U$ are associated with the (bs-)minimal models of $(X,B+G)/U$ when $(X,\Ff,B
;G)$ satisfies Property $(*)$. First we recall the following results in \cite{CHLX23} and \cite{HH20}:

\begin{lem}[{cf. \cite[Lemma 9.2.1]{CHLX23}}]\label{lem: chlx 9.2.1}
    Let $(X,\Ff,B)/U$ be an lc algebraically integrable foliated triple, $G$ a reduced divisor on $X$, and $f: X\rightarrow Z$ a contraction, such that $(X,\Ff,B;G)/Z$ satisfies Property $(*)$ and $K_{\Ff}+B\sim_{\mathbb R,U}K_X+B+G$. Assume that $G$ is super$/Z$. Let $D\geq 0$ be an $\Rr$-divisor on $X$ such that $K_{\Ff}+B+D$ is nef$/U$.
    
    Then any sequence of steps of a $(K_{\Ff}+B)$-MMP$/U$ (with scaling of $D$) is a sequence of steps of a $(K_{X}+B+G)$-MMP$/U$ (with scaling of $D$), and any sequence of steps of a $(K_{X}+B+G)$-MMP$/U$ (with scaling of $D$) is a sequence of steps of a  $(K_{\Ff}+B)$-MMP$/U$ (with scaling of $D$). Moreover, any sequence of steps of a $(K_{\Ff}+B)$-MMP$/U$ or a 
 $(K_{X}+B+G)$-MMP$/U$ is a sequence of steps of an MMP$/Z$.
\end{lem}

\begin{thm}[{\cite[Theorem 1.7]{HH20}}]\label{thm: hh20 1.7}
    Let $(X,B)/U$ be an lc pair and $A$ an ample$/U$ $\Rr$-divisor on $X$ such that $(X,B+A)$ is lc and $K_X+B+A$ is nef$/U$. Assume that $(X,B)/U$ has a $\Qq$-factorial bs-minimal model or $K_X+B$ is not pseudo-effective$/U$. Then there exists a sequence of $(K_X+B)$-MMP$/U$ with scaling of $A$ which terminates with either a minimal model or a Mori fiber space of $(X,B)/U$. 
\end{thm}

In Lemma \ref{lem: chlx 9.2.1}, and many results in \cite{CHLX23}, we will come up with ``MMP$/U$ is always an MMP$/Z$". If we use the language of core models, then it is essentially saying that ``MMP$/U$ is always an MMP$/Z_U$, where $Z_U$ is the core model of $(X\rightarrow U,X\rightarrow Z)$. We have the following lemmas on showing this fact:

\begin{lem}\label{lem: minimal model foliation over zu}
    Let $(X,\Ff,B)/U$ be an lc algebraically integrable foliated triple. Assume that the associated morphism $\pi: X\rightarrow U$ is a contraction, and assume that $\Ff$ is induced by a contraction $f: X\rightarrow Z$. Let $Z_U$ be the core model of $(\pi,f)$. Then:
    \begin{enumerate}
        \item Any sequence of steps of a $(K_{\Ff}+B)$-MMP$/U$ is a step of a $(K_{\Ff}+B)$-MMP$/Z_U$.
        \item If $(X,\Ff,B)$ is $\Qq$-factorial ACSS and $K_{\Ff}+B$ is nef$/U$, then $K_{\Ff}+B$ is nef$/Z_U$.
        \item $(X,\Ff,B)/U$ has a bs-weak lc model if and only if $(X,\Ff,B)/Z_U$ has a bs-weak lc model.
    \end{enumerate}
\end{lem}
\begin{proof}
    (1) By the universal property of the core model (Definition-Lemma \ref{deflem: relative ample model}), we only need to show that any contraction of a $(K_{\Ff}+B)$-negative extremal ray$/U$ is a contraction$/Z$. This follows from the (relative) cone theorem of algebraically integrable foliations \cite[Theorem 3.9]{ACSS21}, \cite[Theorem 2.3.1]{CHLX23}.

    (2)  If $K_{\Ff}+B$ is not nef$/U$, then there exists a $(K_{\Ff}+B)$-negative extremal ray$/U$ $R$. By the (relative) cone theorem of algebraically integrable foliations (\cite[Theorem 3.9]{ACSS21}, \cite[Theorem 2.3.1]{CHLX23}), $R$ is a a $(K_{\Ff}+B)$-negative extremal ray$/Z$. Since $(X,\Ff,B)$ is a $\Qq$-factorial ACSS, there exists a contraction $\cont_R$ of $R$. By Definition-Lemma \ref{deflem: relative ample model}, $\cont_R$ is a contraction$/Z_U$, which is not possible as $K_{\Ff}+B$ is nef$/Z_U$. Therefore, $K_{\Ff}+B$ is nef$/U$.

    (3) First we suppose that $(X,\Ff,B)/U$ has a bs-weak lc model $(X',\Ff',B')/U$. Let $(W,\Ff_W,B_W)$ be a foliated log smooth model of $(X,\Ff,B)/U$ such that the induced map $W\dashrightarrow X'$ is a morphism. By Lemma \ref{lem: foliation lsm has lmm}, we may run a $(K_{\Ff_W}+B_W)$-MMP$/X'$ which terminates with a log minimal model $(Y,\Ff_Y,B_Y)/X'$ such that $(Y,\Ff_Y,B_Y)/U$ is a log minimal model of $(W,\Ff_W,B_W)/U$. By (1), the induced birational map $Y\dashrightarrow Z_U$ is a morphism, so $(Y,\Ff_Y,B_Y)/Z_U$ is a log minimal model of $(W,\Ff_W,B_W)/Z_U$. By Lemma \ref{lem: g-pair version bir12 2.8}, $(Y,\Ff_Y,B_Y)/Z_U$ is a log minimal model of $(X,\Ff,B)/Z_U$. This proves the only if part.

    Next we prove the if part. Assume that $(X,\Ff,B)/Z_U$ has a bs-weak lc model $(X',\Ff',B')/Z_U$. By Lemma \ref{lem: g-pair weak glc imply lmm}, we may assume that $(X',\Ff',B')/Z_U$ is a log minimal model of $(X,\Ff,B)/Z_U$. By Definition-Lemma \ref{deflem: relative ample model}, $Z_U$ is the core model of $(X'\rightarrow U,X'\rightarrow Z)$. By (2), $K_{\Ff'}+B'$ is nef$/U$, so $(X',\Ff',B')/U$ is a bs-weak lc model of $(X,\Ff,B)/U$. This proves the if part.
\end{proof}

\begin{lem}\label{lem: super minimal model over zu}
    Let $(X,B)/U$ be a pair associated with contraction $\pi: X\rightarrow U$. Let $f: X\rightarrow Z$ be a contraction such that $B$ is super$/Z$. Let $Z_U$ be the core model of $(\pi,f)$. Then:
    \begin{enumerate}
    \item If $K_X+B$ is nef$/Z_U$ then $K_X+B$ is nef$/U$.
    \item Any sequence of steps of a $(K_X+B)$-MMP$/U$ is a sequence of steps of a $(K_X+B)$-MMP$/Z_U$.
        \item $(X,B)/U$ has a minimal model if and only if $(X,B)/Z_U$ has a minimal model.
    \end{enumerate}
\end{lem}
\begin{proof}
Let $d:=\dim X$.

(1) Let $R$ be a $(K_X+B)$-negative extremal ray$/U$. Then $R$ which is spanned by a rational curve $C$ such that $0<-(K_{X}+B)\cdot C\leq 2d$. We may assume that $C$ is of minimal degree among all rational curves which span $R$, i.e. for any rational curve $C'$ such that $[C']=R$, $-(K_X+B)\cdot C'\geq -(K_X+B)\cdot C$.

Since $B$ is super$/Z$, $B=\sum_{i=1}^{2d+1}f^*H_i+B_0$ where $H_i$ are ample Cartier divisors on $Z$ and $B_0\geq 0$. If $f(C)$ is not a point, then
$$(K_X+B_0)\cdot C=(K_X+B)\cdot C-\sum_{i=1}^{2d+1}(f^*H_i\cdot C)<-2d,$$
which contradicts the cone theorem. Therefore, $f(C)$ is a point. The contraction of $C$ exists by the usual cone theorem, and it is a contraction$/Z$ and a contraction$/U$. By the universal property of the core models, the contraction of $C$ is a contraction$/Z_U$. 

Therefore, any contraction of a $(K_X+B)$-negative extremal ray$/U$ is a contraction$/Z_U$, so any step of a $(K_X+B)$-MMP$/U$ $(X,B)\dashrightarrow (Y,B_Y)$ is an MMP$/Z_U$. Since $B$ is super$/Z$, $B_Y$ is super$/Z$. By Definition-Lemma \ref{deflem: relative ample model}, $Z_U$ is the core model of $(Y\rightarrow U,Y\rightarrow Z)$. We may replace $(X,B)$ with $(Y,B_Y)$ and continue this process.

(2) Suppose that $(X',B')/Z_U$ is a minimal model of $(X',B')/U$. Since the induced birational map $X\dashrightarrow X'$ does not extract any divisor and is over $Z$, $B'$ is super$/Z$. If $K_{X'}+B'$ is not nef$/U$, then there exists a step of a $(K_{X'}+B')$-MMP$/U$. This step cannot be over $Z_U$ since $K_{X'}+B'$ is nef$/Z_U$. This contradicts (1), so $(X',B')/U$ is a minimal model of $(X,B)/U$.

Suppose that $(X,B)/U$ has a minimal model. By Lemma \ref{lem: g-pair weak glc imply lmm}, $(X,B)/U$ has a log minimal model. By Theorem \ref{thm: hh20 1.7}, we may run a $(K_X+B)$-MMP$/U$ with scaling of an ample divisor which terminates with a minimal model $(X',B')/U$ of $(X,B)/U$. By (1) the induced map $X'\dashrightarrow Z_U$ is a contraction. Therefore, $(X',B')/Z_U$ is a minimal model of $(X,B)/U$.
\end{proof}

The following proposition is crucial for us to prove Theorem \ref{thm: foliation emm equal to tof scaling}.

\begin{prop}\label{prop: eolmm foliation to pair}
    Let $(X,\Ff,B)/U$ be an lc algebraically integrable foliated triple. Assume that $(X,\Ff,B)/U$ has a bs-weak lc model. Then there exists an ACSS modification $h: (X',\Ff',B';G)/Z\rightarrow (X,\Ff,B)$ that is $\Qq$-factorial, strict, and super, and $(X',B'+G)/U$ has a log minimal model.
\end{prop}
\begin{proof}

    Let $(Y,\Ff_Y,B_Y)/U$ be a bs-weak lc model of $(X,\Ff,B)/U$. Let $g: W\rightarrow X$ be a foliated log resolution of $(X,\Ff,B)$ associated with the equidimensional contraction $f_W: W\rightarrow Z$, such that the induced birational map $W\dashrightarrow Y$ is a morphism, $\Ff_W:=g^{-1}\Ff$ is induced by $f_W$, and $B_W:=g^{-1}_*B+\Exc(g)^{\Ff_W-\mathrm{ninv}}$. Then there exists a reduced divisor $G_W\geq 0$ on $W$ such that $G_W$ is super$/Z$, $\Exc(g)\subset\Supp G_W$, and $(W,\Ff_W,B_W;G_W)/Z$ is ACSS. Moreover, $(W,\Ff_W,B_W)$ is foliated log smooth. 

    Let $\pi: X\rightarrow U$ be the associated morphism and $X\rightarrow U'\rightarrow U$ be the Stein factorization of $\pi$. Possibly replacing $U$ with $U'$, we may assume that $\pi$ is a contraction. Let $Z_U$ be the core model of $(\pi\circ g,f_W)$. By Lemma \ref{lem: foliation lsm has lmm}, we may run a $(K_{\Ff_W}+B_W)$-MMP$/U$ which terminates with a log minimal model $(X'',\Ff'',B'')/U$ of $(W,\Ff_W,B_W)/U$. By Lemma \ref{lem: minimal model foliation over zu}, this MMP is a  $(K_{\Ff_W}+B_W)$-MMP$/Z_U$, so $(X'',\Ff'',B'')/Z_U$ is a log minimal model of $(W,\Ff_W,B_W)/Z_U$. By Lemma \ref{lem: chlx 9.2.1}, $K_{\Ff_W}+B_W\sim_{\mathbb R,Z}K_W+B_W+G_W$, so by Lemma \ref{lem: numerical equivalence model}, $(X'',B''+G'')/Z_U$ is a minimal model of $(W,B_W+G_W)/Z_U$, where $G''$ is the image of $G_W$ on $X''$. By Lemma \ref{lem: super minimal model over zu}, $(W,B_W+G_W)/U$ has a minimal model. 

By Theorem \ref{thm: chlx23 9.4.1}, we may run a $(K_{\Ff_W}+B_W)$-MMP$/X$ with scaling of an ample$/X$ divisor which terminates with a log minimal model $(X',\Ff',B')/X$ such that $K_{\Ff'}+B'\sim_{\mathbb R,X}0$. Let $\phi: W\dashrightarrow X'$ be the induced birational map and $h: X'\rightarrow X$ the induced birational morphism, and let $G:=\phi_*G_W$. By our construction,  $h: (X',\Ff',B';G)/Z\rightarrow (X,\Ff,B)$ is an ACSS modification that is $\Qq$-factorial, strict, and super. By Lemma \ref{lem: chlx 9.2.1}, $\phi$ is also a $(K_X+B+G)$-MMP$/X$.

Let $p: V\rightarrow W$ and $q: V\rightarrow X$ be a resolution of indeterminacy such that $p$ is a log resolution of $(W,B_W+G_W)$ and $q$ is a log resolution of $(X',B'+G)$. Let $\Delta_V:=p^{-1}_*(B_W+G_W)+\Exc(p)$. Then $(V,\Delta_V)$ is a (foliated) log smooth model of $(W,B_W+G_W)$ and $(X',B'+G)$. By Lemma \ref{lem: same weak glc model under pullback}, $(V,\Delta_V)/U$ has a bs-weak lc model. By Lemma \ref{lem: g-pair version bir12 2.8}, $(X',B+G)/U$ has a bs-weak lc model. By Lemma \ref{lem: g-pair weak glc imply lmm}, $(X',B+G)/U$ has a log minimal model.
\end{proof}

\section{Existence of polarized log minimal models}\label{sec: eolmm}

The goal of this section is to prove Theorem \ref{thm: take strict simple model run mmp}, which essentially implies Theorem \ref{thm: eolmm+A} and is crucial for the proofs of Theorems \ref{thm: main} and \ref{thm: mmp can run}. We first recall the following results on the MMP for usual pairs.

\begin{lem}[{cf. \cite[Lemma 2.20]{TX24}}]\label{lem: TX24 2.20}
    Let $(X,B+A)/U$ be an lc pair such that $(X,B)$ is lc and $K_X+B+A$ is nef$/U$. Then there exists a positive real number $\epsilon\in (0,1)$ such that any $(K_X+B+(1-\epsilon)A)$-MMP$/U$ is $(K_X+B+A)$-trivial for any $\epsilon\in (0,\epsilon_0)$.
\end{lem}

\begin{thm}\label{thm: bir12 1.9(3)}
     Let $(X,B)/U$ be an lc pair and $H\geq 0$ an $\Rr$-divisor on $X$ such that $K_X+B+H$ is nef$/U$ and $(X,B+H)$ is lc. Assume that there exists an infinite sequence of $(K_X+B)$-MMP$/U$ with scaling of $H$ with scaling numbers $\lambda_i$ such that $\lim_{i\rightarrow+\infty}\lambda_i=\lambda$ and $\lambda\not=\lambda_i$ for any $i$. Then $(X,B+\lambda H)/U$ does not have a bs-minimal model.
\end{thm}
\begin{proof}
By \cite[Theorem 1.9(3)]{Bir12}, $(X,B+\lambda H)/U$ does not have a bs-minimal model that is $\Qq$-factorial dlt. By \cite[Corollary 3.7]{Bir12}, $(X,B+\lambda H)/U$ does not have a bs-minimal model.
\end{proof}

\begin{thm}[{\cite[Corollary 1.4.2]{BCHM10}}]\label{thm: bchm 1.4.2}
    Let $(X,B)/U$ be a $\Qq$-factorial pair and $A\geq 0$ an $\Rr$-divisor on $X$ such that $B$ is big$/U$, $(X,B+A)$ is klt, and $K_X+B+A$ is nef$/U$. Then any $(K_X+B)$-MMP$/U$ with scaling of $A$ terminates with either a minimal model or a Mori fiber space of $(X,B)/U$.
\end{thm}

\begin{lem}\label{lem: gmmp scaling numbers go to 0}
Let $(X,B)/U$ be an lc pair. Let $H\geq 0$ be an $\Rr$-divisor on $X$ such that $(X,B+H)$ is lc and $K_X+B+H$ is nef$/U$. Assume that for any $\mu\in [0,1]$,
\begin{itemize}
    \item either $(X,B+\mu H)/U$ has a log minimal model, or
    \item $K_X+B+\mu H$ is not pseudo-effective$/U$.
\end{itemize}
Then there exists a $(K_X+B)$-MMP$/U$ with scaling of $H$ which terminates after finitely many steps.
\end{lem}
\begin{proof}
Denote by this MMP
$$(X_1,B_1):=(X,B)\dashrightarrow (X_2,B_2)\dashrightarrow\dots\dashrightarrow (X_i,B_i)\dashrightarrow\cdots.$$
Let $H_i$ be the image of $H$ on $X_i$ for each $i$, and let
$$\lambda_i:=\inf\{t\mid t\geq 0, K_{X_i}+B_i+tH_i\text{ is nef/}U\}$$
be the $i$-th scaling number of this MMP for each $i$.

If $\lambda_1=0$ then there is nothing left to prove. So we may assume that $\lambda_1>0$. By Lemma \ref{lem: TX24 2.20}, we may pick $\lambda_1'\in (0,\lambda_1)$ such that any sequence of a $(K_X+B+\lambda_1'H)$-MMP$/U$ is $(K_X+B+\lambda_1H)$-trivial. 

By Theorem \ref{thm: hh20 1.7}, we may run a $(K_X+B+\lambda_1'H)$-MMP$/U$ with scaling of a general ample$/U$ divisor $A$ which terminates. We let
$$(X_1,B_1):=(X,B)\dashrightarrow (X_2,B_2)\dashrightarrow\dots\dashrightarrow (X_{k_1},B_{k_1})$$
be this sequence of the MMP$/U$. Then this sequence consists of finitely many steps of a $(K_X+B)$-MMP$/U$ with scaling of $H$, with scaling numbers $\lambda_1=\lambda_2=\dots=\lambda_{k_1-1}$. If  $K_X+B+\lambda_1'H$ is not pseudo-effective$/U$, then we have already achieved a $(K_{X_{k_1}}+B_{k_1})$-Mori fiber space$/U$ and we are done. Otherwise,
$$K_{X_{k_1}}+B_{k_1}+\lambda'_1H_{k_1}$$
is nef$/U$, so we have $\lambda_{k_1}\leq\lambda_1'<\lambda_1$. 

We may replace $(X,B)/U$ with $(X_{k_1},B_{k_1})/U$ and continue this process. If this MMP does not terminate, then we may let $\lambda:=\lim_{i\rightarrow+\infty}\lambda_i$. Then $\lambda\not=\lambda_i$ for any $i$, and $K_{X_i}+B_i+\lambda_iH_i$ is nef$/U$. Thus $K_X+B+\lambda H$ is pseudo-effective$/U$. By Theorem \ref{thm: bir12 1.9(3)}, $(X,B+\lambda H)$ does not have a log minimal model, which contradicts our assumption. Therefore, this MMP terminates and we are done.
\end{proof}

\begin{thm}[{\cite[Theorem 1.5]{HH20}}]\label{thm: hh20 1.5}
    Let $(X,B)/U$ be an lc pair and $A$ an ample$/U$ $\Rr$-divisor on $X$ such that $(X,B+A)$ is lc. Then $(X,B+A)/U$ has a bs-good minimal model or a bs-Mori fiber space.
\end{thm}

Combining the above theorem with \cite[Theorem 1.7]{HH20}, we can conclude that the existence of a bs-good minimal model or a bs-Mori fiber space is equivalent to the existence of a good minimal model or a Mori fiber space in the setting of the above theorem.

The following theorem is crucial for the proof of our main theorems.

\begin{thm}\label{thm: take strict simple model run mmp}
    Let $(X,\Ff,B)/U$ be an lc algebraically integrable foliated triple and let $A,H$ be two ample$/U$ $\Rr$-divisors on $X$. Let $h: (X',\Ff',B';G)/Z\rightarrow (X,\Ff,B)$ be a simple model of $(X,\Ff,B)$ that is strict and super, $H':=h^*H$, and $A':=h^*A$. Then:
    \begin{enumerate}
        \item We may run a $(K_{\Ff'}+B'+H')$-MMP$/U$ with scaling of $A'$, which terminates with either a minimal model or a Mori fiber space of $(X',\Ff',B'+H')/U$.
        \item If $X$ is potentially klt, any $(K_{\Ff'}+B'+H')$-MMP$/U$ with scaling of $A'$ terminates with either a minimal model or a Mori fiber space of $(X',\Ff',B'+H')/U$..
    \end{enumerate}
\end{thm}
\begin{proof}
Possibly replacing $A$ with a multiple, we may assume that $K_{\Ff}+B+H+A$ is nef$/U$. Let $\pi: X\rightarrow U$ be the induced projective morphism and let $H_U$ be a sufficiently ample $\Rr$-divisor on $U$. Possibly replacing $A$ with $A+\pi^*H_U$ and $H$ with $H+\pi^*H_U$, we may assume that $A$ and $H$ are ample. Possibly replacing $A$ and $H$, we may assume that $A,H$ are general in $|A|_{\mathbb R}$ and $|H|_{\mathbb R}$ respectively. In particular, $(X',B'+H'+G)$ is lc. 
 
Since $G$ is super$/Z$ and $(X',\Ff',B')/U$ is lc, by Lemma \ref{lem: chlx 9.2.1}, any $(K_{\Ff'}+B'+H')$-MMP$/U$ with scaling of $A'$ is a $(K_{X'}+B'+H'+G)$-MMP$/U$ with scaling of $A'$ and is an MMP$/Z$. Then $K_{\Ff'}+B'+H'+A'$ and $K_{X'}+B'+H'+A'+G$ are nef$/U$. Let $d:=\dim X$.

Let $\bar X$ be the core model of $(h,f)$ associated with $(\bar h,\bar f)$. Let $g: X'\rightarrow \bar X$ be the induced birational morphism. By Lemma \ref{lem: existence of core model}(4), there exists a core model $\bar h: (\bar X,\bar\Ff,\bar B;\bar G)/Z\rightarrow (X,\Ff,B)$ that is strict and super. Let $\bar H:=\bar h^*H$ and $\bar A:=\bar h^*A$. By the definition of core models, $\bar H$ and $\bar A$ are ample$/Z$. Since $\bar G$ is super$/Z$, 
$$\bar G\geq\sum_{i=1}^{2\dim X+1}\bar f^*H_i$$
where $H_i$ are ample Cartier divisors on $Z$. Then there exists $0<\epsilon\ll 1$ such that $\epsilon \bar H+\frac{1}{2}\bar f^*H_1$ is ample. Let $\bar L$ be a general element in 
$$\left|\epsilon\bar H+\frac{1}{2}\bar f^*H_1\right|_{\mathbb R},$$ 
$\widehat H:=(1-\epsilon)\bar H$, and $\widehat G:=\bar G-\frac{1}{2}\bar f^*H_1$. Then
$$K_{\bar X}+\bar B+\widehat H+\widehat G+\bar L\sim_{\mathbb R}K_{\bar X}+\bar B+\bar H+\bar G$$
and $(\bar X,\bar B+\widehat H+\widehat G+\bar L)$ is lc.

\medskip

\noindent\textbf{Step 1}. First we prove the theorem when $X$ is potentially klt. By Lemma \ref{lem: existence of klt pair on strict simple models}, $\bar X$ is potentially klt. Since $\bar L$ is ample, by Lemma \ref{lem: gklt is klt}, there exists a klt pair $(\bar X,\bar\Delta)$ such that
$$0\leq\bar\Delta\sim_{\mathbb R,U}\bar B+\widehat H+\widehat G+\frac{1}{2}\bar L.$$
Let $K_{X'}+\tilde\Delta':=g^*(K_{\bar X}+\bar\Delta)$. Then $(X',\tilde\Delta')$ is sub-klt. Let $0<\delta\ll 1$ be a real number. Since $\Supp (G+B')$ contains all $g$-exceptional prime divisors and $(X',B'+H'+G)$ is lc,
$$(X',\widehat\Delta':=\delta\tilde\Delta'+(1-\delta)(B'+H'+G))$$
is klt. Since $g^*\bar L$ is big and nef, there exist ample $\Rr$-divisors $L_n$ and $\Rr$-divisors $E\geq 0$, such that
$$\frac{\delta}{2}g^*\bar L\sim_{\mathbb R}L_n+\frac{1}{n}E$$
for any positive integer $n$. Then for any $n\gg 0$, $(X',\widehat\Delta'+\frac{1}{n}E)$ is klt. Since $L_n$ is ample$/U$, there exists a klt pair $(X',\Delta')$ such that
$$0\leq\Delta'\sim_{\mathbb R,U}\widehat\Delta'+L_n+\frac{1}{n}E.$$
for some $n\gg 0$. By our construction and Lemma \ref{lem: existence of core model}(2),
$$\Delta'\sim_{\mathbb R,U}B'+H'+G.$$
Now any $(K_{\Ff'}+B'+H')$-MMP$/U$ with scaling of $A'$ is a $(K_{X'}+B'+H'+G)$-MMP$/U$ with scaling of $A'$, hence a $(K_{X'}+\Delta')$-MMP$/U$ with scaling of $A'$. By Theorem \ref{thm: bchm 1.4.2}, any such MMP terminates.

\medskip

\noindent\textbf{Step 2}. Now we prove the general case. For any real number $\mu\in [0,1]$ such that $$K_{X'}+B'+H'+G+\mu A'$$ is pseudo-effective$/U$, by Lemma \ref{lem: existence of core model}(2), $$K_{\bar X}+\bar B+\bar H+\bar G+\mu\bar A$$ 
is pseudo-effective$/U$. Therefore, 
$$K_{\bar X}+\bar B+\widehat H+\widehat G+(\bar L+\mu\bar A)$$ 
is pseudo-effective$/U$. Since $\bar L$ is ample and $\bar A$ is big and nef$/U$, $\bar L+\mu\bar A$ is ample$/U$. Since $H,A,\bar L$ are general,
$$(\bar X,\bar B+\widehat H+\widehat G+(\bar L+\mu\bar A))/U$$
is lc. By Theorem \ref{thm: hh20 1.5} and the remark thereafter,
$$(\bar X,\bar B+\widehat H+\widehat G+(\bar L+\mu\bar A))/U$$
has a good minimal model.

Denote by $f$ the contraction $X'\rightarrow Z$. Since
$$K_{X'}+B'+(1-\epsilon)H'+G-\frac{1}{2}f^*H_1+g^*(\bar L+\mu\bar A)=g^*(K_{\bar X}+\bar B+\widehat H+\widehat G+(\bar L+\mu\bar A)),$$
by Lemmas \ref{lem: g-pair weak glc imply lmm} and \ref{lem: same weak glc model under pullback}, 
$$\left(X',B'+(1-\epsilon)H'+G-\frac{1}{2}f^*H_1+g^*(\bar L+\mu\bar A)\right)\Bigg/U$$
has a good log minimal model. By Theorem \ref{thm: hh20 1.7}, 
$$\left(X',B'+(1-\epsilon)H'+G-\frac{1}{2}f^*H_1+g^*(\bar L+\mu\bar A)\right)\Bigg/U$$
has a minimal model. Since
$$B'+(1-\epsilon)H'+G-\frac{1}{2}f^*H_1+g^*(\bar L+\mu\bar A)\sim_{\mathbb R}B'+H'+G+\mu A',$$
by Lemma \ref{lem: numerical equivalence model}, $(X',B'+H'+G+\mu A')/U$ has a good minimal model. By Lemma \ref{lem: g-pair weak glc imply lmm}, $(X',B'+H'+G+\mu A')/U$ has a log minimal model.

By Lemma \ref{lem: gmmp scaling numbers go to 0}, there exists a $(K_{X'}+B'+H'+G)$-MMP$/U$ with scaling of $A'$ which terminates. By Lemma \ref{lem: chlx 9.2.1}, this MMP$/U$ is also a $(K_{\Ff'}+B'+H')$-MMP$/U$ with scaling of $A'$. The theorem follows.
\end{proof}

\section{A Shokurov-type polytope}\label{sec: sho polytope}

The goal of this section is to prove Theorem \ref{thm: shokurov polytope foliation}. 

\begin{proof}[Proof of Theorem \ref{thm: shokurov polytope foliation}] Let $B(\bm{v}):=\sum_{i=1}^m v_iB_i$ for any $\bm{v}:=(v_1,\dots,v_m)\in\mathbb R^m$. By \cite[Theorem 1.5]{DLM23}, there exists an open subset $U_1\ni\bm{v}_0$ in the rational polytope of $\bm{v}_0$, such that for any  $\bm{v}\in U_1$, $(X,\Ff,B(\bm{v}))$ is lc. We let $c:=\dim U_1$ and let $\bm{v}_1,\dots,\bm{v}_{c+1}$ be vectors in $U_1\cap\mathbb Q^m$ such that $\bm{v}_0$ is contained in the convex hull $U_2$ spanned by
    $\bm{v}_1,\dots,\bm{v}_{c+1}$. Then there exist positive real numbers $a_1,\dots,a_{c+1}$ such that $\sum_{i=1}^{c+1}a_i\bm{v}_i=\bm{v}_0$ and $\sum_{i=1}^{c+1}a_i=1$. We let $I$ be a positive integer such that $I(K_{\Ff}+B(\bm{v}_i))$ is Cartier for each $i$. Let $d:=\dim X$ and $a_0:=\min_{1\leq i\leq c+1}\{a_i\}$. 

    Consider the set
    $$\Ii:=\left\{\sum a_i\gamma_i\mid\gamma_i\in [-2dI,+\infty)\cap\mathbb Z\right\}\cap (0,+\infty).$$
    We have $\gamma_0:=\inf\{\gamma\in\Ii\}>0$. We let $U$ be the interior of the set
    $$\left\{\frac{1}{2d+\gamma_0}(2d\bm{r}+\gamma_0\bm{v})\Bigg| \bm{v}\in U_2\right\}.$$

    We show that $U$ satisfies our requirement. By our construction, $(X,\Ff,B(\bm{v}))$ is lc for any $\bm{v}\in U$ so we only need to show that $K_{\Ff}+\sum_{i=1}^mv_iB_i$ is nef$/Z$ for any $\bm{v}=(v_1,\dots,v_m)\in U$. We let $R$ be an extremal ray in $\overline{NE}(X/U)$. There are three cases.

    \medskip

    \noindent\textbf{Case 1}. $(K_{\Ff}+B)\cdot R=0$. In this case, $(K_{\Ff}+B(\bm{v}))\cdot R=0$ for any $\bm{v}\in U_1$, so $(K_{\Ff}+B(\bm{v}))\cdot R=0$ for any $\bm{v}\in U$.

\medskip

    \noindent\textbf{Case 2}. $(K_{\Ff}+B(\bm{v}_i))\cdot R\geq 0$ for any $i$. In this case, $(K_{\Ff}+B(\bm{v}))\cdot R\geq 0$ for any $\bm{v}\in U_1$, so so $(K_{\Ff}+B(\bm{v}))\cdot R\geq 0$ for any $\bm{v}\in U$.

    \medskip

     \noindent\textbf{Case 3}. $(K_{\Ff}+B)\cdot R>0$ and $(K_{\Ff}+B(\bm{v}_j))\cdot R<0$ for some $j$. In this case, by the relative cone theorem for algebraically integrable foliations (cf. \cite[Theorem 2.2.1]{CHLX23}, \cite[Theorem 3.9]{ACSS21}), $R$ is spanned by a curve $C$ such that $(K_{\Ff}+B(\bm{v}_i)\cdot C\geq -2d$ for any $i$. Thus $$I(K_{\Ff}+B(\bm{v}_i))\cdot C\in [-2dI,+\infty)\cap\mathbb Z,$$
    so
    $$I(K_{\Ff}+B(\bm{v}))\cdot C\in\Ii_0.$$
    Then for any $\bm{v}\in U$, there exists $\bm{v}'\in U_2$ such that $(2d+\gamma_0)\bm{v}=2d\bm{r}+\gamma_0\bm{v}'$. We have
    \begin{align*}
        I(K_{\Ff}+B(\bm{v}))\cdot C&=\frac{\gamma_0}{2d+\gamma_0}I(K_{\Ff}+B(\bm{v}'))\cdot C+\frac{2d}{2d+\gamma_0}I(K_{\Ff}+B(\bm{r}))\cdot C\\
        &\geq \frac{\gamma_0}{2d+\gamma_0}\cdot (-2d)+\frac{\gamma_0}{2d+\gamma_0}\cdot\gamma_0=0,
    \end{align*}
    so $I(K_{\Ff}+B(\bm{v}))\cdot R\geq 0$. The theorem follows.
\end{proof}

\section{Proof of the contraction theorem and the existence of flips}\label{sec: proof flip}

We first prove the contraction theorem when the supporting function is not big, and then prove the contraction and the existence of flips when the supporting function is big.

\begin{prop}\label{prop: contraction not big}
    Let $(X,\Ff,B)/U$ be an lc algebraically integrable foliated triple such that $(X,\Delta)$ is lc for some $B\geq\Delta\geq 0$. Let $R$ be a $(K_{\Ff}+B)$-negative extremal ray$/U$ and $H_R$ a supporting function$/U$ of $R$. Suppose that $H_R$ is not big$/U$. Then $R$ is also a $(K_X+\Delta)$-negative extremal ray$/U$. In particular, there exists a contraction $\cont_R$ of $R$. 
\end{prop}
\begin{proof}
By \cite[Theroem 2.2.1, Lemma 8.4.1]{CHLX23}, we may assume that 
       $$H_R=K_{\Ff}+B+A$$
for some ample$/U$ $\Rr$-Cartier $\Rr$-divisor $A$ on $X$. Let $\pi: X\rightarrow U$ be the induced projective morphism and $X\rightarrow U'\rightarrow U$ the Stein factorization of $\pi$. Possibly replacing $U$ with $U'$, we may assume that $\pi$ is a contraction. 

Let $F$ be a general fiber of $\pi$. Then $H_F:=H_R|_F$ is nef but not big. Let $q:=\dim F$ and $A_F:=A|_F$, then there exists an integer $0\leq k\leq q-1$ such that
$$H_F^k\cdot A_F^{q-k}>H_F^{k+1}\cdot A_F^{q-k-1}=0.$$

Let $D_i:=H_R$ for any $1\leq i\leq k+1$, and let $D_i:=A$ for any $k+2\leq i\leq q$. Then
$$(D_1|_F)\cdot (D_2|_F)\cdots\dots\cdot (D_q|_F)=H_F^{k+1}\cdot A_F^{q-k-1}=0$$
and
$$-(K_{\Ff}+B)|_F\cdot (D_2|_F)\cdots\dots\cdot (D_q|_F)=(A_F-H_F)\cdot H_F^{k}\cdot A_F^{q-k-1}=H_F^{k}\cdot A_F^{q-k}>0.$$
Let $M:=H_R+A=K_{\Ff}+B+2A$. Then $M'$ is nef$/U$. By \cite[Theorem 8.1.1]{CHLX23}, for any general closed point $x\in X$, there exists a rational curve $C_x$ such that $x\in C_x$, $\pi(C_x)$ is a closed point, $C_x$ is tangent to $\Ff$, and 
$$0=D_1\cdot C_x=H_R\cdot C_x.$$
In particular, $C_x$ spans $R$.

By Theorem \ref{thm: eo acss model}, there exists an  ACSS modification $h: (X',\Ff',B';G)/Z\rightarrow (X,\Ff,B)$ that is $\Qq$-factorial and strict. Then $G$ contains any $h$-exceptional $\Ff'$-invariant prime divisor, and $\Supp B'$ contains any $h$-exceptional non-$\Ff'$-invariant prime divisor. In particular, Let $\Delta':=h^{-1}_*\Delta$. Since $(X,\Delta)$ is lc, we may write
    $$K_{X'}+\Delta'+E_+=h^*(K_X+\Delta)+E_-$$
    where $E_+,E_-\geq 0$ are exceptional$/X$, and $E_+\wedge E_-=0$. Then
 $$B'+G\geq \Delta'+E_++\Supp E_-\geq \Delta'+E_+-E_-.$$
Let $x$ be a general closed point in $X$ and let $C_x'$ be the strict transform of $C_x$ on $X'$. Let $A':=h^*A$. Since $x$ is a general closed point in $X$ and $C_x$ is tangent to $\Ff$, $C_x'$ is tangent to $\Ff'$. By Proposition \ref{prop: weak cbf gfq},
\begin{align*}
0&=H_R\cdot C_x=h^*H_R\cdot C_x'=(K_{\Ff'}+B'+A')\cdot C_x'=(K_{X'}+B'+A'+G')\cdot C_x'\\
&\geq (K_{X'}+\Delta'+E_+-E_-+A')\cdot C_x'=h^*(K_X+\Delta+A)\cdot C_x'\\
&=(K_X+\Delta+A)\cdot C_x>(K_X+\Delta)\cdot C_x.
\end{align*}
Therefore, $R$ is a $(K_X+\Delta)$-negative extremal ray. The existence of $\cont_R$ follows from the usual contraction theorem for lc pairs. 
\end{proof}

Finally, we prove the contraction theorem and the existence of flips when the supporting function is big. We remark that generalized foliated quadruples will inevitably be used in the proof of the following theorem. For the convenience of the readers that are not familiar with generalized pairs and/or generalized foliated quadruple, in the following proof, we write footnotes whenever when we have to use generalized foliated quadruples and explain the reasons. We also suggest the readers to consider $\Mm=\Nn=\bm{0}$ throughout the proof.

To prove this theorem, we need to use the concept of generalized foliated quadruples. Nevertheless, we can stick to ``NQC generalized foliated quadruples" as the non-NQC case is harder to prove.

\begin{thm}\label{thm: cont and flip with detail}
    Let $(X,\Ff,B,\Mm)/U$ be an lc algebraically integrable generalized foliated quadruple such that $(X,\Delta,\Nn)/U$ is klt, where $B\geq\Delta\geq 0$ and $\Mm-\Nn$ is nef$/U$. Let $R$ be a $(K_{\Ff}+B+\Mm_X)$-negative extremal ray$/U$ and $A$ an ample$/U$ $\Rr$-divisor on $X$, such that $H_R:=K_\Ff+B+\Mm_X+A$ is a supporting function$/U$ of $R$ and is big$/U$. Then:
    \begin{enumerate}
        \item (Contraction theorem) $H_R$ is semi-ample$/U$. In particular, $H_R$ defines a contraction $\cont_R: X\rightarrow T$. Moreover:
        \begin{enumerate}
            \item If $H_R$ is Cartier, then for any integer $m\gg 0$, $\mathcal{O}_X(mH_R)$ is globally generated over $U$.
            \item For any line bundle $L$ on $X$ such that $L\cdot R=0$, $L\cong\cont_R^*L_T$ for some line bundle $L_T$ on $T$.
        \end{enumerate}
        \item (Existence of flips) The ample model$/T$ $X^+$ of $K_{\Ff}+B+\Mm_X$ exists. Moreover:
        \begin{enumerate}
        \item $(X^+,\phi_*\Delta,\Nn)$ is klt, where $\phi: X\dashrightarrow X^+$ is the induced birational map.
        \item If $\cont_R$ is a small contraction, then the induced morphism $X^+\rightarrow T$ is a $(K_{\Ff}+B+\Mm_X)$-flip$/U$.
        \item If $X$ is $\Qq$-factorial, then:
        \begin{enumerate}
            \item $X^+$ is $\Qq$-factorial.
            \item If $\cont_R$ is a divisorial contraction, then $T=X^+$ and $\rho(X)=\rho(T)+1$. 
            \item If $\cont_R$ is a a small contraction, then $\rho(X)=\rho(X^+)$.
        \end{enumerate}
        \end{enumerate}
    \end{enumerate}
\end{thm}
\begin{proof}
\noindent\textbf{Step 1}. We reduce to the case when $(K_X+\Delta+\Nn_X)\cdot R>0$.

Let $\epsilon\in (0,1)$ be a real number such that $(K_X+\Delta+\Nn_X+\epsilon A)\cdot R\not=0$. By Lemma \ref{lem: gklt is klt}, possibly replacing $\Mm$ with $\Mm+\epsilon\bar A$, $\Nn$ with $\Nn+\epsilon\bar A$\footnote{We remark that this is the first place where we need to use the structure of generalized foliated quadruples. Even if $\Mm=\Nn=\bm{0}$ at the beginning, since it may not be possible for us to get an lc foliated triple $(X,\Ff,B+A)$ even if $A$ is general in $|A|_{\mathbb R/U}$. Nevertheless, if the readers only care about the case when $\Mm=\Nn=\bm{0}$, then the readers may always assume that $\Mm,\Nn$ are NQC$/U$ as this property is preserved throughout the proof.}, and $A$ with $(1-\epsilon)A$, we may assume that $(K_X+\Delta+\Nn_X)\cdot R\not=0$, and there exists a klt pair $(X,\tilde\Delta)$ such that
$$0\leq\tilde\Delta\sim_{\mathbb R,U}\Delta+\Nn_X+A.$$

If $(K_X+\Delta+\Nn_X+A)\cdot R<0$, then $(K_X+\tilde\Delta)\cdot R<0$, and the theorem follows from the  contraction theorem and the existence of flips for klt pairs. Thus we may assume that $(K_X+\Delta+\Nn_X)\cdot R>0$.

\medskip

\noindent\textbf{Step 2}. In this step we construct an ACSS model $(X',\Ff',B',\Mm)$ of $(X,\Ff,B,\Mm)$ and run a sequence of steps of MMP $\phi': X'\dashrightarrow X_n$ for this ACSS model to achieve a model $(X_n,\Ff_n,B_n,\Mm)$.

By Theorem \ref{thm: eo acss model} (actually, we need its generalized pair version \cite[Theorem 2.5.1]{CHLX23}), there exists an ACSS model $h: (X',\Ff',B',\Mm;G)\rightarrow (X,\Ff,B)$ that is $\Qq$-factorial, strict, and super. Then $(X',B'+G,\Mm)$ is $\Qq$-factorial qdlt (cf. Remark \ref{rem: defs for gfq and g-pairs}(4)). Therefore, $(X',B'+G,\Nn)$ is qdlt and and $X'$ is klt. Let 
$$K_{X'}+h^{-1}_*\Delta+\Nn_{X'}+E_1:=h^*(K_X+\Delta+\Nn_X)+E_2$$ for some $E_1\geq 0, E_2\geq 0$ such that $E_1\wedge E_2=0$. Since $(X,\Delta,\Nn)$ is klt, there exists a positive real number $t$ such that all coefficients of 
$$\Delta':=h^{-1}_*\Delta+E_1+t\Exc(h)$$
are strictly less than $1$. Since $B'+G\geq\Exc(h)$ and $h_*(B'+G)\geq B\geq\Delta$, $$B'+G\geq h^{-1}_*\Delta+\Exc(h)\geq\Delta'.$$
Thus $(X',\Delta',\Nn)$ is qdlt. Since $\lfloor\Delta'\rfloor=0$, $(X',\Delta',\Nn)$ is klt. 

Let $A':=h^*A$ and $H_R':=h^*H_R=K_{\Ff'}+B'+\Mm_{X'}+A'$. By the generalized pair version of Theorem \ref{thm: take strict simple model run mmp}(2)\footnote{See Theorem \ref{thm: take strict simple model run mmp gfq}(2).}, we may run a $(K_{\Ff'}+B'+\Mm_{X'}+\frac{1}{2}A')$-MMP$/U$ with scaling of $\frac{1}{2}A'$ which terminates. This MMP is also a sequence of steps of a $(K_{\Ff'}+B'+\Mm_{X'})$-MMP$/U$ with scaling of $A'$. We let
$$(X_0,\Ff_0,B_0,\Mm):=(X',\Ff',B',\Mm)\dashrightarrow (X_1,\Ff_1,B_1,\Mm)\dashrightarrow\dots\dashrightarrow (X_n,\Ff_n,B_n,\Mm)\dashrightarrow\dots$$
be this MMP, let $A_i,H_i$ be the images of $A',H_R'$ on $X_i$ for each $i$, and let
$$\lambda_i:=\inf\{t\geq 0\mid K_{\Ff_i}+B_i+\Mm_{X_i}+\lambda_iA_i'\text{ is nef}/U\}$$
be the scaling numbers. Let $n$ be the smallest index such that $\lambda_i<1$. Then the induced birational map $\phi': X'\dashrightarrow X_n$ is a sequence of steps of a $(K_{\Ff'}+B'+\Mm_{X'}+\lambda_n A')$-MMP$/U$ and is $(K_{\Ff'}+B'+\Mm_{X'}+A')$-trivial for each step.

Let $\Delta_n,G_n$ be the images of $\Delta',G$ on $X_n$ respectively. By the generalized pair version of Lemma \ref{lem: chlx 9.2.1}\footnote{See Lemma \ref{lem: chlx 9.2.1 gfq}.}, $\phi'$ is a sequence of steps of a $(K_{X'}+B'+G+\Mm_{X'})$-MMP$/U$, $(X_n,B_n+G_n,\Mm)$ is $\Qq$-factorial qdlt. Since $B'+G\geq\Delta'$, $B_n+G_n\geq\Delta_n$. Since $\Mm-\Nn$ is nef$/U$, $(X_n,\Delta_n,\Nn)$ is $\Qq$-factorial qdlt. Since $\lfloor\Delta_n\rfloor=0$, $(X_n,\Delta_n,\Nn)$ is klt. 

\medskip

\noindent\textbf{Step 3}. In this step we construct an MMP $\varphi: X_n\dashrightarrow \hat X$.

\begin{claim}\label{claim: second mmp trivialness}
    There exists a positive real number $\delta_0\in (0,\frac{1}{2})$ and a function $\mu: (0,\delta_0)\rightarrow (0,+\infty)$, such that for any $\delta\in (0,\delta_0)$ that is general in $\mathbb R/\mathbb Q$ and any $l>\mu(\delta)$, any sequence of steps of a 
$$((K_{X_n}+\Delta_n+\Nn_{X_n})+l(K_{\Ff_n}+B_n+\Mm_{X_n}+(1-\delta)A_n))\text{-MMP}/U$$
is $(K_{\Ff_n}+B_n+\Mm_{X_n}+(1-\delta)A_n)$-trivial,  $(K_{\Ff_n}+B_n+\Mm_{X_n}+A_n)$-trivial, and $A_n$-trivial.  
\end{claim}
\begin{proof}
By Lemma \ref{lem: general real coefficient b-trivial}, we only need to show that the MMP is  $(K_{\Ff_n}+B_n+\Mm_{X_n}+(1-\delta)A_n)$-trivial. When $\Mm,\Nn$ are NQC$/U$, by the generalized pair version of Theorem \ref{thm: shokurov polytope foliation}\footnote{See Theorem \ref{thm: shokurov polytope foliation gfq}.}, $K_{\Ff_n}+B_n+\Mm_{X_n}+(1-\delta)A_n$ is NQC$/U$ for any $\delta\in (0,1-\lambda_n)$, and the claim follows from \cite[Lemma 4.4(3)]{BZ16}\footnote{We remark if we start with $\Mm=\Nn=\bm{0}$ then the NQC case is enough. We also remark that this is the second place where we need generalized quadruples even if $\Mm=\Nn=\bm{0}$. This is because we need to apply Theorem \ref{thm: shokurov polytope foliation} and Lemma \ref{lem: bz16 4.4(3)} for $K_{\Ff_n}+B_n+\Mm_{X_n}+(1-\delta)A_n$. Even if $\Mm=\bm{0}$, we still need to consider the generalized foliated quadruple structure $(X_n,\Ff_n,B_n,(1-\delta)\bar A)$ as $(X_n,\Ff_n,B_n+(1-\delta)A_n)$ may not be lc.}.

When $\Mm$ and $\Nn$ are not necessarily NQC$/U$, the proof is more complicated and the claim follows from  Proposition \ref{prop: kf+k_x mmp trivial}. We remark that Proposition \ref{prop: kf+k_x mmp trivial} shows that we do not need $\delta$ to be general in $\mathbb R/\mathbb Q$ but we do not need this fact.
\end{proof}

\noindent\textit{Proof of Theorem \ref{thm: cont and flip with detail} continued}.
By the generalized pair version of Lemma \ref{lem: lc+ample nef=nqc}\footnote{See Lemma \ref{lem: lc+ample nef=nqc gfq}. This is the third time when we need to use generalized foliated quadruples. Even if $\Mm=\bm{0}$, we need to consider $(X,\Ff,B,\bar A)$ as it is possible that $(X,\Ff,B+A)$ is not lc.} $K_{\Ff}+B+\Mm_X+A$ is NQC$/U$. By the generalized pair version of Lemma \ref{lem: trivial ray when -delta triple}\footnote{See Lemma \ref{lem: trivial ray when -delta gfq}.}, there exists a real number $\delta_1\in (0,\delta_0)$ satisfying the following:
\begin{itemize}
\item $\delta_1$ is general in $\mathbb R/\mathbb Q$.
\item $\delta_1<1-\lambda_n$.
    \item $K_{\Ff}+B+\Mm_X+(1-\delta_1)A$ is big$/U$. 
    \item $\textbf{B}_-(H_R-2\delta_1 A/U)=\textbf{B}_+(H_R/U)$. 
    \item $R$ is an $(H_R-2\delta_1 A)$-negative extremal ray$/U$, and is the only $(H_R-2\delta_1 A)$-non-positive extremal ray$/U$.
\end{itemize}
Since $\delta_1<1-\lambda_n$, $\phi'$ is a sequence of steps of a $(K_{\Ff'}+B'+\Mm_{X'}+(1-\delta_1)A')$-MMP$/U$. Thus there exists a positive real number $l_1>\mu(\delta_1)$ satisfying the following:
\begin{itemize}
\item $\delta_1 A+\frac{1}{l_1}(K_X+\Delta+\Nn_X)$ and $\delta_1 A-\frac{1}{l_1}(K_X+\Delta+\Nn_X)$ are ample$/U$,
\item $\phi'$ is a sequence of steps of a 
    $$\left((K_{\Ff'}+B'+\Mm_{X'}+(1-\delta_1)A')+\frac{1}{l_1}(K_{X'}+\Delta'+\Nn_{X'})\right)\text{-MMP}/U,$$
    hence a sequence of steps of a
    $$((K_{X'}+\Delta'+\Nn_{X'})+l_1(K_{\Ff'}+B'+\Mm_{X'}+(1-\delta_1)A'))\text{-MMP}/U.$$
    \item $l_1$ is general in $\mathbb R/\mathbb Q$.
\end{itemize}
Let 
$$\Pp:=\overline{K_{\Ff_n}+B_n+\Mm_{X_n}+(1-\delta_1)A_n}=\overline{H_n-\delta_1A_n}.$$ Since $(X_n,\Delta_n,\Nn)$ is klt, $(X_n,\Delta_n,\Nn+l_1\Pp)$ is klt. By our choice of $\delta_1$ and $l$,
$$K_{X_n}+\Delta_n+\Nn_{X_n}+l_1\Pp_{X_n}$$
is big. By \cite[Lemma 4.4(2)]{BZ16}, we may run a 
$$(K_{X_n}+\Delta_n+\Nn_{X_n}+l_1\Pp_{X_n})\text{-MMP}/U$$
with scaling of an ample$/U$ $\Rr$-divisor which terminates with a good minimal model 
$$(\widehat X,\widehat\Delta,\Nn+l_1\Pp)/U$$
of $(X_n,\Delta_n,\Nn+l_1\Pp)/U$, and $(\widehat X,\widehat\Delta,\Nn+l_1\Pp)$ is klt\footnote{This is the fourth time when a generalized pair or a generalized foliated quadruple structure is used.}. We let $\varphi: X_n\dashrightarrow\widehat X$ be the induced birational map. By our construction, $\varphi$ is $\Pp_{X_n}$-trivial.

\medskip

\noindent\textbf{Step 4}. In this step we show that the induced birational map $X\dashrightarrow\widehat X$ does not extract any divisor and is $H_R$-trivial.

Let $\psi:=X'\dashrightarrow\widehat X$ and $\alpha: X\dashrightarrow\widehat X$ be the induced birational maps. By our construction, $\psi$ is a sequence of steps of a 
$$((K_{X'}+\Delta'+\Nn_{X'}+l_1(K_{\Ff'}+B'+\Mm_{X'}+(1-\delta_1)A'))\text{-MMP}/U.$$
Let $$N:=N_{\sigma}(X'/U,K_{X'}+\Delta'+\Nn_{X'}+l_1(H_R'-\delta_1 A')).$$ Since
$$K_{X'}+\Delta'+\Nn_{X'}+l_1(H_R'-\delta_1 A')=h^*(l_1(H_R-\delta_1 A)+(K_X+\Delta+\Nn_X))+(E_2+t\Exc(h)),$$
By Lemma \ref{lem: nz keep under pullback}(1), $\Exc(h)\subset\Supp N$. By Lemma \ref{lem: nz keep under pullback}(2)(3),
\begin{align*}
    \Supp f_*N&=\Supp N_{\sigma}(X/U,l_1(H_R-\delta_1 A)+(K_X+\Delta+\Nn_X))\\
    &=\Bb_-(X/U,l_1(H_R-\delta_1 A)+(K_X+\Delta+\Nn_X)).
\end{align*}
Since $l_1\delta A-(K_X+\Delta+\Nn_X)$ and  $l_1\delta_1 A+(K_X+\Delta+\Nn_X)$ are ample$/U$,
\begin{align*}
    &\Bb_-(X/U,l_1(H_R-\delta A)+(K_X+\Delta+\Nn_X))\\=&\Bb_-(X/U,l_1H_R-(l_1\delta_1 A-(K_X+\Delta+\Nn_X)))\\
    \subset&\Bb_+(X/U,H_R),
\end{align*}
and
\begin{align*}
   &\Bb_-(X/U,l_1(H_R-\delta A)+(K_X+\Delta+\Nn_X))\\
   =&\Bb_-(X/U,l_1(H_R-2\delta_1A)+(l_1\delta_1 A+(K_X+\Delta+\Nn_X)))\\
   \supset&\Bb_-(X/U,H_R-2\delta_1A)=\Bb_+(X/U,H_R).
\end{align*}
Therefore, $\Supp f_*N=\Bb_+(X/U,H_R)$, so 
$$\Supp N=\Exc(h)\cup h^{-1}_*\Supp\Bb_+(X/U,H_R).$$
By Lemma \ref{lem: nz for lc divisor}, $\Exc(\psi)=\Supp N$. In particular, the induced birational map $\alpha: X\dashrightarrow \widehat X$ does not extract any divisor. Since $\psi$ and $h$ are $H_R'$-trivial, $\alpha$ is $H_R$-trivial. 

\medskip

\noindent\textbf{Step 5}. In this step we construct the contraction $\cont_R: X\rightarrow T$ and prove (1).

Recall that by \textbf{Step 1}, we have $(K_X+\Delta+\Nn_X)\cdot R>0$. Thus there exists a positive real number $c$ such that
$$cA\cdot R=(K_X+\Delta+\Nn_X)\cdot R.$$
Let $$L_R:=l_1H_R-(K_X+\Delta+\Nn_X)+cA.$$ Then $L_R\cdot R=0$. Since
$$L_R=l_1(H_R-\delta_1A)+l_1\left(\delta_1A-\frac{1}{l_1}(K_X+\Delta+\Nn_X)\right)+cA,$$
$R$ is the only $(H_R-\delta_1A)$-negative extremal ray$/U$, and $\delta_1A+\frac{1}{l_1}(K_X+\Delta+\Nn_X)$ and $A$ are ample$/U$, we have that $L_R$ is a supporting function$/U$ of $R$.  By Lemma \ref{lem: two supporting function trivial}, $\alpha$ is $L_R$-trivial. 

Let $\widehat A, \widehat H_R$, and $\widehat L_R$ be the images of $A, H_R$ and $L_R$ on $\widehat X$ respectively. Then $\widehat L_R$ is big and nef, and
$$\left(l_1+\frac{c}{\delta_1}\right)\widehat H_R=(K_{\widehat X}+\widehat\Delta+\Nn_{\widehat X})+\widehat L_R+\frac{c}{\delta_1}(\widehat H_R-\delta_1\widehat A)=(K_{\widehat X}+\widehat\Delta+\Nn_{\widehat X})+\widehat L_R+\frac{c}{\delta_1}\Pp_{\widehat X}.$$

Since $\Pp$ descends to $X_n$, $\Pp_{X_n}$ is nef. Since $\varphi$ is $\Pp_{X_n}$-trivial, $\Pp_{\widehat X}$ is nef. Since $(\widehat X,\widehat\Delta,\Nn+l_1\Pp)$ is klt, $(\widehat X,\widehat\Delta,\Nn)$ is klt. Since $L_R$ is big$/U$ and nef$/U$ and $\alpha$ is $L_R$-trivial, $\widehat L_R$ is  big$/U$ and nef$/U$. Thus $\widehat L_R+\frac{c}{\delta_1}\Pp_{\widehat X}$ is big$/U$ and nef$/U$. By the base-point-freeness theorem\footnote{Actually, we need the base-point-freeness theorem for generalized pairs (cf. \cite[Theorem 2.2.6]{CHLX23}) here.}, $\widehat H_R$ is semi-ample$/U$, and $\mathcal{O}_{\widehat X}(mH_R)$ is globally generated$/U$ if $\widehat H_R$ is Cartier. We let $\cont_R: X\rightarrow T$ and $\widehat{\cont_R}: \widehat X\rightarrow T$ be the contractions$/U$ induced by $H_R$ and $\widehat H_R$ respectively. 

If $H_R$ is Cartier, then $H_R'$ is Cartier. Since $\phi'$ is a sequence of steps of an MMP of a $\Qq$-factorial ACSS generalized foliated quadruples and $\varphi$ is a sequence of steps of an MMP of a klt pair, $\widehat H_R$ is Cartier. Therefore, $\mathcal{O}_{\widehat X}(m\widehat H_R)$ is globally generated$/U$ for any integer $m\gg 0$. Since $mH_R=\cont_R^*(\widehat{\cont_R})_*(m\widehat H_R)$, $\mathcal{O}_{X}(mH_R)$ is globally generated$/U$ for any integer $m\gg 0$. This implies (1.a).

We prove (1.b). Since $L-(K_{\Ff}+B+\Mm_X)$ is ample$/T$, by (1.a), $\mathcal{O}_X(mL)$ is globally generated$/T$ for any $m\gg 0$. Thus $mL\cong\cont_R^*L_{T,m}$ and  $(m+1)L\cong\cont_R^*L_{T,m+1}$ for line bundles $L_{T,m}$ and $L_{T,m+1}$ for any $m\gg 0$. We may let $L_T:=L_{T,m+1}-L_{T,m}$.

\medskip

\noindent\textbf{Step 6}. We prove (2) and conclude the proof of the theorem.

Since $K_{\widehat X}+\widehat\Delta+\Nn_{\widehat X}+l_1\Pp_{\widehat X}$ is semi-ample$/U$, it is semi-ample$/T$. Let $\beta: \widehat X\rightarrow X^+$ be the ample model$/T$ of $K_{\widehat X}+\widehat\Delta+\Nn_{\widehat X}+l_1\Pp_{\widehat X}$. Since $\widehat X$ is a minimal model of $$K_{X'}+\Delta'+\Nn_{X'}+l_1(H_R'-\delta_1A'),$$ $X^+$ is the ample model$/T$ of
$$K_{X'}+\Delta'+\Nn_{X'}+l_1(H_R'-\delta_1A'),$$
so the induced birational map $X^+$ is the ample model$/T$ of 
$$(K_X+\Delta+\Nn_X)+l_1(H_R-\delta_1A).$$
Since $R$ is the only $((K_X+\Delta+\Nn_X)+l_1(H_R-\delta_1A))$-negative extremal ray$/U$, $X^+$ is also the ample model$/T$ of $K_{\Ff}+B+\Mm_X$.

Since $\beta$ is $(K_{\widehat X}+\widehat\Delta+\Nn_{\widehat X}+l_1\Pp_{\widehat X})$-trivial and $l_1$ is general in $\mathbb R/\mathbb Q$, by Lemma \ref{lem: general real coefficient b-trivial}, $\beta$ is $(K_{\widehat X}+\widehat\Delta+\Nn_{\widehat X})$-trivial. Since $(\bar X,\bar\Delta,\Nn)$ is klt, $(X^+,\Delta^+,\Nn)$ is klt, where $\Delta^+$ is the image of $\widehat\Delta$ on $X^+$. This implies (2.a). (2.b) follows from the definition of a flip.

The proof of (2.c) for the divisorial contraction case is similar to the proof of \cite[Corollary 3.17]{KM98}, and the proof of (2.c) for the flipping contraction case is similar to \cite[Theorem 6.1, Step 3]{HL23}. We omit these proofs.
\end{proof}

\section{Lifting of the MMP}\label{sec: lifting}

We emphasize that the following lemma is not a direct consequence of the cone theorem in \cite{ACSS21}, since $A_n$ may not be ample$/U$ and Bertini-type theorems fail for foliations.

\begin{lem}\label{lem: can run mmp scaling}
        Let $(X,\Ff,B)/U$ be an lc algebraically integrable foliated triple and $A$ an ample$/U$ $\Rr$-divisor on $X$. Let $\mathcal{P}:$
    $$(X_0,\Ff_0,B_0):=(X,\Ff,B)\dashrightarrow (X_1,\Ff_1,B_1)\dashrightarrow\dots\dashrightarrow (X_n,\Ff_n,B_n)$$
    be a sequence of steps of a $(K_{\Ff}+B)$-MMP$/U$ with scaling of $A$ and let $A_i$ be the image of $A$ on $X_i$ for each $i$. Let
    $$\lambda_n:=\inf\{t\geq 0\mid K_{\Ff_n}+B_n+tA_n\text{ is nef}/U\}.$$
    Suppose that $\lambda_n>0$. Then there exists a $(K_{\Ff_n}+B_n)$-negative extremal ray$/U$ $R$ such that $(K_{\Ff_n}+B_n+\lambda_nA_n)\cdot R=0$.
\end{lem}
\begin{proof}
    $\mathcal{P}$ is also a sequence of steps of a $(K_{\Ff}+B+\lambda_nA)$-MMP$/U$ with scaling of $A$. By \cite[Lemma 16.1.1]{CHLX23}, there exists an lc algebraically integrable generalized foliated quadruple $(X_n,\Ff_n,B_n',\Mm')/U$ and an ample$/U$ $\Rr$-divisor $A_n'$ such that $$K_{\Ff_n}+B_n'+A_n'+\Mm'_{X_n}\sim_{\mathbb R,U}K_{\Ff_n}+B_n+\lambda_nA_n.$$
     By the generalized pair version of Lemma \ref{lem: lc+ample nef=nqc}\footnote{See Lemma \ref{lem: lc+ample nef=nqc gfq}.}, $K_{\Ff_n}+B_n'+A_n'+\Mm'_{X_n}$ is NQC$/U$. Thus there exists a positive real number $\epsilon_0$, such that for any curve $C$ on $X_n$, either $(K_{\Ff_n}+B_n+\lambda_nA_n)\cdot C\geq\epsilon_0$ or $(K_{\Ff_n}+B_n+\lambda_nA_n)\cdot C=0$. Since 
        $$\lambda_n:=\inf\{t\geq 0\mid K_{\Ff_n}+B_n+tA_n\text{ is nef}/U\},$$
        there exists a positive real number $\delta\in\left(0,\frac{\epsilon_0}{2\dim X+\epsilon_0}\right)$ and a 
        $(K_{\Ff_n}+B_n+(1-\delta)\lambda_nA_n)$-negative extremal ray$/U$ $R$. If $(K_{\Ff_n}+B_n+\lambda_nA_n)\cdot R\not=0$, then $R$ is spanned by a curve $C$ such that $$0>(K_{\Ff_n}+B_n)\cdot C\ge -2\dim X$$
    and $$(K_{\Ff_n}+B_n+\lambda_nA_n)\cdot C\geq\epsilon_0,$$ 
    so 
\begin{align*}        
    &\left(K_{\Ff_n}+B_n+\lambda_n\frac{2\dim X}{2\dim X+\epsilon_0}A_n\right)\cdot C\\=&\frac{\epsilon_0}{2\dim X+\epsilon_0}(K_{\Ff_n}+B_n)\cdot C+\frac{2\dim X}{2\dim X+\epsilon_0}(K_{\Ff_n}+B_n+\lambda_nA_n)\cdot C\\\geq&(\frac{\epsilon_0}{2\dim X+\epsilon_0})(-2\dim X)+(\frac{2\dim X}{2\dim X+\epsilon_0})\epsilon_0=0,
\end{align*}
    
    which is not possible. Therefore, $(K_{\Ff_n}+B_n+\lambda_nA_n)\cdot R=0$ and we are done.
\end{proof}

\begin{prop}\label{prop: lift mmp}
    Let $(X,\Ff,B)/U$ be an lc algebraically integrable foliated triple. Let $\mathcal{P}:$
$$(X_0,\Ff_0,B_0):=(X,\Ff,B)\dashrightarrow (X_1,\Ff_1,B_1)\dashrightarrow\dots\dashrightarrow (X_n,\Ff_n,B_n)\dashrightarrow\dots$$
be a (possibly infinite) sequence of $(K_{\Ff}+B)$-MMP$/U$. For each $i\geq 0$, we let $\psi_i: X_i\rightarrow T_i$ and $\psi_i^+:X_{i+1}\rightarrow T_{i}$ be the $(i+1)$-th step of this MMP and let $\phi_i:=(\psi_{i}^+)^{-1}\circ\psi_i: X_i\dashrightarrow X_{i+1}$ be the induced birational map. Let $h: (Y,\Ff_Y,B_Y;G)/Z\rightarrow (X,\Ff,B)$ be an ACSS modification of $(X,\Ff,B)$ that is $\Qq$-factorial, strict, and super. Let $A$ be an ample$/U$ $\Rr$-divisor on $X$ and let $A_i$ be the image of $A$ on $X_i$ for each $i$.

Then there exist a (possibly infinite) sequence $\mathcal{P}_Y$ of birational maps 
$$(Y_0,\Ff_{Y_0},B_{Y_0}):=(Y,\Ff_Y,B_Y)\dashrightarrow (Y_1,\Ff_{Y_1},B_{Y_1})\dashrightarrow\dots\dashrightarrow (Y_n,\Ff_{Y_n},B_{Y_n})\dashrightarrow\dots$$
satisfying the following. Let $\phi_{i,Y}: Y_i\dashrightarrow Y_{i+1}$ be the birational map in the $\mathcal{P}_Y$ above. Then:
\begin{enumerate}
\item For any $i\geq 0$, there exist an ACSS modification $h_i: (Y_i,\Ff_{Y_i},B_{Y_i};G_i)/Z\rightarrow (X_i,\Ff_i,B_i)$ that is $\Qq$-factorial, strict, and super, such that $h_0=h$ and $G_i$ is the image of $G$ on $Y_i$.
\item For any $i\geq 0$, $h_{i+1}\circ\phi_{i,Y}=\phi_i\circ h_i$.
\item For any $i\geq 0$, $\phi_{i,Y}$ is a sequence of steps of $(K_{\Ff_i}+B_{Y_i})$-MMP$/T_i$ and $(Y_{i+1},\Ff_{Y_{i+1}},B_{Y_{i+1}})/T_i$ is the output of this MMP, such that $\phi_{i,Y}$ is not the identity map.
\item $\mathcal{P}_Y$ is a sequence of steps of a $(K_{\Ff_Y}+B_Y)$-MMP$/U$.
\item Assume that $\mathcal{P}$ is a sequence of steps of MMP$/U$ with scaling of $A$. Let $A_{Y}:=h^*A$ and let $A_{Y_i}$ the image of $A_Y$ on $Y_i$ for each $i$. Let
$$\lambda_i:=\inf\{t\geq 0\mid K_{\Ff_i}+B_i+tA_i\text{ is nef}/U\}$$
be the $(i+1)$-th scaling number. Then:
\begin{enumerate}
    \item $\phi_{i,Y}$ is a sequence of steps of a $(K_{\Ff_{Y_i}}+B_{Y_i})$-MMP$/U$ with scaling of $A_{Y_i}$, and the scaling number of each step of $\phi_{i,Y}$ is $\lambda_i$.
    \item $\mathcal{P}_Y$ is sequence of steps of a $(K_{\Ff_Y}+B_Y)$-MMP$/U$ with scaling of $A_Y$.
\end{enumerate}
\end{enumerate}
\end{prop}
\begin{proof}
Since (4) follows from (3) and (5.b) follows from (5.a), we only need to prove (1)(2)(3) and (5.a).

Let $n$ be a non-negative integer. We prove the proposition by induction on $n$ under the induction hypothesis that we have already constructed $(Y_i,\Ff_i,B_i;G_i)/U$ and $h_i$ for any $i\leq n$ and $\phi_{i,Y}$ for any $i\leq n-1$ which satisfy (1)(2)(3)(5). When $n=0$, this follows from our assumption, so we may assume that $n>0$. We need to construct $\phi_{n,Y},h_{n+1}$, and $(Y_{n+1},\Ff_{n+1},B_{n+1};G_{n+1})/U$.

We let $H_n$ be a supporting function of the extremal ray$/U$ contracted by $\psi_n$ and let  
$$L_n:=H_n-(K_{\Ff_n}+B_n),$$
such that $L_n=\lambda_nA_n$ if $\mathcal{P}$ is an MMP$/U$ with scaling of $A$. Then $L_n$ is ample$/T_n$. Now we run a $(K_{\Ff_{Y_n}}+B_{Y_n})$-MMP$/T_n$ with scaling of $h_n^*L_n$, then this MMP is also a $(K_{\Ff_{Y_n}}+B_{Y_n}+\frac{1}{2}h_n^*L_n)$-MMP$/T_i$ with scaling of $h_n^*L_n$. By Theorem \ref{thm: take strict simple model run mmp}, we may choose such an MMP which terminates with a good minimal model $(Y_{n+1},\Ff_{Y_{n+1}},B_{Y_{n+1}})/T_n$ of $(Y_{n},\Ff_{Y_{n}},B_{Y_{n}})/T_n$. Since $X_{n+1}$ is the ample model$/T_n$ of $K_{\Ff_n}+B_n$, $X_{n+1}$ is also the ample model$/T_n$ of $K_{\Ff_{Y_{n+1}}}+B_{Y_{n+1}}$, so there exists an induced birational morphism $h_{n+1}: Y_{n+1}\rightarrow X_{n+1}$. Since $K_{Y_n}+B_{Y_n}$ is not nef$/T_n$ and $K_{Y_{n+1}}+B_{Y_{n+1}}$ is nef$/T_n$, $\phi_{i,Y}$ is not the identity map.

Let $\phi_{n,Y}: Y_n\dashrightarrow Y_{n+1}$ be the induced birational map. (1) for $n+1$ immediately follows by our construction and \cite[Lemma 9.1.4]{CHLX23}. (2)(3) for $n+1$ follow immediately from our construction. Since $H_n\sim_{\mathbb R,T_n}0$, $\phi_{n,Y}$ is $(h_n^*H_n)$-trivial, so (5.a) for $n+1$ immediately follows. Thus (1)(2)(3) and (5.a) follow from induction on $n$ and the proposition follows.
\end{proof}

\section{Proof of the main theorems}\label{sec: proof of main theorem}

\begin{proof}[Proof of Theorem \ref{thm: main}]
     It follows from Proposition \ref{prop: contraction not big} and Theorem \ref{thm: cont and flip with detail}.
\end{proof}

\begin{proof}[Proof of Theorem \ref{thm: mmp can run}]
    By Proposition \ref{prop: contraction not big} and Theorem \ref{thm: cont and flip with detail}, we can run a step of a $(K_{\Ff}+B)$-MMP$/U$. By Theorem \ref{thm: cont and flip with detail}(2.a), after a step of the MMP $\phi: X\dashrightarrow X'$ that is not a Mori fiber space, $(X',\Delta':=\phi_*\Delta)$ is klt. Thus we may continue this process. 
\end{proof}

\begin{proof}[Proof of Theorem \ref{thm: main theorem smooth projective}]
    It is a special case of Theorems \ref{thm: main} and \ref{thm: mmp can run}.
\end{proof}

\begin{thm}\label{thm: can run mmp scaling}
      Let $(X,\Ff,B)/U$ be an lc algebraically integrable foliated triple such that $(X,\Delta)$ is klt for some $B\geq\Delta\geq 0$. Let $A$ be an ample$/U$ $\Rr$-divisor on $X$. Then we may run a $(K_{\Ff}+B)$-MMP$/U$ with scaling of $A$. 
\end{thm}
\begin{proof}
    It follows from Theorem \ref{thm: mmp can run} and Lemma \ref{lem: can run mmp scaling}.
\end{proof}

\begin{proof}[Proof of Theorem \ref{thm: eomfs}]
By Theorem \ref{thm: can run mmp scaling}, we can run a $(K_{\Ff}+B)$-MMP$/U$ with scaling of any ample$/U$ $\Rr$-divisor $A$. Let $\epsilon$ be a positive real number such that $(K_{\Ff}+B+\epsilon A)$ is not pseudo-effective$/U$. 

Let $h: (X',\Ff',B';G)/Z\rightarrow (X,\Ff,B)$ be an ACSS modification of $(X,\Ff,B)$ that is $\Qq$-factorial, strict, and super. By Proposition \ref{prop: lift mmp}, any infinite sequence of steps of a $(K_{\Ff}+B+\epsilon A)$-MMP$/U$ with scaling of $A$ induces an infinite sequence of steps of a $(K_{\Ff'}+B'+\epsilon h^*A)$-MMP$/U$ with scaling of $h^*A$. By Theorem \ref{thm: take strict simple model run mmp}, any $(K_{\Ff'}+B'+\epsilon h^*A)$-MMP$/U$ with scaling of $h^*A$ terminates with a Mori fiber space$/U$. Thus any $(K_{\Ff}+B+\epsilon A)$-MMP$/U$ with scaling of $A$ terminates with a Mori fiber space$/U$, and the theorem follows.    
\end{proof}

\begin{proof}[Proof of Theorem \ref{thm: eolmm+A 1}(1)]
By the generalized pair version of Theorem \ref{thm: can run mmp scaling}\footnote{See Theorem \ref{thm: can run mmp scaling gfq}.}, we can run a $(K_{\Ff}+B+A)$-MMP$/U$ with scaling of any ample$/U$ $\Rr$-divisor $H$. Let $h: (X',\Ff',B';G)/Z\rightarrow (X,\Ff,B)$ be an ACSS modification of $(X,\Ff,B)$ that is $\Qq$-factorial, strict, and super. By Proposition \ref{prop: lift mmp}, any infinite sequence of steps of a $(K_{\Ff}+B+A)$-MMP$/U$ with scaling of $H$ induces an infinite sequence of steps of a $(K_{\Ff'}+B'+h^*A)$-MMP$/U$ with scaling of $h^*H$. By Theorem \ref{thm: take strict simple model run mmp}, any $(K_{\Ff'}+B'+h^*A)$-MMP$/U$ with scaling of $h^*H$ terminates with a minimal model$/U$. Thus any $(K_{\Ff}+B+A)$-MMP$/U$ with scaling of $H$ terminates with a minimal model of $(X,\Ff,B+A)/U$, and the theorem follows.
\end{proof}

\begin{proof}[Proof of Theorem \ref{thm: bpf intro}]
    Let $H:=K_{\Ff}+B+A$. By Lemma \ref{lem: lc+ample nef=nqc}, $H$ is NQC$/U$. Thus we have $H=\sum a_iH_i$ for some nef$/U$ Cartier divisors $H_i$ on $X$ and $a_i>0$ for each $i$. Let $\epsilon_0:=\min\{a_i\}$ and let $l>\frac{2\dim X}{\epsilon_0}$ be an integer.

    Let $0<e\ll 1$ be a real number such that $$\widehat{A}:=A+e(K_{\Ff}+B)-e(K_X+\Delta)$$
    is ample$/U$. Let $K:=(1-e)(K_{\Ff}+B)+e(K_X+\Delta)$, then $H=K+\widehat{A}$.

    \begin{claim}\label{claim: K_F+hat A is pe}
        $K_{\Ff}+B+\widehat{A}$ is pseudo-effective$/U$.
    \end{claim}
    Assume Claim \ref{claim: K_F+hat A is pe}, then $K_{\Ff}+B+\widehat{A}+lH$ is also pseudo-effective$/U$. By Theorem \ref{thm: eolmm+A 1 gfq} and Lemma \ref{lem: single step enef trivial}, there exists a positive integer $l\gg 0$ which does not depend on $e$, such that we may run a $(K_{\Ff}+B+\widehat{A}+lH)$-MMP$/U$, each step of this MMP is $H$-trivial, and the MMP terminates with a minimal model$/U$. Let $Y$ be the output of this MMP and let $B_Y,\Delta_Y,H_Y,A_Y$ be the images of $B,\Delta,H,A$ on $Y$ respectively. Let $\Aa:=\overline{\widehat{A}}$. By Theorem \ref{thm: mmp can run gfq}, $(Y,\Delta_Y,\Aa)$ is klt. Since
    $$(l+1-el)H=K+\widehat{A}+(l-el)H=(1-e)(K_{\Ff}+B+\widehat{A}+lH)+e(K_X+\Delta+\widehat{A}),$$
    we have
    $$K_Y+\Delta_Y+\Aa_Y+\frac{1-e}{e}(K_{\Ff_Y}+B_Y+\Aa_Y+lH_Y)=\frac{l+1-el}{e}H_Y$$
    is nef$/U$. Let $$\Pp:=\Aa+\overline{\frac{1-e}{e}(K_{\Ff_Y}+B_Y+\Aa_Y+lH_Y)},$$
    then $(Y,\Delta_Y,\Pp)/U$ is klt. Since $\widehat{A}$ is ample$/U$, $\Aa_Y$ is big$/U$. By \cite[Lemma 4.4(2)]{BZ16}, $(Y,\Delta_Y,\Pp)/U$ is a good minimal model of itself, so $K_Y+\Delta_Y+\Pp_Y$ is semi-ample$/U$. Thus $H_Y$ is semi-ample$/U$. Since $X\dashrightarrow Y$ is $H$-trivial, $H$ is semi-ample$/U$. Moreover, if $H$ is Cartier, then $H_Y$ is Cartier. By Lemma \ref{lem: klt bpf}, $\mathcal{O}_Y(nH_Y)$ is globally generated$/U$ for any integer $n\gg 0$, and so $\mathcal{O}_X(nH)$ is globally generated$/U$ for any integer $n\gg 0$.
    
\vspace{.5em}

   Finally we give the proof of Claim \ref{claim: K_F+hat A is pe}: 
    
Let $h: (X',\Ff',B';G)/Z\rightarrow (X,\Ff,B)$ be an ACSS modification that is $\Qq$-factorial, strict and super. Suppose $K_{\Ff}+B+\widehat A$ is not pseudo-effective$/U$, then $K_{\Ff'}+B'+\Aa_{X'}$ is also not pseudo-effective$/U$. By \cite[Proposition 9.3.2]{CHLX23}, we may run a $(K_{\Ff'}+B'+\Aa_{X'})$-MMP$/U$ which terminates with a Mori fiber space$/U$, and this MMP is also an MMP$/Z$. Since
$$K_{X'}+B'+G+\Aa_{X'}\sim_{\mathbb R,Z}K_{\Ff'}+B'+\Aa_{X'},$$
we have
$$(1-e)K_{\Ff'}+eK_{X'}+B'+eG+\Aa_{X'}\sim_{\mathbb R,Z}K_{\Ff'}+B'+\Aa_{X'},$$
so this MMP is also a $((1-e)K_{\Ff'}+eK_{X'}+B'+eG+\Aa_{X'})$-MMP$/U$. In particular, $(1-e)(K_{\Ff'}+B'+\Aa_{X'})+e(K_{X'}+B'+G+\Aa_{X'})$ is not pseudo-effective$/U$.

Since $(X,\Delta)$ is klt and $B\ge\Delta$, from the definition of ACSS modification we see that
$$K_{X'}+B'+G\ge K_{X'}+\Delta'=h^*(K_X+\Delta),$$
hence $(1-e)(K_{\Ff'}+B'+\Aa_{X'})+e(K_{X'}+\Delta'+\Aa_{X'})$ is not pseudo-effective$/U$ either. However, this is not possible as
\[
(1-e)(K_{\Ff'}+B'+\Aa_{X'})+e(K_{X'}+\Delta'+\Aa_{X'})=h^*(K_{\Ff}+B+A)=h^*H.
\]  
\end{proof}

\begin{proof}[Proof of Theorem \ref{thm: eolmm+A 1}(2)]
By Theorem \ref{thm: eolmm+A 1}(1) we can run a $(K_\Ff+B+A)$-MMP$/U$ and it terminate with a minimal model $\phi:X\dashrightarrow X'$. In particular the corresponding birational transform $K_{\Ff'}+B'+A'$ is nef$/U$. Notice that $(X',\Delta':=\phi_*\Delta)$ is klt by Theorem \ref{thm: mmp can run gfq}. By \cite[Lemma 16.1.1]{CHLX23}, there exists a nef$/U$ {\bf b}-divisor $\Nn$ and an ample$/U$ $\Rr$-divisor $\tilde A'$ such that 
\begin{enumerate}
    \item $(X',\Ff',B',\Nn)$ is lc, and
    \item $A'\sim_{\Rr,U}\tilde A'+\Nn_{X'}.$
\end{enumerate}
Then it suffices to show $K_{\Ff'}+B'+\Nn_{X'}+\tilde A'$ is semi-ample$/U$, which follows from Theorem \ref{thm: bpf intro gfq}.
\end{proof}

\begin{proof}[Proof of Theorem \ref{thm: fg intro}]
This is an immediate consequence of Theorem \ref{thm: eolmm+A 1}(2).  
\end{proof}

\begin{proof}[Proof of Theorem \ref{thm: nt eogmm}]
By Theorem \ref{thm: can run mmp scaling}, we can run a $(K_{\Ff}+B)$-MMP with scaling of any ample $\Rr$-divisor $A$. Let $h: (X',\Ff',B';G)/Z\rightarrow (X,\Ff,B)$ be an ACSS modification of $(X,\Ff,B)$ that is $\Qq$-factorial, strict, and super. By Proposition \ref{prop: lift mmp}, any infinite sequence of steps of a $(K_{\Ff}+B)$-MMP with scaling of $A$ induces an infinite sequence of steps of a $(K_{\Ff'}+B')$-MMP with scaling of $h^*A$. 

Suppose that this MMP does not terminate. Let $\lambda_i$ be the scaling numbers of this MMP and let $\lambda:=\lim_{i\rightarrow+\infty}\lambda_i$. If $\lambda>0$, then we have an infinite sequence of steps of a $(K_{\Ff'}+B'+\lambda A')$-MMP with scaling of $h^*A$, which contradicts Theorem \ref{thm: take strict simple model run mmp}. Therefore, $\lambda=0$. Let
$$(X_0,\Ff_0,B_0):=(X',\Ff',B')\dashrightarrow (X_1,\Ff_1,B_1)\dashrightarrow\dots\dashrightarrow (X_n,\Ff_n,B_n)\dashrightarrow\dots$$
be this MMP and let $A_i$ be the image of $h^*A$ on $X_i$ for each $i$. Then there exists $n>0$ such that the induced birational map $\phi_i: X_n\dashrightarrow X_i$ is small for any $i\geq n$. Therefore,
$$K_{\Ff_n}+B_n=\lim_{i\rightarrow+\infty}(\phi_i^{-1})_*(K_{\Ff_i}+B_i+\lambda_i A_i)$$
is movable. Since $\kappa_{\sigma}(K_{\Ff}+B)=0$, $\kappa_{\sigma}(K_{\Ff_n}+B_n)=0$. By \cite[Lemma 4.2.4]{CHLX23}, $K_{\Ff_n}+B_n\equiv 0$, a contradiction. Therefore, any $(K_{\Ff}+B)$-MMP with scaling of $A$ terminates with a minimal model $(X_{\min},\Ff_{\min},B_{\min})$ of $(X,\Ff,B)$ such that $K_{\Ff_{\min}}+B_{\min}\equiv 0$. By \cite[Theorem 1.4]{DLM23}, $K_{\Ff_{\min}}+B_{\min}\sim_{\mathbb R}0$.
\end{proof}

\begin{proof}[Proof of Theorem \ref{thm: weak mfs}(1)]
    Let $\epsilon$ be a positive real number such that $H:=-K_{\Ff}+\epsilon D$ is ample. Then any $D$-MMP is a $(K_{\Ff}+H)$-MMP. Since $(X,\Ff,0,\bar H)$ is an lc generalized foliated quadruple and $X$ is klt, by the generalized pair version of Theorems \ref{thm: eomfs} and \ref{thm: eolmm+A 1} (see Theorems \ref{thm: eomfs gfq} and \ref{thm: eolmm+A 1 gfq}), we may run a $(K_{\Ff}+H)$-MMP with scaling of an ample divisor which terminates with either a good minimal model or a Mori fiber space. 
\end{proof}

\begin{lem}\label{lem: lift fano foliation over proj bundle}
Let $(X,\Ff,B)/U$ be an lc algebraically integrable foliated triple such that $(X,\Delta)$ is klt for some $B\geq\Delta\geq 0$. Assume that $-(K_{\Ff}+B)$ is ample$/U$ and $L_1,...,L_m$ be Cartier divisors on $X$. Let $\mathcal{E}:=\oplus_{i=1}^m\Oo_X(L_i)$ and $\pi: Y=\mathbb{P}_X(\mathcal{E})\to X$ be the corresponding projective bundle. Then $(Y,\pi^*\Delta)$ is klt and there exists $B_Y\ge \pi^*\Delta$ such that
\begin{enumerate}
    \item $(Y,\pi^{-1}\Ff,B_Y)$ is lc,
    \item $-(K_{\pi^{-1}\Ff}+B_Y)$ is ample$/U$.
\end{enumerate}
\end{lem}
\begin{proof}
  $(Y,\pi^*\Delta)$ is klt by the smoothness of $\pi$. Let $T$ be the union of the sections of $\pi$ corresponding to $\mathcal{E}\twoheadrightarrow L_i, 1\le i\le m$.
We will show $B_Y:=(1-\epsilon)T+\pi^{*}B$ satisfies our requirements for any $0<\epsilon\ll1$.

By taking a foliated log resolution of $(X,\Ff,B)$ and considering the base change of $\pi$, we can easily see by definition that $(Y,\pi^{-1}\Ff,T+\pi^{*}B)$ is lc if and only if $(X,\Ff,B)$ is lc. Hence $(Y,\pi^{-1}\Ff,B_Y)$ is lc as well.

Since $K_{\pi^{-1}\Ff}+T=\pi^*K_\Ff$ and $T$ is $\pi$-ample, we have
  \begin{align*}
      -(K_{\pi^{-1}\Ff}+B_Y)\sim_{\Rr,U}-\pi^*(K_\Ff+B)+\epsilon T
  \end{align*}
is ample$/U$ for any $0<\epsilon\ll1$.
\end{proof}

\begin{proof}[Proof of Theorem \ref{thm: weak mfs}(2)]
By Theorem \ref{thm: bpf intro} we have $\Pic(X)_\Qq\simeq N^1(X)_\Qq$. Choose a basis $L_1,...,L_m$ of $\Pic(X)_\Qq$ such that each $L_i$
is a Cartier divisor on $X$ and the convex hull of $L_1,...,L_m$ in $N^1(X)_\Rr$ contains the effective cone Eff(X). Let $Y:=\mathbb{P}_X(\oplus_{i=1}^m\Oo_X(L_i))$ and $H$ be the Cartier divisor which corresponds to the tautological line bundle $\Oo_Y(1)$. Note that $H$ is a big Cartier divisor on $Y$. It is easy to check that the Cox ring of $X$ is finitely generated if and only if the section ring $R(Y, H)$ is finitely generated. By Lemma \ref{lem: lift fano foliation over proj bundle}, there is a boundary divisor $B_Y$ on $Y$ such
that $(Y,\pi^{-1}\Ff,B_Y)$ is lc and $-(K_{\pi^{-1}\Ff}+B_Y)$ is ample. We can
choose an ample divisor $A$ on $Y$ such that $K_{\pi^{-1}\Ff}+B_Y+A\sim_\Qq\delta H$ for a sufficiently small rational number $\delta>0$. By
Theorem \ref{thm: fg intro}, $R(Y,\delta H)$ is finitely generated, and so is $R(Y,H)$. Therefore, the Cox ring of $X$ is also finitely generated and hence $X$ is a Mori dream space.
\end{proof}

\begin{proof}[Proof of Theorem \ref{thm: eolmm+A}]

Let $h: (X',\Ff',B';G)/Z\rightarrow (X,\Ff,B)$ be an ACSS model of $(X,\Ff,B)$ that is $\Qq$-factorial, strict, and super. Let $A':=h^*A$. Let $H$ be an ample $\Rr$-divisor on $X$ and let $H':=h^*H$. By Theorem \ref{thm: take strict simple model run mmp}, we may run a $(K_{\Ff'}+B'+A')$-MMP$/U$ $\phi': X'\dashrightarrow X''$ which terminates with either a minimal model $(X'',\Ff'',B''+A'')/U$ or a Mori fiber space $(X'',\Ff'',B''+A'')\rightarrow T$ of $(X',\Ff',B'+A')/U$, where $B''$ and $A''$ are the images of $B'$ and $A'$ on $X''$ respectively. Let $\phi: X\dashrightarrow X''$ be the induced birational map.

For any prime divisor $D$ that is extracted by $\phi^{-1}$, $D$ is also extracted by $h$, so $$-\mult_D(B''+A'')=a(D,\Ff,B+A)\leq a(D,\Ff,B)=-\epsilon_{\Ff}(D).$$ 
Therefore, $(X'',\Ff'',B''+A'')$ is a log birational model of $(X,\Ff,B)$. For any prime divisor $D$ on $X$ that is exceptional$/X''$, $D$ is also a prime divisor on $X'$ that is exceptional$/X''$, so 
$$a(D,\Ff,B+A)=a(D,\Ff',B'+A')<a(D,\Ff'',B''+A'').$$
Thus either $(X'',\Ff'',B''+A'')/U$ is a bs-minimal model of $(X,\Ff,B+A)/U$, or $(X'',\Ff'',B''+A'')\rightarrow T$ is a bs-Mori fiber space of $(X,\Ff,B+A)/U$. In particular, $(X,\Ff,B+A)/U$ has a bs-minimal model or a bs-Mori fiber space.
\end{proof}

\begin{proof}[Proof of Theorem \ref{thm: foliation emm equal to tof scaling}]
There exists an ACSS modification $h: (X',\Ff',B';G)/Z\rightarrow (X,\Ff,B)$ 
that is $\Qq$-factorial, strict, and super, and $(X',B'+G)/U$ has a log minimal model if $(X,\Ff,B)/U$ has a bs-minimal model by Proposition \ref{prop: eolmm foliation to pair}. Let $f: X'\rightarrow Z$ be the associated contraction and let $A':=h^*A$. By Proposition \ref{prop: lift mmp}, any infinite sequence of $(K_{\Ff}+B)$-MMP$/U$ with scaling of $A$ induces an infinite sequence of $(K_{\Ff'}+B')$-MMP$/U$ with scaling of $A'$. 

First we prove (1). Suppose that the MMP does not terminate. We let $\lambda_i$ be the scaling numbers of the $(K_{\Ff}+B)$-MMP$/U$ with scaling of $A$ and let $\lambda:=\lim_{i\rightarrow+\infty}\lambda_i$. By Lemma \ref{lem: chlx 9.2.1}, any $(K_{\Ff'}+B')$-MMP$/U$ with scaling of $A'$ is also a $(K_{X'}+B'+G)$-MMP$/U$ with $A'$, and $\lambda$ is also the limit of the scaling numbers of the $(K_{X'}+B'+G)$-MMP$/U$ with $A'$. By Theorem \ref{thm: take strict simple model run mmp}, we have $\lambda=0$. In particular, $\lambda\not=\lambda_i$ for any $i$ and $K_{\Ff}+B$ is pseudo-effective$/U$. Thus $(X,\Ff,B)/U$ has a bs-minimal model, so $(X',B'+G)/U$ has a log minimal model. This contradicts Theorem \ref{thm: bir12 1.9(3)}.

(2) follows from (1) and Theorem \ref{thm: can run mmp scaling}.
\end{proof}

\section{Further discussions}\label{sec: problem}

\subsection{New definition of foliated klt singularities} \begin{rem}\label{rem: new definition klt}
       We remark that our definitions of lc and klt singularities in Definition \ref{defn: foliation singularity} have some differences with the classical definitions \cite[Definition I.1.5]{McQ08}, where the $-1$ is replaced with $-\epsilon_{\Ff}(E)$. The other definition is used in most literature (e.g. \cite{CS20,ACSS21,CS21,CHLX23}). Lemma \ref{lem: equivalence definition lc} shows that our definition of ``lc" coincides with the classical definition of ``lc". We briefly explain why we change the definition of ``klt". This is with several reasons:
       \begin{enumerate}
      \item  The classical ``klt" is an empty condition in many scenarios. For any non-trivial algebraically integrable foliation $\Ff$ ($\Ff\not=T_X$), there are a lot of $\Ff$-invariant divisors, and each of them is an lc place. Therefore, the condition
         \begin{center}
         $a(E,\Ff,0)>-\epsilon_{\Ff}(E)$ for any prime divisor $E$ over $X$
         \end{center}
         as in \cite[Definition I.1.5(3)]{McQ08} is a condition that cannot be satisfied by any non-trivial algebraically integrable foliation. Similar issues may appear for foliations with non-trivial algebraic parts.
         \item ``Plt" is missing. For example, \cite[Theorem 1.1]{CS25b} established the correspondence via adjunction to non-invariant divisors for lc singularities. But the ``plt-klt" correspondence is missing. This prevents us to prove a lot of things, e.g. the existence of ``pl-flips" for foliations in dimension $4$. 
         \item ``Terminal" is also missing. \cite[Theorems 1.1, 3.16]{CS25b} can only show that ``adjunction of canonical (resp. terminal) singularities to non-invariant divisors is canonical (resp. terminal)". However, for usual pairs, we know that ``adjunction of canonical singularities to divisors is terminal". One reason for this is that the definition of ``terminal" for foliations requires that the discrepancies of $\Ff$-invariant divisors are $>0$. Thus it is also natural to ask whether we can establish the ``canonical-terminal" correspondence if we ignore the invariant lc places.
         \item There are substantial differences between non-invariant divisors and invariant divisors and it is very important to use non-invariant divisors to lift sections. This is because we usually have exact sequences of the form
         $$\mathcal{O}_X(L-S)\rightarrow\mathcal{O}_X(L)\rightarrow\mathcal{O}_S(L|_S)$$
         which allow us to lift sections. Here $L$ usually has an lc structure $K+B$ and $S$ is a component of  $\lfloor B\rfloor$. However, if $K=K_{\Ff}$ and $S$ is an $\Ff$-invariant divisor, then $(X,\Ff,B)$ will not be lc by \cite[Remark 2.3]{CS21}. Therefore, non-invariant lc places behave better than lc places in this scenario.
       \end{enumerate}
\end{rem}
With the above discussion, we also propose the definition of ``plt":
\begin{defn}
    Let $(X,\Ff,B)$ be a foliated triple. We say that $(X,\Ff,B)$ is \emph{plt} if  $a(E,\Ff,B)>-\epsilon_{\Ff}(E)$ for any prime divisor $E$ that is exceptional over $X$.
\end{defn}

We propose the following questions on foliations with klt singularities.
\begin{ques}[cf. {\cite[Theorem 1.1]{CS25b}}]
   Let $(X,\Ff,B)$ be a ($\mathbb Q$-factorial) plt foliated triple and $S$ a component of $\lfloor B\rfloor$ with normalization $S^\nu$. Let $\Ff_{S^\nu}$ be the restricted foliation of $\Ff$ on $S^\nu$ and let
   $$K_{\Ff_{S^\nu}}+\Diff_{S^\nu}(\Ff,B):=(K_{\Ff}+B)|_{S^\nu}.$$
   Is
   $(S^\nu,\Ff_{S^\nu},\Diff_{S^\nu}(\Ff,B))$ klt?
\end{ques}

\begin{ques}[cf. {\cite[Conjecture 4.2(2)]{CS25a}}]
    Let $(X,\Ff,B)$ be a $\Qq$-factorial klt algebraically integrable foliated triple such that $(X,B)$ is klt. Is $\Ff$ induced by a contraction?
\end{ques}

\begin{ques}[cf. {\cite[Theorem 11.3]{CS21}}]
    Let $(X,\Ff,B)$ be a klt foliated triple such that $\dim X=3$ and $\rk\Ff=2$. Is $\Ff$ non-dicritical?
\end{ques}

Finally, we remark that the existence of pl-flips is one crucial step towards the existence of flips for usual varieties. With this in mind, we ask the following:

\begin{ques}[Pl-flip]
    Let $(X,\Ff,B)$ be a $\Qq$-factorial plt projective foliated triple and $f: X\rightarrow Z$ a small contraction such that $\rho(X/Z)=1$, $-(K_\Ff+B)$ is ample$/Z$, $S:=\lfloor B\rfloor$ is irreducible, and $-S$ is ample$/Z$. Assume that $\dim X=4$. Does the flip $X^+\rightarrow Z$ of $f$ exist?
\end{ques}

\subsection{MMP when the ambient variety is not klt} 
We still expect that the minimal model program holds for lc algebraically integrable foliations even if the ambient variety is not necessarily klt. We propose the following conjecture:
\begin{conj}[MMP for algebraically integrable foliations]\label{conj: mmp foliation}
    Let $(X,\Ff,B)/U$ be an lc algebraically integrable foliated triple and $R$ a $(K_{\Ff}+B)$-negative extremal ray. Then:
    \begin{enumerate}
        \item (Contraction theorem) There exists a contraction$/U$ $\cont_R: X\rightarrow T$ of $R$.
        \item (Existence of flips) If $\cont_R$ is a flipping contraction, then the flip$/U$ $X^+\rightarrow T$ associated to $R$ exists.
        \item (MMP) We can run a $(K_{\Ff}+B)$-MMP$/U$.
    \end{enumerate}
\end{conj}
Theorem \ref{thm: eolmm+A} provides some positive evidence towards Conjecture \ref{conj: mmp foliation}. In fact, if the ``minimal model in the sense of Birkar-Shokurov" in Theorem \ref{thm: eolmm+A} is replaced with ``\emph{good} minimal model in the sense of Birkar-Shokurov", then Conjecture \ref{conj: mmp foliation} will immediately follow.

Another positive evidence for Conjecture \ref{conj: mmp foliation} is the case when $\Ff$ is induced by a locally stable family $f: (X,B)\rightarrow Z$. In this case, Conjecture \ref{conj: mmp foliation} is essentially settled in \cite[Theorem 1.5]{MZ23} although it is not written in the language of foliations. We reinterpret \cite[Theorem 1.5]{MZ23} in the following way:

\begin{thm}[{\cite[Theorem 1.5]{MZ23}}]\label{thm: mmp locally stable variation}
    Let $f: (X,B)\rightarrow Z$ be a locally stable family over a normal variety with normal generic fiber and $\Ff$ the foliation induced by $f$. Then  $K_{\Ff}=K_{X/Z}$ and $(X,\Ff,B)$ is lc. Moreover, we may run a $(K_{\Ff}+B)$-MMP$/Z$.
\end{thm}
\begin{proof}
It follows from the definition of locally stable families \cite[Theorem-Definition 4.7]{Kol23} and \cite[Theorem 1.5]{MZ23}.
\end{proof}
We also have the following result for lc algebraically integrable foliations of co-rank $1$.

\begin{thm}\label{thm: corank 1}
    Let $(X,\Ff,B)/U$ be an lc algebraically integrable foliated triple such that $\rk\Ff=\dim X-1$ and let $A,H$ be two ample$/U$ $\Rr$-divisors on $X$. Then:
    \begin{enumerate}
        \item If $(X,\Ff,B)$ is klt, then $\Ff$ is induced by a contraction $f: X\rightarrow Z$.
        \item Suppose that $\Ff$ is induced by a contraction $f: X\rightarrow Z$. Then:
        \begin{enumerate}
            \item We may run a $(K_{\Ff}+B)$-MMP$/U$.
            \item We may run a $(K_{\Ff}+B+A)$-MMP$/U$ with scaling of $H$ which terminates with either a good minimal model or a Mori fiber space of $(X,\Ff,B+A)/U$.
        \end{enumerate}
    \end{enumerate}
\end{thm}
\begin{proof}
    (1) By \cite[Theorem 2.1.10]{CHLX23}, $\Ff$ is induced by an almost holomorphic map $f: X\dashrightarrow Z$ where $Z$ is a curve. By the rigidity lemma, $f$ is a contraction.

    (2) Let $h: (X',\Ff',B';G')/Z'\rightarrow (X,\Ff,B)$ be a super ACSS modification of $(X,\Ff,B)$ whose existence is guaranteed by Theorem \ref{thm: eo acss model}. Since $Z'$ is a curve, $Z'\cong Z$. Therefore, $X'\rightarrow Z'$ factors through $X$, so $X$ is the core model of $(h, X'\rightarrow Z)$. Thus $(X,B+G)$ is lc, where $G:=h_*G'$, and $K_X+B+G\sim_{\mathbb R,Z}K_\Ff+B$ by Lemma \ref{lem: existence of core model}. (2.a) follows from \cite[Lemma 9.1.4]{CHLX23} and (2.b) follows from \cite[Theorem 16.1.4]{CHLX23}.
\end{proof}

Finally, we conjecture the following:

\begin{conj}[Base-point-freeness]\label{conj: base-point-freeness theorem}
    Let $(X,\Ff,B)/U$ be an lc algebraically integrable foliated triple. Let $A$ be an ample$/U$ $\Rr$-divisor on $X$ such that $K_{\Ff}+B+A$ is nef$/U$. Then:
    \begin{enumerate}
        \item $K_{\Ff}+B+A$ is semi-ample$/U$.
        \item If $K_{\Ff}+B+A$ is Cartier, then $\mathcal{O}_X(m(K_{\Ff}+B+A))$ is globally generated$/U$ for any integer $m\gg 0$.
    \end{enumerate}
\end{conj}

\appendix

\section{Generalized foliated quadruples}\label{sec: gfq case}

In this appendix, we discuss the generalized foliated quadruple version of our main theorems. Due to technicality, we shall omit or only sketch the proofs of these theorems in this appendix. This appendix is organized in the following way: first we shall define generalized foliated quadruples and its related concepts. Then we shall state the generalized foliated quadruple version of the main theorems of the paper and provide sketch of their proofs which relies on results in the rest part of this appendix. Finally, we shall provide the generalized foliated quadruple version of all other results in this paper.

\subsection{Definitions}
\begin{defn}
A \emph{generalized foliated quadruple} $(X,\Ff,B,\Mm)/U$ consists of a normal quasi-projective variety $X$, a foliation $\Ff$ on $X$, an $\Rr$-divisor $B\geq 0$ on $X$, a projective morphism $X\rightarrow U$, and a nef$/U$ $\bb$-divisor $\Mm$, such that $K_{\Ff}+B+\Mm_X$ is $\mathbb R$-Cartier. We say that $(X,\Ff,B,\Mm)/U$ is \emph{NQC} if $\Mm$ is NQC$/U$. If $\Ff=T_X$ then we say that $(X,B,\Mm)/U$ is a \emph{generalized pair}. 
\end{defn}

\begin{rem}\label{rem: defs for gfq and g-pairs} We briefly remark the definitions of other concepts of  generalized foliated quadruples and generalized pairs.
\begin{enumerate}
\item   Related descriptions of generalized foliated quadruples are defined in the same way as in Definition \ref{defn: foliated triple} (e.g. generalized foliated sub-quadruple).
\item     Singularities of generalized foliated quadruples are defined in the same way as in Definition \ref{defn: foliation singularity}. We do not define ``potentially generalized klt" as in Definition \ref{defn: potentially klt} because it is equivalent to ``potentially klt". See Lemma \ref{lem: gklt is klt gfq} below.
\item Foliated log resolution (Definition \ref{defn: log resolution}) and foliated log smooth model (Definition \ref{defn: log smooth models}) can be defined for  generalized foliated quadruples in the same way by requiring that $\Mm$ descends $X'$.
\item Property $(*)$ (Definition \ref{defn: foliation property *}) and ACSS (Definition \ref{defn: ACSS f-triple}), can be defined for generalized foliated quadruples in a similar way. The ``qdlt" used for Definition \ref{defn: ACSS f-triple}(3) shall be replaced by ``qdlt" for generalized pairs defined in \cite[Definition 7.1.1]{CHLX23}.
\item Simple, core, ACSS modification/models defined in Definitions \ref{defn: simple model} and \ref{defn: core model} can be defined for generalized foliated quadruples in the same way.
\item Log birational model (Definition \ref{defn: log birational model}), all types of models in Definition \ref{defn: minimal model}, and different types of Mori fiber spaces in Definition \ref{defn: mfs}, can be defined in the same way for generalized foliated quadruples by following the principle that when taking these models, the nef part $\Mm$ does not change.
\end{enumerate}
\end{rem}

\subsection{Generalized foliated quadruple version of the main theorems}

\begin{thm}[Theorem \ref{thm: main}]\label{thm: main gfq}
Let $(X,\Ff,B,\Mm)/U$ be an lc algebraically integrable generalized foliated quadruple such that $(X,\Delta,\Nn)/U$ is klt, where $B\geq\Delta\geq 0$ and $\Mm-\Nn$ is nef$/U$.  Let $R$ be a $(K_{\Ff}+B+\Mm_X)$-negative extremal ray$/U$. Then:
    \begin{enumerate}
        \item (Contraction theorem) There exists a contraction$/U$ $\cont_R: X\rightarrow Z$ of $R$.
        \item (Existence of flips) If $\cont_R$ is a flipping contraction, then the flip$/U$ $X^+\rightarrow T$ associated to $R$ exists.
    \end{enumerate}   
\end{thm}
\begin{proof}
    It follows from Proposition \ref{prop: contraction not big gfq} and Theorem \ref{thm: cont and flip with detail}.
\end{proof}

\begin{thm}[Theorem \ref{thm: mmp can run}]\label{thm: mmp can run gfq}
Let $(X,\Ff,B,\Mm)/U$ be an lc algebraically integrable generalized foliated quadruple such that $(X,\Delta,\Nn)/U$ is klt, where $B\geq\Delta\geq 0$ and $\Mm-\Nn$ is nef$/U$. Then we may run a $(K_{\Ff}+B+\Mm_X)$-MMP$/U$. Moreover, for any birational map $\phi: X\dashrightarrow X^+$ that is a sequence of steps of a $(K_{\Ff}+B+\Mm_X)$-MMP$/U$, $(X,\Delta^+:=\phi_*\Delta,\Nn)$ is klt.
\end{thm}
\begin{proof}[Proof of Theorem \ref{thm: mmp can run}]
    By Proposition \ref{prop: contraction not big gfq} and Theorem \ref{thm: cont and flip with detail}, we can run a step of a $(K_{\Ff}+B)$-MMP$/U$. By Theorem \ref{thm: cont and flip with detail}(2.a), after a step of the MMP $\phi: X\dashrightarrow X'$ that is not a Mori fiber space, $(X',\Delta':=\phi_*\Delta)$ is klt. Thus we may continue this process. 
\end{proof}

\begin{thm}[Theorem \ref{thm: can run mmp scaling}]\label{thm: can run mmp scaling gfq}
Let $(X,\Ff,B,\Mm)/U$ be an lc algebraically integrable generalized foliated quadruple such that $(X,\Delta,\Nn)/U$ is klt, where $B\geq\Delta\geq 0$ and $\Mm-\Nn$ is nef$/U$. Let $A$ be an ample$/U$ $\Rr$-divisor on $X$. Then we may run a $(K_{\Ff}+B+\Mm_X)$-MMP$/U$ with scaling of $A$. 
\end{thm}
\begin{proof}
    It follows from Theorem \ref{thm: mmp can run gfq} and Lemma \ref{lem: can run mmp scaling gfq}.
\end{proof}

\begin{thm}[Theorem \ref{thm: eomfs}]\label{thm: eomfs gfq}
Let $(X,\Ff,B,\Mm)/U$ be an lc algebraically integrable generalized foliated quadruple such that $(X,\Delta,\Nn)/U$ is klt, where $B\geq\Delta\geq 0$ and $\Mm-\Nn$ is nef$/U$. Assume that $K_{\Ff}+B+\Mm_X$ is not pseudo-effective$/U$. Then we may run a $(K_{\Ff}+B+\Mm_X)$-MMP$/U$ with scaling of an ample$/U$ $\Rr$-divisor and any such MMP terminates with a Mori fiber space$/U$.
\end{thm}
\begin{proof}
    It follows from the same lines of the proof of Theorem \ref{thm: eomfs} except that we replace Theorems \ref{thm: take strict simple model run mmp}, \ref{thm: can run mmp scaling} and Proposition \ref{prop: lift mmp} with Theorems \ref{thm: take strict simple model run mmp gfq}, \ref{thm: can run mmp scaling gfq} and Proposition \ref{prop: lift mmp gfq} respectively.
\end{proof}

\begin{thm}[Theorem \ref{thm: eolmm+A 1}]\label{thm: eolmm+A 1 gfq} 
Let $(X,\Ff,B,\Mm)/U$ be an lc algebraically integrable generalized foliated quadruple such that $(X,\Delta,\Nn)/U$ is klt, where $B\geq\Delta\geq 0$ and $\Mm-\Nn$ is nef$/U$. Let $A$ be an ample$/U$ $\Rr$-divisor on $X$. Then: 
\begin{enumerate}
    \item We may run a $(K_{\Ff}+B+A+\Mm_X)$-MMP$/U$ with scaling of an ample$/U$ $\Rr$-divisor and any such MMP terminates with a minimal model of $(X,\Ff,B+A,\Mm)/U$.
    \item The minimal model in (1) is a good minimal model.
\end{enumerate}
\end{thm}
\begin{proof}
    It follows from the same lines of the proof of Theorem \ref{thm: eolmm+A 1} except that we replace Theorems \ref{thm: take strict simple model run mmp}, \ref{thm: can run mmp scaling} and Proposition \ref{prop: lift mmp} with Theorems \ref{thm: take strict simple model run mmp gfq}, \ref{thm: can run mmp scaling gfq} and Proposition \ref{prop: lift mmp gfq} respectively.
\end{proof}

\begin{thm}[Theorem \ref{thm: bpf intro}]\label{thm: bpf intro gfq}
  Let $(X,\Ff,B,\Mm)/U$ be an lc algebraically integrable generalized foliated quadruple such that $(X,\Delta,\Nn)/U$ is klt, where $B\geq\Delta\geq 0$ and $\Mm-\Nn$ is nef$/U$. Let $A$ be an ample$/U$ $\Rr$-divisor on $X$ such that $K_{\Ff}+B+A+\Mm_X$ is nef$/U$. Then:
  \begin{enumerate}
      \item $K_{\Ff}+B+A+\Mm_X$ is semi-ample$/U$.
      \item If $K_{\Ff}+B+A+\Mm_X$ is Cartier, then $\mathcal{O}_X(n(K_{\Ff}+B+A+\Mm_X))$ is globally generated$/U$ for any integer $n\gg 0$.
  \end{enumerate}
\end{thm}
\begin{proof}
    It follows the same proof of Theorem \ref{thm: bpf intro}, except we replace Lemma \ref{lem: lc+ample nef=nqc} with Lemma \ref{lem: lc+ample nef=nqc gfq}.
\end{proof}

\begin{thm}[Theorem \ref{thm: fg intro}]\label{thm: fg intro gfq}
  Let $(X,\Ff,B,\Mm)/U$ be an lc algebraically integrable generalized foliated quadruple such that $(X,\Delta,\Nn)$ is klt, where $B\geq\Delta\geq 0$ and $\Mm-\Nn$ is nef$/U$. Let $A$ be an ample$/U$ $\Rr$-divisor on $X$ such that $B+A+\Mm_X$ is a $\Qq$-divisor. Then the log canonical ring
$$R(X,K_{\Ff}+B+A+\Mm_X):=\oplus_{m=0}^{+\infty}\pi_*\mathcal{O}_X(\lfloor m(K_{\Ff}+B+A+\Mm_X)\rfloor)$$
  is a finitely generated $\mathcal{O}_U$-algebra.
\end{thm}
\begin{proof}
    This is an immediate consequence of Theorem \ref{thm: eolmm+A 1 gfq}(2).
\end{proof}

\begin{thm}[Theorem \ref{thm: nt eogmm}]\label{thm: nt eogmm gfq}
Let $(X,\Ff,B,\Mm)/U$ be an lc algebraically integrable generalized foliated quadruple such that $(X,\Delta,\Nn)/U$ is klt, where $B\geq\Delta\geq 0$ and $\Mm-\Nn$ is nef$/U$. Assume that $\kappa_{\sigma}(K_{\Ff}+B+\Mm_X)=0$.

    The we may run a $(K_{\Ff}+B+\Mm_X)$-MMP with scaling of an ample $\Rr$-divisor and any such MMP terminates with a minimal model $(X_{\min},\Ff_{\min},B_{\min},\Mm)$ of $(X,\Ff,B,\Mm)$ such that $K_{\Ff_{\min}}+B_{\min}+\Mm_{X_{\min}}\equiv 0$. Moreover, if $\kappa_{\iota}(K_{\Ff}+B+\Mm_X)=0$, then $K_{\Ff_{\min}}+B_{\min}+\Mm_{X_{\min}}\sim_{\mathbb R}0$.
\end{thm}
\begin{proof}
Except the last sentence of the proof where \cite[Theorem 1.4]{DLM23} is applied to show that $K_{\Ff_{\min}}+B_{\min}\sim_{\mathbb R}0$, the proof of Theorem \ref{thm: nt eogmm gfq} follows from the same lines of the proof of Theorem \ref{thm: eolmm+A 1} by replacing Theorems \ref{thm: take strict simple model run mmp}, \ref{thm: can run mmp scaling} and Proposition \ref{prop: lift mmp} with Theorems \ref{thm: take strict simple model run mmp gfq}, \ref{thm: can run mmp scaling gfq} and Proposition \ref{prop: lift mmp gfq} respectively. In this case, we get a minimal model $(X_{\min},\Ff_{\min},B_{\min},\Mm)$ of $(X,\Ff,B,\Mm)$ such that $K_{\Ff_{\min}}+B_{\min}+\Mm_{X_{\min}}\equiv 0$. The moreover part is obvious.
\end{proof}

\begin{thm}[Theorem \ref{thm: weak mfs}(1)]\label{thm: weak mfs gfq}
Let $(X,\Ff,B,\Mm)/U$ be an lc algebraically integrable generalized foliated quadruple such that $(X,\Delta,\Nn)/U$ is klt, where $B\geq\Delta\geq 0$ and $\Mm-\Nn$ is nef$/U$. Assume that $-(K_{\Ff}+B+\Mm_X)$ is ample$/U$. Let $D$ be an  $\Rr$-Cartier $\Rr$-divisor on $X$. Then we may run a $D$-MMP which terminates with either a good minimal model$/U$ of $D$ or a Mori fiber space$/U$ of $D$.
\end{thm}
\begin{proof}
It follows from the same lines of the proof of Theorem \ref{thm: weak mfs} except we replace Theorems \ref{thm: eomfs} and \ref{thm: eolmm+A 1} with Theorems \ref{thm: eomfs gfq} and \ref{thm: eolmm+A 1 gfq} respectively.
\end{proof}

\begin{thm}[Theorem \ref{thm: eolmm+A}]\label{thm: eolmm+A gfq} 
Let $(X,\Ff,B,\Mm)/U$ be an lc algebraically integrable generalized foliated quadruple and $A$ an ample$/U$ $\Rr$-divisor on $X$. Assume that either $X$ is potentially klt or $\Mm$ is NQC$/U$. Then $(X,\Ff,B,\Mm+\bar A)/U$ has either a minimal model or a Mori fiber space in the sense of Birkar-Shokurov.
\end{thm}
\begin{proof}
    It is an immediate consequence of Theorem \ref{thm: take strict simple model run mmp gfq}. Note that the proof is even simpler comparing to Theorem \ref{thm: eolmm+A gfq} since $A$ does not contribute to any singularity if it is in the nef part rather than the boundary part.
\end{proof}

\begin{thm}[Theorem \ref{thm: foliation emm equal to tof scaling}]\label{thm: foliation emm equal to tof scaling gfq}
Let $(X,\Ff,B,\Mm)/U$ be an NQC lc algebraically integrable generalized foliated quadruple. Assume that $(X,\Ff,B,\Mm)/U$ has a minimal model or a Mori fiber space in the sense of Birkar-Shokurov and $X$ is potentially klt. Let $A$ be an ample$/U$ $\Rr$-divisor on $X$. Then:
\begin{enumerate}
    \item  Any  $(K_{\Ff}+B+\Mm_X)$-MMP$/U$ with scaling of $A$ terminates.
    \item If there exists a klt generalized  pair $(X,\Delta,\Nn)$ such that $B\geq\Delta\geq 0$ and $\Mm-\Nn$ is nef$/U$, then $(X,\Ff,B,\Mm)/U$ has a minimal model or a Mori fiber space.
\end{enumerate}
\end{thm}
\begin{proof}
    It follows from the same lines of the proof of Theorem \ref{thm: foliation emm equal to tof scaling} except that we replace Lemma \ref{lem: chlx 9.2.1}, Propositions \ref{prop: eolmm foliation to pair} and \ref{prop: lift mmp}, Theorems \ref{thm: bir12 1.9(3)}, \ref{thm: take strict simple model run mmp}, and \ref{thm: can run mmp scaling} with Lemma \ref{lem: chlx 9.2.1 gfq}, Propositions \ref{prop: eolmm foliation to pair gfq} and \ref{prop: lift mmp gfq}, Theorems \ref{thm: bir12 1.9(3) gfq}, \ref{thm: take strict simple model run mmp gfq}, and \ref{thm: can run mmp scaling gfq} respectively. 
\end{proof}

\begin{thm}[Theorem \ref{thm: shokurov polytope foliation}]\label{thm: shokurov polytope foliation gfq}
    Let $(X,\Ff,B:=\sum_{i=1}^mv^0_iB_i,\Mm=\sum_{i=1}^n\mu^0_i\Mm_i)/Z$ be an lc algebraically integrable generalized foliated quadruple such that $K_{\Ff}+B+\Mm_X$ is nef$/Z$, each $B_i\geq 0$ is a Weil divisor, and each $\Mm_i$ is a nef$/Z$ Cartier $\bb$-divisor.
    
    Let $\bm{v}_0:=(v^0_1,\dots,v^0_m,\mu^0_1,\dots,\mu^0_n)$. Then there exists an open subset $U$ of the rational envelope of $\bm{v}_0$ in $\mathbb R^{m+n}$, such that $(X,\Ff,\sum_{i=1}^mv_iB_i,\sum_{i=1}^{n}\mu_i\Mm_i)$ is lc and $K_{\Ff}+\sum_{i=1}^mv_iB_i+\sum_{i=1}^{n}\mu_i\Mm_{i,X}$ is nef$/Z$ for any $(v_1,\dots,v_m,\mu_1,\dots,\mu_n)\in U$.
\end{thm}
\begin{proof}
    It follows from the same lines of the proof of Theorem \ref{thm: shokurov polytope foliation} except we replace \cite[Theorem 1.5]{DLM23} with \cite[Theorem 2.4.7]{CHLX23}.
\end{proof}

\subsection{Generalized foliated quadruple version of other results}

\begin{lem}[Lemma \ref{lem: equivalence definition lc}]\label{lem: equivalence definition lc gfq}
    Let $(X,\Ff,B,\Mm)$ be a generalized foliated sub-quadruple. The following two conditions are equivalent:
    \begin{enumerate}
        \item $(X,\Ff,B,\Mm)$ is sub-lc.
        \item $a(E,\Ff,B,\Mm)\geq -\epsilon_{\Ff}(E)$ for any prime divisor $E$ over $X$.
    \end{enumerate}
\end{lem}
\begin{proof}
It follows from the same lines of the proof of Lemma \ref{lem: equivalence definition lc}.
\end{proof}

\begin{lem}[Lemma \ref{lem: gklt is klt}]\label{lem: gklt is klt gfq}
    Let $(X,B,\Mm)/U$ be an lc g-pair such that $(X,\Delta_0,\Nn)/U$ is klt for some $\Delta_0,\Nn$, and let $A$ be an ample$/U$ $\Rr$-divisor on $X$. Then $X$ is potentially klt, and there exists a klt pair $(X,\Delta)$ such that $\Delta\sim_{\mathbb R,U}B+\Mm_X+A$.
\end{lem}
\begin{proof}
    It is \cite[Lemma 3.4]{HL22}.
\end{proof}

\begin{prop}[Proposition \ref{prop: weak cbf gfq}]\label{prop: weak cbf gfq gfq}
Let $(X,\Ff,B,\Mm)$ be a generalized foliated quadruple. Let $G\geq 0$ be a reduced divisor on $X$ and $f: X\rightarrow Z$ an equidimensional contraction, such that $(X,\Ff,B,\Mm;G)/Z$ satisfies Property $(*)$ and $B$ is horizontal$/Z$. Then
$$K_{\Ff}+B+\Mm_X\sim_{Z}K_X+B+G+\Mm_X.$$
\end{prop}
\begin{proof}
    It follows from \cite[Proposition 7.3.6]{CHLX23}.
\end{proof}

\begin{lem}[Lemma \ref{lem: bz16 4.4(3)}]\label{lem: bz16 4.4(3) gpair}
    Let $(X,B,\Mm)/U$ be a $\Qq$-factorial lc g-pair and $L$ an NQC$/U$ $\Rr$-divisor on $X$ such that $X$ is klt. Then there exists a positive real number $l_0$ such that any sequence of steps of a $(K_X+B+\Mm_X+lL)$-MMP$/U$ is $L$-trivial for any $l>l_0$.
\end{lem}
\begin{proof}
    It is \cite[Lemma 3.22]{HL22}.
\end{proof}

\begin{lem}[Lemma \ref{lem: lc+ample nef=nqc}]\label{lem: lc+ample nef=nqc gfq}
    Let $(X,\Ff,B,\Mm)/U$ be an lc algebraically integrable generalized foliated quadruple and let $D$ be an nef$/U$ $\Rr$-divisor on $X$ such that $D-(K_{\Ff}+B+\Mm_X)$ is ample$/U$. Then $D$ is NQC$/U$.
\end{lem}
\begin{proof}
    It follows from the same lines of the proof of Lemma \ref{lem: lc+ample nef=nqc}. Note that \cite[Theorem 2.3.1]{CHLX23} is applicable to generalized foliated quadruples.
\end{proof}

\begin{lem}[Lemma \ref{lem: trivial ray when -delta triple}]\label{lem: trivial ray when -delta gfq}
Let $(X,\Ff,B,\Mm)/U$ be an lc algebraically integrable generalized foliated quadruple and $D$ an $\Rr$-divisor on $X$, such that  $K_{\Ff}+B+\Mm_X+D$ is NQC$/U$. Then there exists $\delta_0\in (0,1)$, such that for any $\delta\in (0,\delta_0)$, any $(K_{\Ff}+B+\Mm_X+(1-\delta)D)$-non-positive extremal ray$/U$ is a $(K_{\Ff}+B+\Mm_X+D)$-trivial extremal ray$/U$. 
\end{lem}\begin{proof}
    It follows from the same lines of the proof of Lemma \ref{lem: trivial ray when -delta triple}. Note that \cite[Theorem 2.3.1]{CHLX23} is applicable to generalized foliated quadruples.
\end{proof}

\begin{thm}[Theorem \ref{thm: eo acss model}]\label{thm: eo acss model gfq}
    Let $(X,\Ff,B,\Mm)$ be an lc algebraically integrable generalized foliated quadruple. Then $(X,\Ff,B)$ has an ACSS model $h: (X',\Ff',B',\Mm;G)/Z\rightarrow (X,\Ff,B,\Mm)$ that is $\Qq$-factorial, strict, and super.
\end{thm}
\begin{proof}
    It is \cite[Theorem 2.5.1]{CHLX23}.
\end{proof}

\begin{lem}[Lemma \ref{lem: existence of core model}]\label{lem: existence of core model gfq}
    Let $(X,\Ff,B,\Mm)$ be an lc algebraically integrable generalized foliated quadruple and let $h: (X',\Ff',B',\Mm;G)/Z\rightarrow (X,\Ff,B,\Mm)$ be a simple model. Let $f: X'\rightarrow Z$ the associated contraction, and let $\bar X$ be the core model of $(h,f)$ associated with $(\bar h,\bar f)$. Let $g: X'\rightarrow\bar X$ be the induced birational morphism, $\bar\Ff:=g_*\Ff',\bar B:=g_*B'$, and $\bar G:=g_*G$. 
    
    Assume that $f$ is equidimensional. Then:
    \begin{enumerate}
    \item $K_{\Ff'}+B'+\Mm_{X'}=g^*(K_{\bar\Ff}+\bar B+\Mm_{\bar X})$.
        \item $K_{X'}+B'+G+\Mm_X=g^*(K_{\bar X}+\bar B+\bar G+\Mm_{\bar X})$.
        \item $\bar h: (\bar X,\bar\Ff,\bar B,\Mm;\bar G)/Z\rightarrow (X,\Ff,B,\Mm)$
        is a core model.
        \item If $h: (X',\Ff',B',\Mm;G)/Z\rightarrow (X,\Ff,B,\Mm)$ is strict (resp. super), then $\bar h: (\bar X,\bar\Ff,\bar B,\Mm;\bar G)/Z\rightarrow (X,\Ff,B,\Mm)$ is strict (resp. super).
    \end{enumerate}
\end{lem}
\begin{proof}
    It follows from the same lines of the proof of Lemma \ref{lem: existence of core model} except that we use Proposition \ref{prop: weak cbf gfq gfq} instead of Proposition \ref{prop: weak cbf gfq}.
\end{proof}

\begin{lem}[Lemma \ref{lem: existence of klt pair on strict simple models}]\label{lem: existence of klt pair on strict simple models gfq}
       Let $(X,\Ff,B,\Mm)$ be an lc algebraically generalized foliated quadruple and let $h: (X',\Ff',B',\Mm;G)/Z\rightarrow (X,\Ff,B,\Mm)$ be a strict simple model of $(X,\Ff,B,\Mm)$. If $X$ is potentially klt, then $X'$ is potentially klt.
\end{lem}
\begin{proof}
    By the same lines of the proof of Lemma \ref{lem: existence of klt pair on strict simple models}, we can show that there exists a klt generalized pair $(X',\Delta',\Nn)$. The lemma follows from Lemma \ref{lem: bz16 4.4(3) gpair}.
\end{proof}

\begin{lem}[Lemma \ref{lem: hl23 2.6}]\label{lem: hl23 2.6 gfq}
Let $(X,\Ff,B,\Mm)/U$ be a generalized foliated quadruple and let $(X',\Ff',B',\Mm)/U$ a bs-weak lc model of $(X,\Ff,B,\Mm)/U$ associated with the birational map $\phi: X\dashrightarrow X'$. Let $p: W\rightarrow X$ and $q: W\rightarrow X'$ be birational morphisms such that $q=\phi\circ p$. Assume that
$$p^*(K_\Ff+B+\Mm_X)=q^*(K_{\Ff'}+B'+\Mm_{X'})+E,$$
then $E\geq 0$ and is exceptional$/X'$.
\end{lem}
\begin{proof}
    It follows from the same lines of the proof of Lemma \ref{lem: hl23 2.6}.
\end{proof}

\begin{lem}[Lemma \ref{lem: g-pair version bir12 2.7}]\label{lem: g-pair version bir12 2.7 gfq}
Let $(X,\Ff,B,\Mm)/U$ be a generalized foliated quadruple. Let $(X_1,\Ff_1,B_1,\Mm)/U$ and $(X_2,\Ff_2,B_2,\Mm)/U$ be two bs-weak lc models of $(X,\Ff,B,\Mm)/U$ with induced birational maps $\phi: X_1\dashrightarrow X_2$. Let $h_1: W\rightarrow X_1$ and $h_2: W\rightarrow X_2$ be two birational morphisms such that $\phi\circ h_1=h_2$. Then:
\begin{enumerate}
    \item $$h_1^*(K_{\Ff_1}+B_1+\Mm_{X_1})=h_2^*(K_{\Ff_2}+B_2+\Mm_{X_2}).$$
    \item If $K_{\Ff_2}+B_2+\Mm_{X_2}$ is semi-ample$/U$, then $K_{\Ff_1}+B_1+\Mm_{X_1}$ is semi-ample$/U$.
    \item If $K_{\Ff_2}+B_2+\Mm_{X_2}$ is ample$/U$, then $\phi$ is a morphism.
\end{enumerate}
\end{lem}
\begin{proof}
    It follows from the same lines of the proof of Lemma \ref{lem: g-pair version bir12 2.7} except we replace  Lemma \ref{lem: hl23 2.6} with  Lemma \ref{lem: hl23 2.6 gfq}.
\end{proof}

\begin{lem}[Lemma \ref{lem: numerical equivalence model}]\label{lem: numerical equivalence model gfq}
    Let $r$ be a positive real number. Let $(X,\Ff_1,B_1,\Mm_1)/U$ and $(X,\Ff_2,B_2,\Mm_2)/U$ be two generalized foliated quadruples such that
    $$K_{\Ff_2}+B_2+\Mm_{2,X}\equiv_U r(K_{\Ff_1}+B_1+\Mm_{1,X})$$
    Let $(X',\Ff'_1,B_1',\Mm_1)/U$ be a weak lc model (resp. minimal model) of $(X,\Ff_1,B_1,\Mm_1)/U$ with induced birational map $\phi: X\dashrightarrow X'$. Let $\Ff_2':=\phi_*\Ff$ and $B_2':=\phi_*B_2$. Then $(X',\Ff_2',B_2',\Mm_2)/U$ is a weak lc model (resp. minimal model) of $(X,\Ff_2,B_2,\Mm_2)/U$.

    If  $(X',\Ff'_1,B_1',\Mm_1)/U$ is a semi-ample model (resp. good minimal model) of $(X,\Ff_1,B_1,\Mm_1)/U$ and
    $$K_{\Ff_2}+B_2+\Mm_{2,X}\sim_{\mathbb R,U} r(K_{\Ff_1}+B_1+\Mm_{1,X}),$$
        $(X',\Ff'_2,B_2',\Mm_2)/U$ is a semi-ample model (resp. good minimal model) of $(X,\Ff_2,B_2,\Mm_2)/U$.
\end{lem}
\begin{proof}
    It follows from the same lines of the proof of Lemma \ref{lem: numerical equivalence model} except we replace  Lemma \ref{lem: hl23 2.6} with  Lemma \ref{lem: hl23 2.6 gfq}.
\end{proof}

\begin{thm}[Theorem \ref{thm: chlx23 9.4.1}]\label{thm: chlx23 9.4.1 gfq}
    Let $(X,\Ff,B,\Mm)/U$ be a $\Qq$-factorial ACSS algebraically integrable generalized foliated quadruple such that $K_{\Ff}+B+\Mm_X\sim_{\mathbb R,U}E\geq 0$ and $E$ is very exceptional$/U$. Then we may run a $(K_{\Ff}+B+\Mm_X)$-MMP$/U$ with scaling of an ample$/U$ $\Rr$-divisor $A$ and any such MMP terminates with a good log minimal model $(X',\Ff',B',\Mm)/U$ such that $K_{\Ff'}+B'+\Mm_X\sim_{\mathbb R,U}0$.
\end{thm}
\begin{proof}
    It follows from \cite[Theorem 9.4.1]{CHLX23}.
\end{proof}

\begin{lem}[Lemma \ref{lem: g-pair version bir12 2.8}]\label{lem: g-pair version bir12 2.8 gfq}
Let $(X,\Ff,B,\Mm)/U$ be an lc algebraically integrable generalized foliated quadruple. Let $(W,\Ff_W,B_W,\Mm)$ be a foliated log smooth model of $(X,\Ff,B,\Mm)$. 

Then any bs-weak lc model (resp. bs-minimal model, bs-semi-ample model, bs-good minimal model, log minimal model, good log minimal model) of $(W,\Ff_W,B_W,\Mm)/U$ is a bs-weak lc model (resp. bs-minimal model, bs-semi-ample model, bs-good minimal model, log minimal model, good log minimal model) of $(X,\Ff,B,\Mm)/U$. 
\end{lem}
\begin{proof}
    It follows from the same lines of the proof of Lemma \ref{lem: g-pair version bir12 2.8} except we replace  Lemma \ref{lem: hl23 2.6} with  Lemma \ref{lem: hl23 2.6 gfq}.
\end{proof}

\begin{lem}[Lemma \ref{lem: foliation lsm has lmm}]\label{lem: foliation lsm has lmm gfq}
Let $(X,\Ff,B,\Mm)/U$ be an lc algebraically integrable generalized foliated quadruple and $(X',\Ff',B',\Mm)/U$ a bs-weak lc model of $(X,\Ff,B,\Mm)/U$. Let $(W,\Ff_W,B_W,\Mm)$ be a foliated log smooth model of $(X,\Ff,B,\Mm)$ such that the induced birational map $\phi_W: W\dashrightarrow X'$ is a morphism. 

Then we may run a $(K_{\Ff_W}+B_W+\Mm_W)$-MMP$/X'$ with scaling of an ample$/X'$ $\Rr$-divisor which terminates with a good minimal model $(Y,\Ff_Y,B_Y,\Mm)/X'$ of $(W,\Ff_W,B_W,\Mm)/X'$ such that $$K_{\Ff_Y}+B_Y+\Mm_Y=q^*(K_{\Ff'}+B'+\Mm_{X'}).$$
where $q: Y\rightarrow X'$ is the induced morphism. In particular, $(Y,\Ff_Y,B_Y,\Mm)/U$ is a log minimal model of $(W,\Ff_W,B_W,\Mm)/U$.
\end{lem}
\begin{proof}
    It follows from the same lines of the proof of Lemma \ref{lem: foliation lsm has lmm} except we replace Lemma \ref{lem: hl23 2.6} with  Lemma \ref{lem: hl23 2.6 gfq} and replace Theorem \ref{thm: chlx23 9.4.1} with Theorem \ref{thm: chlx23 9.4.1 gfq}.
\end{proof}

\begin{lem}[Lemma \ref{lem: g-pair weak glc imply lmm}]\label{lem: g-pair weak glc imply lmm gfq}
Let $(X,\Ff,B,\Mm)/U$ be an lc algebraically integrable foliated quadruple. If $(X,\Ff,B,\Mm)/U$ has a bs-weak lc model (resp. bs-semi-ample model), then $(X,\Ff,B,\Mm)/U$ has a log minimal model (resp. good log minimal model).
\end{lem}
\begin{proof}
By Lemma \ref{lem: g-pair version bir12 2.7 gfq} we only need to prove the bs-weak lc model case. The lemma follows immediately from Lemmas \ref{lem: g-pair version bir12 2.8 gfq} and \ref{lem: foliation lsm has lmm gfq}.
\end{proof}

\begin{lem}[Lemma \ref{lem: same weak glc model under pullback}]\label{lem: same weak glc model under pullback gfq}
Let $(X,\Ff,B,\Mm)/U$ and $(Y,\Ff_Y,B_Y,\Mm)/U$ be two lc algebraically integrable generalized foliated quadruples, and let $f: Y\rightarrow X$ be a birational morphism such that
$$K_{\Ff_Y}+B_Y+\Mm_Y=f^*(K_\Ff+B+\Mm_X)+E$$
for some $E\geq 0$ that is exceptional$/X$ and $f_*\Ff_Y=\Ff$. Then:
\begin{enumerate}
    \item Any bs-weak lc model of $(X,\Ff,B,\Mm)/U$ is a bs-weak lc model of $(Y,\Ff_Y,B_Y,\Mm)/U$.
    \item If $(X,\Ff,B,\Mm)/U$ has a bs-weak lc model (resp. bs-semi-ample model), then
    
    $(Y,\Ff_Y,B_Y,\Mm)/U$ has a log minimal model (resp. good log minimal model).
\end{enumerate}
\end{lem}
\begin{proof}
    It follows from the same lines of the proof of Lemma \ref{lem: same weak glc model under pullback} except that we replace Lemmas \ref{lem: hl23 2.6}, \ref{lem: g-pair version bir12 2.7}, and \ref{lem: g-pair weak glc imply lmm} with Lemmas \ref{lem: hl23 2.6 gfq}, \ref{lem: g-pair version bir12 2.7 gfq}, and \ref{lem: g-pair weak glc imply lmm gfq} respectively.
\end{proof}

\begin{lem}[Lemma \ref{lem: chlx 9.2.1}]\label{lem: chlx 9.2.1 gfq}
    Let $(X,\Ff,B,\Mm)/U$ be an lc algebraically integrable generalized foliated quadruple, $G$ a reduced divisor on $X$, and $f: X\rightarrow Z$ a contraction, such that $(X,\Ff,B,\Mm;G)/Z$ satisfies Property $(*)$ and $K_{\Ff}+B+\Mm_X\sim_{\mathbb R,U}K_X+B+G+\Mm_X$. Assume that $G$ is super$/Z$. Let $D\geq 0$ be an $\Rr$-divisor on $X$ such that $K_{\Ff}+B+D+\Mm_X$ is nef$/U$.
    
    Then any sequence of steps of a $(K_{\Ff}+B+\Mm_X)$-MMP$/U$ (with scaling of $D$) is a sequence of steps of a $(K_{X}+B+G+\Mm_X)$-MMP$/U$ (with scaling of $D$), and any sequence of steps of a $(K_{X}+B+G+\Mm_X)$-MMP$/U$ (with scaling of $D$) is a sequence of steps of a  $(K_{\Ff}+B+\Mm_X)$-MMP$/U$ (with scaling of $D$). Moreover, any sequence of steps of a $(K_{\Ff}+B+\Mm_X)$-MMP$/U$ or a 
 $(K_{X}+B+G+\Mm_X)$-MMP$/U$ is a sequence of steps of an MMP$/Z$.
\end{lem}
\begin{proof}
    It follows from \cite[Lemma 9.2.1]{CHLX23}.
\end{proof}

We remark that the condition ``NQC" is needed for the next theorem.

\begin{thm}[Theorem \ref{thm: hh20 1.7}]\label{thm: hh20 1.7 gfq}
    Let $(X,B,\Mm)/U$ be an NQC lc generalized pair and $A$ an ample$/U$ $\Rr$-divisor on $X$ such that $(X,B+A,\Mm)$ is lc and $K_X+B+A+\Mm_X$ is nef$/U$. Assume that $(X,B,\Mm)/U$ has a $\Qq$-factorial bs-minimal model or $K_X+B+\Mm_X$ is not pseudo-effective$/U$. Then there exists a sequence of $(K_X+B+\Mm_X)$-MMP$/U$ with scaling of $A$ which terminates with either a minimal model or a Mori fiber space of $(X,B)/U$. 
\end{thm}
\begin{proof}
It follows from \cite[Theorem A]{TX24}.
\end{proof}

\begin{lem}[Lemma \ref{lem: minimal model foliation over zu}]\label{lem: minimal model foliation over zu gfq}
    Let $(X,\Ff,B,\Mm)/U$ be an lc algebraically integrable generalized foliated quadruple. Assume that the associated morphism $\pi: X\rightarrow U$ is a contraction, and assume that $\Ff$ is induced by a contraction $f: X\rightarrow Z$. Let $Z_U$ be the core model of $(\pi,f)$. Then:
    \begin{enumerate}
        \item Any sequence of steps of a $(K_{\Ff}+B+\Mm_X)$-MMP$/U$ is a step of a $(K_{\Ff}+B+\Mm_X)$-MMP$/Z_U$.
        \item If $(X,\Ff,B,\Mm)$ is $\Qq$-factorial ACSS and $K_{\Ff}+B+\Mm_X$ is nef$/U$, then $K_{\Ff}+B+\Mm_X$ is nef$/Z_U$.
        \item $(X,\Ff,B,\Mm)/U$ has a bs-weak lc model if and only if $(X,\Ff,B,\Mm)/Z_U$ has a bs-weak lc model.
    \end{enumerate}
\end{lem}
\begin{proof}
    It follows from the same lines of the proof of Lemma \ref{lem: minimal model foliation over zu} except that we replace Lemmas \ref{lem: g-pair version bir12 2.8}, \ref{lem: foliation lsm has lmm}, \ref{lem: g-pair weak glc imply lmm} with Lemmas \ref{lem: g-pair version bir12 2.8 gfq}, \ref{lem: foliation lsm has lmm gfq}, \ref{lem: g-pair weak glc imply lmm gfq} respectively.
\end{proof}

\begin{lem}[Lemma \ref{lem: super minimal model over zu}]\label{lem: super minimal model over zu gfq}
    Let $(X,B,\Mm)/U$ be a generalized pair associated with contraction $\pi: X\rightarrow U$. Let $f: X\rightarrow Z$ be a contraction such that $B$ is super$/Z$. Let $Z_U$ be the core model of $(\pi,f)$. Then:
    \begin{enumerate}
    \item If $K_X+B+\Mm_X$ is nef$/Z_U$ then $K_X+B+\Mm_X$ is nef$/U$.
    \item Any sequence of steps of a $(K_X+B+\Mm_X)$-MMP$/U$ is a sequence of steps of a $(K_X+B+\Mm_X)$-MMP$/Z_U$.
    \item If $(X,B,\Mm)/Z_U$ has a minimal model then $(X,B,\Mm)/U$ has a minimal model.
    \item If $(X,B,\Mm)/U$ has a minimal model and $\Mm$ is NQC$/U$, then $(X,B,\Mm)/Z_U$ has a minimal model.
    \end{enumerate}
\end{lem}
\begin{proof}
It follow from the same lines of the proof of  Lemma \ref{lem: super minimal model over zu} except that the length of extremal ray control for pair is replaced by \cite[Theorem 2.2.1(2)]{CHLX23} for generalized pairs, and Lemma \ref{lem: g-pair weak glc imply lmm}, Theorem \ref{thm: hh20 1.7} are replaced with \ref{lem: g-pair weak glc imply lmm gfq} and Theorem \ref{thm: hh20 1.7 gfq} respectively.
\end{proof}

\begin{prop}[Proposition \ref{prop: eolmm foliation to pair}]\label{prop: eolmm foliation to pair gfq}
    Let $(X,\Ff,B,\Mm)/U$ be an lc algebraically integrable generalized foliated quadruple. Assume that $(X,\Ff,B,\Mm)/U$ has a bs-weak lc model. Then there exists an ACSS modification $h: (X',\Ff',B',\Mm;G)/Z\rightarrow (X,\Ff,B,\Mm)$ that is $\Qq$-factorial, strict, and super, and $(X',B'+G,\Mm)/U$ has a log minimal model.
\end{prop}
\begin{proof}
    It follows from the same lines of the proof of Proposition \ref{prop: eolmm foliation to pair gfq} will the following modifications: Lemmas \ref{lem: numerical equivalence model}, \ref{lem: g-pair version bir12 2.8}, \ref{lem: foliation lsm has lmm}, \ref{lem: g-pair weak glc imply lmm}, \ref{lem: same weak glc model under pullback}, \ref{lem: chlx 9.2.1}, \ref{lem: minimal model foliation over zu}, and Theorem \ref{thm: chlx23 9.4.1} are replaced by Lemmas \ref{lem: numerical equivalence model gfq}, \ref{lem: g-pair version bir12 2.8 gfq}, \ref{lem: foliation lsm has lmm gfq}, \ref{lem: g-pair weak glc imply lmm gfq},  \ref{lem: same weak glc model under pullback gfq}, \ref{lem: chlx 9.2.1 gfq}, \ref{lem: minimal model foliation over zu gfq}, and Theorem \ref{thm: chlx23 9.4.1 gfq} respectively, and Lemma \ref{lem: super minimal model over zu} is replaced with Lemma \ref{lem: super minimal model over zu gfq}(3).
\end{proof}

\begin{lem}[Lemma \ref{lem: TX24 2.20}]\label{lem: TX24 2.20 gfq}
    Let $(X,B+A,\Mm)/U$ be an NQC lc generalized pair such that $(X,B,\Mm)$ is lc and $K_X+B+A+\Mm_X$ is NQC$/U$. Then there exists a positive real number $\epsilon\in (0,1)$ such that any $(K_X+B+(1-\epsilon)A+\Mm_X)$-MMP$/U$ is $(K_X+B+A+\Mm_X)$-trivial for any $\epsilon\in (0,\epsilon_0)$.
\end{lem}
\begin{proof}
    It follows from \cite[Lemma 2.20]{TX24}.
\end{proof}

\begin{thm}[Theorem \ref{thm: bir12 1.9(3)}]\label{thm: bir12 1.9(3) gfq}
     Let $(X,B,\Mm)/U$ be a $\Qq$-factorial NQC lc generalized pair and $H\geq 0$ an $\Rr$-divisor on $X$ such that $K_X+B+H+\Mm_X$ is nef$/U$ and $(X,B+H,\Mm)$ is lc. Assume that $X$ is klt, and there exists an infinite sequence of $(K_X+B+\Mm_X)$-MMP$/U$ with scaling of $H$ with scaling numbers $\lambda_i$ such that $\lim_{i\rightarrow+\infty}\lambda_i=\lambda$ and $\lambda\not=\lambda_i$ for any $i$.
     
     Then $(X,B+\lambda H,\Mm)/U$ does not have a bs-minimal model.
\end{thm}
\begin{proof}
    By \cite[Theorem 4.1]{HL22}, $(X,B+\lambda H,\Mm)/U$ does not have a log minimal model that is $\Qq$-factorial dlt. By \cite[Lemma 3.8]{HL23}, $(X,B+\lambda H,\Mm)/U$ does not have a  bs-minimal model.
\end{proof}

\begin{lem}[Lemma \ref{lem: gmmp scaling numbers go to 0}]\label{lem: gmmp scaling numbers go to 0 gfq}
Let $(X,B,\Mm)/U$ be a $\Qq$-factorial NQC lc generalized pair such that $X$ is klt. Let $H\geq 0$ be an $\Rr$-divisor on $X$ such that $(X,B+H,\Mm)$ is lc and $K_X+B+H+\Mm_X$ is nef$/U$. Assume that for any $\mu\in [0,1]$,
\begin{itemize}
    \item either $(X,B+\mu H,\Mm)/U$ has a log minimal model, or
    \item $K_X+B+\mu H+\Mm_X$ is not pseudo-effective$/U$.
\end{itemize}
Then there exists a $(K_X+B+\Mm_X)$-MMP$/U$ with scaling of $H$ which terminates after finitely many steps.
\end{lem}
\begin{proof}
  It follows from the same lines of the proof of Lemma \ref{lem: gmmp scaling numbers go to 0} except that we replace Lemma \ref{lem: TX24 2.20}, Theorem \ref{thm: hh20 1.7},  Theorem \ref{thm: bir12 1.9(3)} with Lemma \ref{lem: TX24 2.20 gfq}, Theorem \ref{thm: hh20 1.7 gfq},  Theorem \ref{thm: bir12 1.9(3) gfq} respectively.
\end{proof}

\begin{thm}[Theorem \ref{thm: hh20 1.5}]\label{thm: hh20 1.5 gfq}
    Let $(X,B,\Mm)/U$ be an NQC lc generalized pair and $A$ an ample$/U$ $\Rr$-divisor on $X$ such that $(X,B+A,\Mm)$ is lc. Then $(X,B+A,\Mm)/U$ has a bs-good minimal model or a bs-Mori fiber space.
\end{thm}
\begin{proof}
    It follows from \cite[Theorems A,F]{TX24}.
\end{proof}

\begin{thm}[Theorem \ref{thm: take strict simple model run mmp}]\label{thm: take strict simple model run mmp gfq}
    Let $(X,\Ff,B,\Mm)/U$ be an lc algebraically integrable generalized foliated quadruple and let $A,H$ be two ample$/U$ $\Rr$-divisors on $X$. Let $h: (X',\Ff',B,\Mm;G)/Z\rightarrow (X,\Ff,B,\Mm)$ be a simple model of $(X,\Ff,B,\Mm)$ that is strict and super, $H':=h^*H$, and $A':=h^*A$. Assume that
    \begin{itemize}
        \item either $X$ is potentially klt, or
        \item $X'$ is $\Qq$-factorial klt and $\Mm$ is NQC$/U$.
    \end{itemize}
    Then:
    \begin{enumerate}
        \item We may run a $(K_{\Ff'}+B'+H')$-MMP$/U$ with scaling of $A'$, say $\mathcal{P}$, such that $\mathcal{P}$ terminates with either a minimal model or a Mori fiber space of $(X',\Ff',B'+H')/U$.
        \item If $X$ is potentially klt, then $\mathcal{P}$ can be any $(K_{\Ff'}+B'+H')$-MMP$/U$ with scaling of $A'$.
    \end{enumerate}
\end{thm}
\begin{proof}
     It follows from the same lines of the proof of Theorem \ref{thm: take strict simple model run mmp} except the following differences: we replace Lemmas \ref{lem: gklt is klt}, \ref{lem: existence of core model}, \ref{lem: existence of klt pair on strict simple models}, \ref{lem: numerical equivalence model}, \ref{lem: g-pair weak glc imply lmm}, \ref{lem: same weak glc model under pullback}, \ref{lem: chlx 9.2.1}, \ref{lem: gmmp scaling numbers go to 0} and Theorem \ref{thm: hh20 1.5} with  Lemmas \ref{lem: gklt is klt gfq}, \ref{lem: existence of core model gfq}, \ref{lem: existence of klt pair on strict simple models gfq}, \ref{lem: numerical equivalence model gfq}, \ref{lem: g-pair weak glc imply lmm gfq}, \ref{lem: same weak glc model under pullback gfq}, \ref{lem: chlx 9.2.1 gfq}, \ref{lem: gmmp scaling numbers go to 0 gfq} and Theorem \ref{thm: hh20 1.5 gfq} respectively. We remark that  Lemma \ref{lem: gklt is klt gfq} may need to be applied in \textbf{Step 1} again in order to get the boundary $\Delta'$.
\end{proof}

\begin{prop}[Proposition \ref{prop: contraction not big}]\label{prop: contraction not big gfq}
    Let $(X,\Ff,B,\Mm)/U$ be an lc algebraically integrable generalized foliated quadruple such that $(X,\Delta,\Nn)/U$ is klt, where $B\geq\Delta\geq 0$ and $\Mm-\Nn$ is nef$/U$. Let $R$ be a $(K_{\Ff}+B+\Mm_X)$-negative extremal ray$/U$ and $H_R$ a supporting function$/U$ of $R$. Suppose that $H_R$ is not big$/U$. Then $R$ is also a $(K_X+\Delta+\Nn_X)$-negative extremal ray$/U$. In particular, there exists a contraction $\cont_R$ of $R$.
\end{prop}
\begin{proof}
    It follows from the same lines of the proof of Proposition \ref{prop: contraction not big} except that we replace Proposition \ref{prop: weak cbf gfq} and Theorem \ref{thm: eo acss model} with Proposition \ref{prop: weak cbf gfq gfq} and Theorem \ref{thm: eo acss model gfq} respectively, and use \cite[Theorem 2.2.1(4)]{CHLX23} instead of the contraction theorem for lc pairs.
\end{proof}

\begin{lem}\label{lem: klt bpf}
    Let $(X,B,\Mm)/U$ be a klt generalized pair and $L$ a nef$/U$ $\Rr$-divisor on $X$ such that $aL-(K_X+B+\Mm_X)$ is big$/U$ and nef$/U$ for some $a>0$. Then $L$ is semi-ample$/U$. Moreover, if $L$ is Cartier, then $\mathcal{O}_X(mL)$ is globally generated$/U$ for any $m\gg 0$.
\end{lem}
\begin{proof}
We have $aL-(K_X+B+\Mm_X)=L_n+\frac{1}{n}E$ for some ample$/U$ $\Rr$-divisor $L_n$ and $E\geq 0$. Let $n\gg 0$ be an integer, then $(X,\Delta:=B+\frac{1}{n}E,\Mm)$ is klt. and $aL-(K_X+\Delta+\Mm_X)$ is ample$/U$. By \cite[Theorem 2.2.7]{CHLX23}, $L$ is semi-ample$/U$.

Assume that $L$ is Cartier. Let $p>a,q>a$ be two different prime numbers. Then $pL-(K_X+\Delta+\Mm_X)$ and $qL-(K_X+\Delta+\Mm_X)$ are ample$/U$. By \cite[Theorems 2.2.6]{CHLX23}, $\mathcal{O}_X(p^rL)$ and $\mathcal{O}_X(q^sL)$ are globally generated$/U$ for some positive integers $r,s$. For any integer $m\gg0$, $m=bp^r+cq^s$ for some non-negative integers $b,c$. Thus $\mathcal{O}_X(mL)$ is globally generated$/U$ for any $m\gg 0$.
\end{proof}

\begin{lem}[Lemma \ref{lem: can run mmp scaling}]\label{lem: can run mmp scaling gfq}
        Let $(X,\Ff,B,\Mm)/U$ be an lc algebraically integrable generalized foliated quadruple and $A$ an ample$/U$ $\Rr$-divisor on $X$. Let $\mathcal{P}:$
    $$(X_0,\Ff_0,B_0,\Mm):=(X,\Ff,B,\Mm)\dashrightarrow (X_1,\Ff_1,B_1,\Mm)\dashrightarrow\dots\dashrightarrow (X_n,\Ff_n,B_n,\Mm)$$
    be a sequence of steps of a $(K_{\Ff}+B+\Mm_X)$-MMP$/U$ with scaling of $A$ and let $A_i$ be the image of $A$ on $X_i$ for each $i$. Let
    $$\lambda_n:=\inf\{t\geq 0\mid K_{\Ff_n}+B_n+tA_n+\Mm_X\text{ is nef}/U\}.$$
    Suppose that $\lambda_n>0$. Then there exists a $(K_{\Ff_n}+B_n+\Mm_{X_n})$-negative extremal ray$/U$ $R$ such that $(K_{\Ff_n}+B_n+\lambda_nA_n+\Mm_{X_n})\cdot R=0$.
\end{lem}
\begin{proof}
    It follows from the same lines of the proof of Lemma \ref{lem: can run mmp scaling}.
\end{proof}

\begin{prop}[Proposition \ref{prop: lift mmp}]\label{prop: lift mmp gfq}
    Let $(X,\Ff,B,\Mm)/U$ be an lc algebraically integrable generalized foliated quadruple such that either $\Mm$ is NQC$/U$, or $X$ is potentially klt. Let $\mathcal{P}:$
$$(X_0,\Ff_0,B_0,\Mm):=(X,\Ff,B,\Mm)\dashrightarrow (X_1,\Ff_1,B_1,\Mm)\dashrightarrow\dots\dashrightarrow (X_n,\Ff_n,B_n,\Mm)\dashrightarrow\dots$$
be a (possibly infinite) sequence of $(K_{\Ff}+B+\Mm_X)$-MMP$/U$. For each $i\geq 0$, we let $\psi_i: X_i\rightarrow T_i$ and $\psi_i^+:X_{i+1}\rightarrow T_{i}$ be the $(i+1)$-th step of this MMP and let $\phi_i:=(\psi_{i}^+)^{-1}\circ\psi_i: X_i\dashrightarrow X_{i+1}$ be the induced birational map. Let $h: (Y,\Ff_Y,B_Y,\Mm;G)/Z\rightarrow (X,\Ff,B,\Mm)$ be an ACSS modification of $(X,\Ff,B,\Mm)$ that is $\Qq$-factorial, strict, and super. Let $A$ be an ample$/U$ $\Rr$-divisor on $X$ and let $A_i$ be the image of $A$ on $X_i$ for each $i$.

Then there exist a (possibly infinite) sequence $\mathcal{P}_Y$ of birational maps 
$$(Y_0,\Ff_{Y_0},B_{Y_0},\Mm):=(Y,\Ff_Y,B_Y,\Mm)\dashrightarrow (Y_1,\Ff_{Y_1},B_{Y_1},\Mm)\dashrightarrow\dots\dashrightarrow (Y_n,\Ff_{Y_n},B_{Y_n},\Mm)\dashrightarrow\dots$$
satisfying the following. Let $\phi_{i,Y}: Y_i\dashrightarrow Y_{i+1}$ be the induced birational map. Then:
\begin{enumerate}
\item For any $i\geq 0$, there exist an ACSS modification $h_i: (Y_i,\Ff_{Y_i},B_{Y_i},\Mm;G_i)/Z\rightarrow (X_i,\Ff_i,B_i)$ that is $\Qq$-factorial, strict, and super, such that $h_0=h$ and $G_i$ is the image of $G$ on $Y_i$.
\item For any $i\geq 0$, $h_{i+1}\circ\phi_{i,Y}=\phi_i\circ h_i$.
\item For any $i\geq 0$, $\phi_{i,Y}$ is a $(K_{\Ff_i}+B_{Y_i}+\Mm_{Y_i})$-MMP$/T_i$ and $(Y_{i+1},\Ff_{Y_{i+1}},B_{Y_{i+1}},\Mm)/T_i$ is the output of this MMP, such that $\phi_{i,Y}$ is not the identity map.
\item $\mathcal{P}_Y$ is a sequence of steps of a $(K_{\Ff_Y}+B_Y+\Mm_Y)$-MMP$/U$.
\item Suppose that $\mathcal{P}$ is an MMP$/U$ with scaling of $A$. Let $A_{Y}:=h^*A$ and let $A_{Y_i}$ the image of $A_Y$ on $Y_i$ for each $i$. Let
$$\lambda_i:=\inf\{t\geq 0\mid K_{\Ff_i}+B_i+tA_i+\Mm_{X_i}\text{ is nef}/U\}$$
be the $(i+1)$-th scaling number. Then:
\begin{enumerate}
    \item $\phi_{i,Y}$ is a sequence of steps of a $(K_{\Ff_{Y_i}}+B_{Y_i}+\Mm_{Y_i})$-MMP$/U$ with scaling of $A_{Y_i}$, and the scaling number of each step of $\phi_{i,Y}$ is $\lambda_i$.
    \item $\mathcal{P}_Y$ is sequence of steps of a $(K_{\Ff_Y}+B_Y+\Mm_Y)$-MMP$/U$ with scaling of $A_Y$.
\end{enumerate}
\end{enumerate}
\end{prop}
\begin{proof}
    It follows from the same lines of the proof of Proposition \ref{prop: lift mmp} except that we replace Theorem \ref{thm: take strict simple model run mmp} with Theorem \ref{thm: take strict simple model run mmp gfq}.
\end{proof}

\section{Nef divisors with real coefficients}\label{sec: enef}

In \textbf{Step 3} of the proof of Theorem \ref{thm: cont and flip with detail}, we use Shokurov-type polytopes to construct an MMP $\varphi: X_n\dashrightarrow\bar X$ which is $(K_{\Ff_n}+B_n+\Mm_{X_n}+A_n)$-trivial, which is no longer valid if $\Mm$ is not NQC$/U$. This prevents us from proving the (possibly non-NQC) generalized foliated quadruple variations of our main theorems. We still want to prove these variations, not only for completeness but also for potential applications to the minimal model program on K\"ahler varieties (cf. \cite{DH23,DHY23}). In this appendix, we prove a result, Proposition \ref{prop: kf+k_x mmp trivial}, as a substitution of Shokurov-type polytopes for non-NQC generalized foliated quadruples, which is applied to the proof of Theorem \ref{thm: cont and flip with detail} in \textbf{Step 3}. The key idea for our proof is to introduce the notation of \emph{$\epsilon$-nef} $\Rr$-divisors and study its behavior.

\begin{defn}
Let $\epsilon$ be a positive real number and $\Ii$ a set of positive real numbers. Let $\pi: X\rightarrow U$ be a projective morphism from a normal quasi-projective variety to a variety and let $D$ be a nef$/U$ $\Rr$-divisor on $X$.

\begin{enumerate}
    \item ($\epsilon$-nef) We say that $D$ is \emph{$\epsilon$-nef}$/U$ if
    $$D\cdot C\geq\epsilon$$
    for any curve $C$ on $X$ such that 
    \begin{itemize}
        \item $D\cdot C>0$, and
        \item $C$ spans an extremal ray in $\overline{NE}(X/U)$.
    \end{itemize}
    \item ($\Ii$-NQC) We say that $D$ is \emph{$\Ii$-NQC}$/U$ if we can write $D=\sum a_iD_i$, such that each $a_i\in\Ii$ and each $D_i$ is nef$/U$ Cartier. In addition, if each $a_i\geq\epsilon$, then we say that $D$ is \emph{$\epsilon$-NQC$/U$}.
\end{enumerate}
For any nef$/U$ $\bb$-divisor $\Mm$ on $X$, we say that $\Mm$ is \emph{$\epsilon$-nef}$/U$ (resp. $\Ii$-NQC$/U$, $\epsilon$-NQC$/U$) if there exists a birational morphism $Y\rightarrow X$, such that $\Mm$ descends to $Y$ and $\Mm_Y$ is $\epsilon$-nef$/U$ (resp. $\Ii$-NQC$/U$, $\epsilon$-NQC$/U$).

The following results are clear and we are free to use it in the rest of this appendix.
\begin{itemize}
\item  Any $\epsilon$-NQC$/U$ $\Rr$-divisor and $\bb$-divisor is $\epsilon$-nef$/U$.
\item If $\min\{\gamma\in\Ii\}\geq\epsilon$, then any $\Ii$-NQC$/U$ $\Rr$-divisor and $\bb$-divisor is $\epsilon$-NQC$/U$ and $\epsilon$-nef$/U$.
\item Any NQC$/U$ $\Rr$-divisor and $\bb$-divisor is $\epsilon$-NQC$/U$ for some positive real number $\epsilon$.
\end{itemize}
\end{defn}

\begin{lem}\label{lem: nqc good decomposition}
    Let $X\rightarrow U$ be a projective morphism from a normal quasi-projective variety to a variety and let $\Mm$ be an NQC$/U$ $\bb$-divisor on $X$. Then there exist positive real numbers $a_1,\dots,a_k$ that are linearly independent over $\mathbb Q$ such that $\Mm$ is $\{a_1,\dots,a_k\}$-NQC$/U$.
\end{lem}
\begin{proof}
    We may write $\Mm=\sum_{i=1}^mv_i\Nn_i$ where each $\Nn_i$ is a nef$/U$ $\bb$-Cartier $\bb$-divisor and each $v_i>0$. Let $\bm{v}_0:=(v_{0,1},\dots,v_{0,m})$ and let $V$ be the rational polytope of $\bm{v}_0$ in $\mathbb R^m$. Let $\Mm(\bm{v}):=\sum_{i=1}^mv_i\Nn_i$ for any $\bm{v}=(v_1,\dots,v_m)$ in $\mathbb R^m$. 
    
    Suppose that $\dim V=n$, then we may take rational points $\bm{v}_1,\dots,\bm{v}_{n+1}$ in $V$ such that $\bm{v}$ is contained in the interior of the convex hull of $\bm{v}_1,\dots,\bm{v}_{n+1}$. Then there exist unique positive real numbers $b_1,\dots b_{n+1}$ such that $\sum_{i=1}^{n+1}b_i=1$ and $\sum_{i=1}^{n+1}b_i\bm{v}_i=\bm{v}_0$. Moreover, $b_1,\dots,b_{n+1}$ are linearly independent over $\mathbb Q$. Let $N$ be a positive integer such that $N\Mm(\bm{v}_i)$ is Cartier for each $i$. Then we have
    $$\Mm=\sum_{i=1}^{n+1}\frac{b_i}{N}(N\Mm(\bm{v}_i)).$$
    We may take $k:=n+1$ and $a_i:=\frac{b_i}{N}$ for each $i$. 
\end{proof}

\begin{lem}\label{lem: inf=min length}
    Let $\pi: X\rightarrow U$ be a projective morphism from a normal quasi-projective variety to a variety. Let $C$ be a curve on $X$ such that $\pi(C)$ is a point. Let
    $$\lambda_0:=\inf\{\lambda>0\mid \text{there exists a curve }C'\equiv_U\lambda C, \pi(C')=\{pt\}\}.$$
    Then
    $$\lambda_0:=\min\{\lambda>0\mid \text{there exists a curve }C'\equiv_U\lambda C, \pi(C')=\{pt\}\}.$$
\end{lem}
\begin{proof}
    Suppose not, then there exists a sequence of curves $C_i$ on $X$ and a strictly decreasing sequence of real numbers $\{\lambda_i\}_{i=1}^{+\infty}$, such that $\lim_{i\rightarrow+\infty}\lambda_i=\lambda_0$, $\pi(C_i)=\{pt\}$, and $C_i\equiv\lambda_iC$ for each $i$. Let $H$ be an ample Cartier divisor on $X$, then 
    $$\lambda_i=\frac{H_i\cdot C}{H\cdot C},$$
    so $(H\cdot C)\lambda_i$ is an integer for each $i$, which is not possible.
\end{proof}

\begin{lem}\label{lem: weak nqc only check smallest}
      Let $\pi: X\rightarrow U$ be a projective morphism from a normal quasi-projective variety to a variety and let $D$ be a nef$/U$ $\Rr$-divisor on $X$. Let $\epsilon$ be a positive real number. Suppose that
    $$D\cdot C\geq\epsilon$$
        for any curve $C$ on $X$ such that 
        \begin{itemize}
            \item $D\cdot C>0$,
            \item  $C$ spans an extremal ray in $\overline{NE}(X/U)$, and
            \item for any curve $C'$ on $X$ so that $C'\equiv_U\lambda C$ for some real number $\lambda>0$, we have $\lambda\geq 1$.
        \end{itemize}
        Then $D$ is $\epsilon$-nef$/U$.
\end{lem}
\begin{proof}
Suppose that $D$ is not $\epsilon$-nef$/U$. Then there exists a curve $C'$ on $X$ such that $D\cdot C'>0$, $C'$ spans an extremal ray $R$ in $\overline{NE}(X/U)$, and $D\cdot C'<\epsilon$. We let
$$\lambda_0:=\inf\{\lambda>0\mid \text{there exists a curve }C''\equiv_U\lambda C', \pi(C')=\{pt\}\}.$$
Then $\lambda_0\leq 1$. By Lemma \ref{lem: inf=min length}, there exists a curve $C_0$ on $X$ such that $C_0\equiv_U\lambda_0C'$ and $\pi(C_0)=\{pt\}$. Then $C_0$ spans $R$, and for any curve $C''$ on $X$ so that $C''\equiv_U\lambda C_0$ for some real number $\lambda>0$, we have $\lambda\geq 1$. By our assumption, $D\cdot C_0\geq\epsilon$. Therefore,
$$D\cdot C'=\frac{1}{\lambda_0}D\cdot C_0\geq\frac{\epsilon}{\lambda_0}\geq\epsilon,$$
a contradiction.
\end{proof}

\begin{lem}\label{lem: nqc perturb weak nqc}
    Let $d$ be a positive integer and $\epsilon$ a positive real number. Let $a_1,\dots,a_k$ be positive real numbers that are linearly independent over $\mathbb Q$. Let $\delta_0:=\frac{\epsilon}{2(2d+\epsilon)}$. Then there exists function $\tau: (0,\delta_0]\rightarrow\mathbb R_{>0}$ depending only on $d,\epsilon$, and $a_1,\dots,a_k$ satisfying the following.
    
    Let $(X,\Ff,B,\Mm)/U$ be an lc algebraically integrable generalized foliated quadruple and $\Nn$ an NQC$/U$ $\bb$-divisor on $X$, such that 
    \begin{enumerate}
      \item   $(X,\Ff,B,\Mm+\Nn)$ is lc,
      \item   $K_{\Ff}+B+\Mm_X$ is nef$/U$,
      \item   $K_{\Ff}+B+\Mm_X+\Nn_X$ is $\epsilon$-NQC$/U$, and
      \item $\Nn=\sum a_i\Nn_i$, where each $\Nn_i$ is a nef$/U$ Cartier $\bb$-divisor and each $\Nn_{i,X}$ is Cartier.
    \end{enumerate}
Then $$K_{\Ff}+B+\Mm_X+(1-\delta)\Nn_X$$ is $\tau(\delta)$-nef$/U$ for any $\delta\in (0,\delta_0)$.
\end{lem}
\begin{proof}
Let $M:=\max\left\{\frac{2d}{a_i}\Big| 1\leq i\leq k\right\}$. Consider the set
$$\Ii_0:=\left\{-\sum_{i=1}^ka_i\gamma_i\Bigg| \gamma_i\in\mathbb Z\cap (-\infty,M]\right\}.$$
It is easy to see that $\Ii_0$ is a set whose only accumulation point is $+\infty$. In particular,
$$\gamma_0:=\inf\{\gamma\in\Ii_0\mid \gamma>0\}=\min\{\gamma\in\Ii_0\mid\gamma>0\}>0.$$
In the following, we shall show that
$$\tau: \delta\rightarrow\min\left\{\frac{\epsilon}{2},\delta\gamma_0\right\}$$
satisfies our requirements.

Fix $\delta\in (0,\delta_0)$. By Lemma \ref{lem: weak nqc only check smallest}, we only need to show that $(K_{\Ff}+B+\Mm_X+(1-\delta)\Nn_X)\geq\tau(\delta)$ for any curve $C$ on $X$ satisfying the following:
\begin{itemize}
    \item $(K_{\Ff}+B+\Mm_X+(1-\delta)\Nn_X)\cdot C>0$.
    \item $C$ spans an extremal ray in $\overline{NE}(X/U)$.
    \item For any curve $C'$ on $X$ such that $C'\equiv_U\lambda C$ for some real number $\lambda>0$, we have $\lambda\geq 1$. 
\end{itemize}
For any such curve $C$, there are two possibilities.

\medskip

\noindent\textbf{Case 1}. $(K_{\Ff}+B+\Mm_X+\Nn_X)\cdot C>0$. Since $K_{\Ff}+B+\Mm_X+\Nn_X$ is $\epsilon$-NQC$/U$, $$(K_{\Ff}+B+\Mm_X+\Nn_X)\cdot C\geq\epsilon.$$ 

Suppose that $(K_{\Ff}+B+\Mm_X+(1-\delta)\Nn_X)\cdot C<\frac{\epsilon}{2}$. Then
\begin{align*}
    &(K_{\Ff}+B+\Mm_X)\cdot C\\
    =&\frac{1}{\delta}((K_{\Ff}+B+\Mm_X+(1-\delta)\Nn_X)\cdot C-(1-\delta)(K_{\Ff}+B+\Mm_X+\Nn_X)\cdot C)\\
    <&\frac{1}{\delta}\left(\frac{\epsilon}{2}-(1-\delta)\epsilon\right)=\epsilon-\frac{\epsilon}{2\delta}<\epsilon-\frac{\epsilon}{2\delta_0}=-2d.
\end{align*}
Let $R$ be the extremal ray spanned by $C$. Then $R$ is a $(K_{\Ff}+B+\Mm_X)$-negative extremal ray. By \cite[Theorem 2.3.1]{CHLX23}, $R$ is spanned by a curve $C'$ such that 
$$(K_{\Ff}+B+\Mm_X)\cdot C'\geq-2d.$$
Therefore, $C'\equiv \lambda C$ for some $\lambda\in (0,1)$, which is not possible. Therefore, 
$$(K_{\Ff}+B+\Mm_X+(1-\delta)\Nn_X)\cdot C\geq\frac{\epsilon}{2}\geq\tau(\delta).$$

\medskip

\noindent\textbf{Case 2}. $(K_{\Ff}+B+\Mm_X+\Nn_X)\cdot C=0$. 

Suppose that $\Nn_{i,X}\cdot C>\frac{2d}{a_i}$ for some $i$. Then
    $$(K_X+B+\Mm_X+\Nn_X-a_i\Nn_{i,X})\cdot C<-2d.$$
    Let $R$ be the extremal ray spanned by $C$. Then $R$ is a $(K_{\Ff}+B+\Mm_X+\Nn_X-a_i\Nn_{i,X})$-negative extremal ray. By \cite[Theorem 2.3.1]{CHLX23}, $R$ is spanned by a curve $C'$ such that 
$$(K_{\Ff}+B+\Mm_X+\Nn_X-a_i\Nn_{i,X})\cdot C'\geq-2d.$$
Therefore, $C'\equiv \lambda C$ for some $\lambda\in (0,1)$, a contradiction.

Therefore, $\Nn_{i,X}\cdot C\leq\frac{2d}{a_i}\leq M$ for each $i$. We have
$$0<(K_{\Ff}+B+\Mm_X+(1-\delta)\Nn_X)\cdot C=-\delta\sum a_i(\Nn_{i,X}\cdot C)\in \left\{\sum(-\delta a_i)\gamma_i\Big| \gamma_i\in\mathbb Z,\gamma_i\leq M\right\}.$$
Therefore,
$$(K_{\Ff}+B+\Mm_X+(1-\delta)\Nn_X)\cdot C\geq\delta\gamma_0\geq\tau(\delta).$$
\end{proof}

\begin{lem}\label{lem: single step enef trivial}
    Let $d$ be a positive integer and $\epsilon$ a positive real number. Let $(X,\Ff,B,\Mm)/U$ be an lc algebraically integrable generalized foliated quadruple of dimension $d$ and $D$ an $\epsilon$-nef $\Rr$-divisor on $X$. Then for any real number $l>\frac{2d}{\epsilon}$, any single step of a
    $$(K_{\Ff}+B+\Mm_X+lD)\text{-MMP}/U$$
    is $D$-trivial.
\end{lem}
\begin{proof}
    Let $R$ be a $(K_{\Ff}+B+\Mm_X+lD)$-negative extremal ray$/U$. Then $D$ is a $(K_{\Ff}+B+\Mm_X)$-negative extremal ray$/U$. By \cite[Theorem 2.3.1]{CHLX23}, $D$ is spanned by a curve $C$ such that 
    $$0<-(K_{\Ff}+B+\Mm_X)\cdot C\leq 2d.$$
    Therefore, $lD\cdot C\leq 2d$, so $D\cdot C<\epsilon$. Therefore, $D\cdot C=0$, and the lemma follows.
\end{proof}

\begin{lem}\label{lem: trivial ray when -delta}
    Let $d$ be a positive integer and $\epsilon$ a positive real number. Then $\delta_0:=\frac{\epsilon}{2d+\epsilon}$ satisfies the following. Let $(X,\Ff,B,\Mm)/U$ be an lc algebraically integrable generalized foliated quadruple and $D$ an $\Rr$-divisor on $X$, such that  $K_{\Ff}+B+\Mm_X+D$ is $\epsilon$-NQC$/U$. Then for any $\delta\in (0,\delta_0)$, any $(K_{\Ff}+B+\Mm_X+(1-\delta)D)$-non-positive extremal ray$/U$ is a $(K_{\Ff}+B+\Mm_X+D)$-trivial extremal ray$/U$. 
\end{lem}
\begin{proof}
Fix $\delta\in (0,\delta_0)$ and let $R$ be a $(K_{\Ff}+B+\Mm_X+(1-\delta)D)$-non-positive extremal ray$/U$. If $R$ is not $(K_{\Ff}+B+\Mm_X+D)$-trivial, then $R$ is $(K_{\Ff}+B+\Mm_X+D)$-positive, hence $(K_{\Ff}+B+\Mm_X)$-negative. By \cite[Theorem 2.3.1]{CHLX23}, $R$ is spanned by a curve $C$ such that $0<-(K_{\Ff}+B+\Mm_X)\cdot C\leq 2d$. Since $K_{\Ff}+B+\Mm_X+D$ is $\epsilon$-NQC$/U$, $(K_{\Ff}+B+\Mm_X+D)\cdot C\geq\epsilon$. Thus
\begin{align*}
   0&=(K_{\Ff}+B+\Mm_X+(1-\delta)D)\cdot C\\
   &=(1-\delta)(K_{\Ff}+B+\Mm_X+D)\cdot C+\delta(K_{\Ff}+B+\Mm_X)\cdot C\\
   &\geq (1-\delta)\epsilon-2d\delta>\epsilon-(2d+\epsilon)\delta_0=0,
\end{align*}
which is not possible. The lemma follows.
\end{proof}

\begin{lem}\label{lem: cartier preseve under mmp gfq}
Let $(X,\Ff,B,\Mm)/U$ be a $\Qq$-factorial ACSS algebraically integrable foliated quadruple and $D$ a Cartier divisor on $X$. Let $\phi: X\dashrightarrow X'$ be a birational map that is a step of a $(K_{\Ff}+B+\Mm_X)$-MMP$/U$ such that $\phi$ is $D$-trivial. Then $\phi_*D$ is Cartier.
\end{lem}
\begin{proof}
    Let $X\xrightarrow{f} T\xleftarrow{g} X'$ be this step of the MMP. By \cite[Theorem 16.1.3]{CHLX23}, $\mathcal{O}_X(mD)$ is globally generated over $T$ for any $m\gg 0$. Therefore, $mD=f^*L_m$ and $(m+1)D=f^*L_{m+1}$ for some Cartier divisors $L_m,L_{m+1}$ for any $m\gg 0$, so $D=f^*(L_{m+1}-L_m)$. Therefore, $\phi_*D=g^*(L_{m+1}-L_m)$ is Cartier.
\end{proof}

\begin{prop}\label{prop: kf+k_x mmp trivial}
    Let $(X,\Ff,B,\Mm)/U$ be an lc algebraically integrable generalized foliated quadruple and let $\Nn$ be an NQC$/U$ $\bb$-divisor on $X$, such that
    \begin{itemize}
        \item $(X,\Ff,B,\Mm+\Nn)$ is lc,
        \item $K_{\Ff}+B+\Mm_X$ is nef$/U$, and
        \item  $K_{\Ff}+B+\Mm_X+\Nn_X$ is NQC$/U$.
    \end{itemize}
    Then there exists a real number $\delta_0\in (0,1)$ and a function $\mu: (0,\delta)\rightarrow (0,+\infty)$ satisfying the following. Assume that
    \begin{enumerate}
        \item $\delta\in (0,\delta_0)$ is a real number,
        \item $l>\mu(\delta)$ is a real number, and
        \item $(X,\Ff',B',\Mm')/U$ is an lc algebraically integrable generalized foliated quadruple such that one of the following conditions hold:
        \begin{enumerate}
            \item $(X',\Ff',B',\Mm')$ is $\Qq$-factorial ACSS.
            \item $\Ff'=T_X$ and $\Mm'$ is NQC$/U$.
            \item $\Ff'=T_X$ and $(X',B',\Mm')$ is klt.
        \end{enumerate}
    \end{enumerate}
    Then any sequence of steps of a 
    $$((K_{\Ff'}+B'+\Mm'_X)+l(K_{\Ff}+B+\Mm_X+(1-\delta)\Nn_X))\text{-MMP}/U$$
    is $(K_{\Ff}+B+\Mm_X+(1-\delta)\Nn_X)$-trivial, $(K_{\Ff}+B+\Mm_X+\Nn_X)$-trivial, and $\Nn_X$-trivial.
\end{prop}
\begin{proof}
Let $d:=\dim X$.

\medskip

\noindent\textbf{Step 1}. In this step we introduce real numbers $\epsilon$, $a_1,\dots,a_k$, and $\bb$-divisors $\Nn_1,\dots,\Nn_k$. We construct  $\delta_0$ and $\mu$ so that they only depend on $d,\epsilon,a_1,\dots,a_k$.

 Since $K_{\Ff}+B+\Mm_X+\Nn_X$ is NQC$/U$, there exists a positive real number $\epsilon$ such that  $K_{\Ff}+B+\Mm_X+\Nn_X$ is $\epsilon$-NQC$/U$.

Since $\Nn$ is an NQC$/U$ $\bb$-divisor, by Lemma \ref{lem: nqc good decomposition}, there exist positive real numbers $a_1,\dots,a_k$ that are linearly independent over $\mathbb Q$, such that $\Nn=\sum_{i=1}^ka_i\Nn_i$, where each $\Nn_i$ is nef$/U$ Cartier. Since $K_{\Ff}+B+\Mm_X+\Nn_X$ and $K_{\Ff}+B+\Mm_X$ are $\Rr$-Cartier, $\Nn_X$ is $\Rr$-Cartier. Therefore, $\Nn_{i,X}$ is $\Qq$-Cartier for each $i$. We let $I$ be a positive integer such that $I\Nn_{i,X}$ is Cartier for each $i$. Possibly replacing each $a_i$ with $\frac{a_i}{I}$ and $\Nn_i$ with $I\Nn_i$, in the following, we shall assume that $\Nn_{i,X}$ is Cartier for each $i$.

We let $\delta_0:=\frac{\epsilon}{2(2d+\epsilon)}$ and let $\tau: (0,\delta_0)\rightarrow (0,+\infty)$ be the function constructed in Lemma \ref{lem: nqc perturb weak nqc} which depends only on $d,\epsilon,a_1,\dots,a_k$.  We define $\mu(\delta):=\frac{2d}{\tau(\delta)}$. 

\medskip

\noindent\textbf{Step 2}. In this step we prove the proposition by induction on the number of steps of the MMP.

\begin{claim}\label{claim: kx+kf mmp trivial}
Let $\delta\in (0,\delta_0)$ and $l>\mu(\delta)$ be two real numbers. Let $n$ be a non-negative integer and let 
$$X_0:=X\dashrightarrow X_1\dashrightarrow\dots\dashrightarrow X_n\dashrightarrow X_{n+1}$$
be a sequence of steps of a
$$((K_{\Ff'}+B'+\Mm'_X)+l(K_{\Ff}+B+\Mm_X+(1-\delta)\Nn_X))\text{-MMP}/U.$$
For each $j$, we let $\Ff_j,\Ff_j'$ be the induced foliations of $\Ff,\Ff'$ on $X_j$, and let $B_j,B_j'$ be the images of $B,B'$ on $X_j$ respectively. Then for any $0\leq j\leq n+1$,
\begin{enumerate}
\item $(X_j,\Ff_j,B_j,\Mm+\Nn)$ is lc,
\item $K_{\Ff_j}+B_j+\Mm_{X_j}$ is nef$/U$,
\item $\Nn_{i,X_j}$ is Cartier for each $i$,
\item $K_{\Ff_j}+B_j+\Mm_{X_j}+\Nn_{X_j}$ is $\epsilon$-NQC$/U$,
\item $(X_j,\Ff'_j,B_j',\Mm')$ is lc, and
\begin{enumerate}
            \item if $(X,\Ff',B',\Mm')$ is $\Qq$-factorial ACSS, then $(X_j,\Ff'_j,B_j',\Mm')$ is $\Qq$-factorial ACSS, and
            \item if $\Ff'=T_X$ and $(X',B',\Mm')$ is klt, then $(X_j,B_j',\Mm')$ is klt,
\end{enumerate}
and
\item if $j\leq n$, then $X_{j}\dashrightarrow X_{j+1}$ is 
    $(K_{\Ff_j}+B_j+\Mm_{X_j}+(1-\delta)\Nn_{X_j})$-trivial and $\Nn_{i,X_j}$-trivial for each $i$.
\end{enumerate}
\end{claim}
\begin{proof}
We prove Claim \ref{claim: kx+kf mmp trivial} by induction on $n$. When $n=0$, (1-5) hold by our construction. By induction, we may assume that (1-5) hold for $j\leq n$ and (6) holds for $j\leq n-1$. By Lemma \ref{lem: nqc perturb weak nqc}, $K_{\Ff_n}+B_n+\Mm_{X_n}+(1-\delta)\Nn_{X_n}$ is $\tau(\delta)$-nef$/U$ for any $\delta\in (0,\delta_0)$. By Lemma \ref{lem: single step enef trivial}, $X_n\dashrightarrow X_{n+1}$ is  $(K_{\Ff_n}+B_n+\Mm_{X_n}+(1-\delta)\Nn_{X_n})$-trivial. By Lemma \ref{lem: trivial ray when -delta}, $X_n\dashrightarrow X_{n+1}$ is  $(K_{\Ff_n}+B_n+\Mm_{X_n}+\Nn_{X_n})$-trivial, hence $\Nn_{X_n}$-trivial. Since $a_1,\dots,a_k$ are linearly independent over $\mathbb Q$, $X_n\dashrightarrow X_{n+1}$ is  $\Nn_{i,X_n}$-trivial for any $i$. This deduces (6) for $j=n$.

We left to prove (1-5) for $j=n+1$. Since $X_n\dashrightarrow X_{n+1}$ is $(K_{\Ff_n}+B_n+\Mm_{X_n}+(1-\delta)\Nn_{X_n})$-trivial, we get (1)(4) by induction hypothesis, and we have that $X_n\dashrightarrow X_{n+1}$ is a step of a $(K_{\Ff'}+B'+\Mm_X')$-MMP$/U$. By (6) for $j=n$, $X_n\dashrightarrow X_{n+1}$ is $(K_{\Ff_n}+B_n+\Mm_{X_n})$-trivial, so (2) follows from the induction hypothesis.  (3) follows from the induction hypothesis and Lemma \ref{lem: cartier preseve under mmp gfq}. (5) follows from induction hypothesis and \cite[Lemma 9.1.4]{CHLX23}.
\end{proof}

\noindent\textit{Proof of Proposition \ref{prop: kf+k_x mmp trivial} continued}. It immediately follows from Claim \ref{claim: kx+kf mmp trivial}(6).
\end{proof}

\end{document}